\documentclass[9pt,twoside]{article}
\usepackage[colorlinks]{hyperref}
\usepackage[margin=1in]{geometry}
\usepackage{times,amsmath,stackengine,amsthm,amsfonts,amssymb,mathtools,color,latexsym,amsxtra,layout}
\usepackage{graphicx}
\usepackage{algorithmicx, algpseudocode}
\usepackage{algorithm}
\usepackage{multirow}
\usepackage{comment, centernot}

\newcommand{\honeinner}[1]{\left \langle #1  \right \rangle_{H^1}}

\newcommand{\ltwoinner}[1]{\left \langle #1  \right \rangle_{2}}
\newcommand{\loneinftynorm}[1]{\left \| #1 \right \|_{{1,\infty}}}

\newcommand{\ltwonorm}[1]{\left \| #1 \right \|_{2}}
\newcommand{\lfonorm}[1]{\left \| #1 \right \|_{4}}
\newcommand{\Mnorm}[1]{\left \| #1 \right \|_{M}}

\newcommand{\honenorm}[1]{\left \| #1 \right \|_{H^1}}

\newcommand{\p}{\mathcal{P}}

\theoremstyle{plain}
\newtheorem{theorem}{Theorem}[section]

\newtheorem{o-thm}[theorem]{Theorem}

\theoremstyle{definition}

\newtheorem{remark}[theorem]{Remark}

\author{Sophie Moufawad \thanks{Sophie Moufawad, Department of Mathematics, American University of Beirut, Lebanon (sm101@aub.edu.lb)}\and Nabil Nassif \thanks{Nabil Nassif, Department of Mathematics, American University of Beirut, Lebanon (nn12@aub.edu.lb)}}

\title{Newton Type Methods for solving a Hasegawa-Mima Plasma Model }

\date{\today}

\begin{document}

\maketitle

\begin{abstract}
\noindent In  \cite{kn}, the non-linear space-time Hasegawa-Mima plasma equation is formulated as a coupled system of two linear PDEs, a solution of which is a pair $(u,w)$, with $w=(I-\Delta)u$. The first equation is of hyperbolic type and the second of elliptic type. Variational frames for obtaining weak solutions to the initial value Hasegawa-Mima problem with periodic boundary conditions were also derived. 
In a more recent work \cite{FEHM}, a numerical approach consisting of a finite element space-domain combined with an Euler-implicit time scheme was used to discretize the coupled variational Hasegawa-Mima model.
A semi-linear version of this implicit nonlinear scheme was tested for several types of initial conditions. This semi-linear scheme proved to lack efficiency for long time, which necessitates imposing a cap on the magnitude of the solution.\\
To circumvent this difficulty, in this paper, we use Newton-type methods (Newton, Chord and an introduced Modified Newton method) to solve numerically the fully-implicit non-linear scheme. Testing these methods in FreeFEM++ indicates significant improvements as no cap needs to be imposed for long time. In the sequel, we demonstrate the validity of these methods by proving several results, in particular the convergence of the implemented methods. 




\end{abstract}
\section{Introduction}	
In this paper, we consider the Hasegawa-Mima equation 
 \cite{hm77,hm78}, 
  given by \eqref{hm}
\begin{equation}\label{hm} 
-\Delta u_t+u_t = \{u,\Delta u\} + \{p,u\}
\end{equation}
\noindent where $\{u,v\}=u_xv_y-u_yv_x$ is the Poisson bracket, $u(x,y,t)$ describes the electrostatic potential, $p(x,y)= \ln \dfrac{n_0}{\omega_{ci}}$ is a function depending on the background particle density $n_0$ and the ion cyclotron frequency $\omega_{ci}$, which in turn depends on the initial magnetic field. 

In \cite{kn}, the Hasegawa-Mima model on a square domain with periodic boundary conditions, is reformulated as a hyperbolic-elliptic coupled system of PDEs, 
where a new variable $w=-\Delta u + u$ is introduced, leading to \eqref{HMC2}
\begin{equation}\label{HMC2}
\left\{\begin{array}{lll}
w_t + \vec{V}(u) \cdot \nabla w = \{p,u\}= \vec{V}(p) \cdot \nabla u & \mbox{on }  \Omega \times  (0,T] &(1)\\
-\Delta u+u=w &\mbox{on }  \Omega \times  (0,T] &(2)\\
\mbox{PBC's on } u,\, u_x,\, u_y,\,w & \mbox{on }  \partial \Omega\times [0,T] & (3)\\
u(0)=u_0 \mbox{ and }w(0)=w_0 & \mbox{ on } \overline{\Omega}.  & (4)\\
\end{array}\right.
\end{equation} 
where $\vec{V}(u)= -u_y \vec{\textbf{i}} + u_x \vec{\textbf{j}}$ is a {\it divergence-free vector field} ($\mbox{div}(\vec{V}(u))=0$).

The full discretization of the coupled system were obtained in \cite{FEHM} where starting with the given initial condition at $t=0$, the subsequent solutions are approximated for a chosen time interval $\tau$, to reach the end time $T$ in a finite number of steps. These can be sumarized in the following three sections. 
\subsection{Time Integral Variational Formulation}\label{sec:intro1}
The skew-symmetry property, $ \ltwoinner{\vec{V}(u) \cdot \nabla v,w} = -  \ltwoinner{\vec{V}(u) \cdot \nabla w,v},$ $\forall u \in  H^2 \cap H_P^1,  v\in W^{1,\infty} \cap H^1_P,$ $w \in H^1_P$,  in addition to  the integration of
$$\ltwoinner{w_t,v} = \ltwoinner{\vec{V}(u) \cdot \nabla v,w} + \ltwoinner{ \vec{V}(p) \cdot \nabla u,v},   $$
on the interval $[t,t+\tau]$ leads to seeking the pair  $$\{u,w\}: [0,T] \rightarrow {H^2 \cap H_P^1 \,\times\, {L^2} }$$
 such that $\forall v \in W^{1,\infty} \cap H^1_P,\quad  0\le t \leq s \leq t+\tau\le T,\,$ and $\,\tau>0$
 \begin{equation}\label{HMC2-tau}
\left\{\begin{array}{ll}
\ltwoinner{w(t+\tau)-w(t),v} =\int_t^{t+\tau} \ltwoinner{\vec{V}(u(s)) \cdot \nabla v,w(s)}+ \ltwoinner{ \vec{V}(p) \cdot \nabla u(s),v} ds&(1)\\
\honeinner{u(s),v}=\ltwoinner{w(s),v}, \ &(2)\\
u(0)=u_0\in H^2\cap H^1_P &\\ 
w(0)=w_0=u_0-\Delta u_0&\\
\end{array}\right.\vspace{-2mm}
\end{equation}

\subsection{ Full $\mathbf{\mathbb{P}_1}$  Finite-Element Space, Euler-Implicit Time Discretizations}\label{sec:intro2}
Let 
$\p_x=\{x_i|i=1,...,n\}$ be a partition of $(0,L)$: $0=x_1<x_2<...<x_n=L$ in the $x$ direction and similarly in the $y$ direction, $\p_y=\{y_j|j=1,...,n\}$. Let now:
$$\mathcal{N}=\{P_I(x_i,y_j)|I=1,2,...,N=n^2\}=\p_x\times\p_y,$$
be a structured set of nodes covering $\overline{\Omega}$. Based on $\mathcal{N}$,
one obtains a conforming (Delaunay) structured triangulation $\mathcal{T}$ of $\overline{\Omega}$, i.e., $\mathcal{T}=\{E_J|J=1,2,...,M\},\,\,\overline{\Omega}=\cup_{J}{E_J}$. The $\mathbb{P}_1$ finite element subspace $X_N$ of $H^1(\Omega)$ is given by:
 $$X_N=\{v\in C(\overline{\Omega})|v\mbox{ restricted to }E_J\in\mathbb{P}_1,\,J=1,2..,M\} \subset W_P^{1,p}, \quad  1 \leq p \leq \infty$$ 
  with ${\bigcup}_{N\geq 1}\{X_N\}$  dense in $H^1(\Omega)$. For that purpose, we let ${B}_N=\{\varphi_I|\,I=1,2,...N\}$ be a finite element basis of functions with compact support in $\Omega$, i.e.,: 
 $$\forall v_N \in X_N:\,v_N(x,y)=\sum_{I=1}^N{V_{I}\varphi_I(x,y)},\, V_I=v_N(x_I,y_I).$$
  
 To obtain a fully discrete scheme, we start by projecting  \eqref{HMC2-tau} on  $X_{N,P}\times X_{N,P}$,  seeking therefor the pair  $$\{u,w\}: [0,T] \rightarrow X_{N,P}\times X_{N,P} $$
  such that 
\begin{equation}\label{HMC-Comp-Mod}
\left\{\begin{array}{ll}
\ltwoinner{w_N(t+\tau) - w_N(t),v} = \int_t^{t+\tau}\ltwoinner{\vec{V}(u_N(s)) \cdot \nabla v,w_N(s)} +\ltwoinner{ \vec{V}(p) \cdot \nabla u_N(s),v}ds&(1)\\
\honeinner{u_N(s),v}=\ltwoinner{w_N(s),v},  &(2)\\
w_N(0)=\pi_N(w_0)&\\
\honeinner{u_N(0),v}=\ltwoinner{w_N(0),v}
\end{array}\right.
\end{equation} 
$\forall v \in X_{N,P},\quad  0\le t \leq s \leq t+\tau\le T,\,$ and $\,\tau>0$,
where $\pi_N(v):=\sum_{I=1}^N \ltwoinner{v,\varphi_I}\varphi_I(x,y)\in X_N$ is the L2 projection of  $v$ on $X_N$.
 In addition, to obtain the Euler-Implicit formulation, we replace the term $\int_t^{t+\tau}{\ltwoinner{\vec{V}(u(s)) \cdot \nabla v,w(s)}}$ with 
${\tau}\ltwoinner{\vec{V}(u(t+\tau)) \cdot \nabla v,w(t+\tau)}$, thus yielding the following fully implicit {Computational Model} \eqref{HMC-Comp}.
  
  Given $(u_N(t),w_N(t)) \in X_{N,P}\times X_{N,P}$, one seeks $(u_N(t+\tau),w_N(t+\tau))\in  X_{N,P}\times X_{N,P}$, such that:

\begin{equation}\label{HMC-Comp}
\left\{\begin{array}{ll}
\ltwoinner{w_N(t+\tau)-w_N(t),v} = \tau\ltwoinner{\vec{V}(u_N(t+\tau)) \cdot \nabla v,w_N(t+\tau)} &\\
\qquad \qquad \qquad \qquad \qquad \qquad \qquad +\; {\tau}\ltwoinner{ \vec{V}(p) \cdot \nabla u_N(t+\tau),v}, &(1)\\
\honeinner{u_N(s),v} = \ltwoinner{w_N(s),v},  \quad \forall s\in\{t,t+\tau\},  &(2)\\
\end{array}\right.
\end{equation} 
$\forall v\in  X_{N,P}, $ and $\forall  t\in [0,T]$.

\subsection{ The Non-Linear Algebraic system}\label{sec:intro3}
When implementing system \eqref{HMC-Comp} one takes periodicity into account, reducing the degrees of freedom from $N = n^2$ to $N = (n-1)^2$. Thus, in matrix notations and using the expressions: $$w_N(t)=\sum_{I=1}^N{W_I(t)\varphi_I(x,y)},\, \mbox{ and } u_N(x,y,t)=\sum_{J=1}^N{U_J(t)\varphi_J(x,y)}, $$  where $W_I(t)=w_N(x_I,y_I,t)$, and $U_J(t)=w_N(x_J,y_J,t),$ then (\ref{HMC-Comp}) can be rewritten as follows:\\
Given $(U(t),W(t)) \in \mathbb{R}^{N}\times  \mathbb{R}^{N}$, seek $(U(t+\tau),W(t+\tau))\in  \mathbb{R}^{N}\times  \mathbb{R}^{N}$, such that:
\begin{equation}\label{HMC-Disc-Comp}
\left\{\begin{array}{ll}
(M+ \tau\,S(U(t+\tau))\;W(t+\tau)-\tau\, R\;U(t+\tau)=M\,W(t) &(1)\\
KU(s)=MW(s), \;\;\; \forall s \in \{t,t+\tau\} &(2)\\
\end{array}\right.
\end{equation}
with $M$, $K$, $S(U)$ and $R$, $N\times N$ matrices, whose entries are defined as follows for $1\le I,J\le N$:
\begin{itemize}
\item $M_{I,J}=\ltwoinner{\varphi_I,\varphi_J}$,\quad $M$ is the well-known Mass matrix for periodic boundary conditions. 
\item $K_{I,J}=\honeinner{\varphi_I,\varphi_J}$, \quad $K= M+A$, where $A$ is the stiffness matrix for periodic boundary conditions.
\item { $R_{I,J}= \ltwoinner{\vec{V}(p).\nabla\varphi_{J},\varphi_I}$.}
\item { $S_{I,J}(U)= -\ltwoinner{\vec{V}(u_N) \cdot \nabla\varphi_I,\varphi_J} =  \ltwoinner{\vec{V}(u_N) \cdot \nabla\varphi_J,\varphi_I}$.} 

\end{itemize}
In \cite{FEHM}, we prove the existence of a solution to system \eqref{HMC-Disc-Comp}  for 
$\tau \leq \dfrac{1}{2||p||_{1,\infty}},$
 with a stronger restriction for uniqueness, 
\begin{equation} \label{existTau}
\tau \leq \min\left\{\dfrac{1}{2||p||_{1,\infty}}, \dfrac{h^2}{16{c}_{0,inv}\Mnorm{W(t)}}\right\}
\end{equation}
 where ${c}_{0,inv}$ is an inverse inequality constant, as obtained in Ciarlet \cite{ciarlet} Theorem 3.2.6  $$\forall v\in X_N, \honenorm{v} \leq \dfrac{{c}_{0,inv}}{h} \ltwonorm{v}.$$
 \begin{remark}
\noindent However, note that in our computations we did not use the restrictive condition \eqref{existTau}, as we took $\tau = O(h)$.
  \end{remark}

\noindent The nonlinearity of the problem originates from $S(U)$, that must be computed at each iteration. The derivation of $S(U)$ and $R$ is detailed in \cite{FEHMX}, where they have the same block sparsity patterns as that of $M$ and $K$.%
 
\noindent First results to solve \eqref{HMC-Disc-Comp} where obtained in \cite{FEHM} using a simple semi-linear practical approach: 
\begin{equation}\label{HMC-Disc-Comp-semi}
\left\{\begin{array}{ll}
(M+ \tau\,S(U(t))\;W(t+\tau)=M\,W(t)+\tau\, R\;U(t) &(1)\\
KU(t+\tau)=MW(t+\tau), \;\;\;  &(2)\\
\end{array}\right.
\end{equation}

\noindent However, this approach fails to simulate accurately the wave phenomena that is supposed to remain bounded. As a matter of fact, in the semilinear approach, one has to put a cap on the amplitute of the wave that stops the algorithms once this cap value is reached. 
\subsection{Results}
\noindent To remedy the ill-behaviored semi-linear approach \eqref{HMC-Disc-Comp-semi}, we propose in this paper Newton-type algorithms
that are based on the well-known Newton's method for solving the full discrete system \eqref{HMC-Disc-Comp}.
System \eqref{HMC-Disc-Comp} is equivalent to finding $(U,W) \in  \mathbb{R}^{N}\times \mathbb{R}^{N}$ such that:
\begin{equation}\label{eq:18}\left\{
\begin{array}{ll} F_1(U,W) = (M+\tau S(U)) W-\tau RU -Z = 0 \\
F_2(U,W) = KU-MW=0
\end{array}
\right.
\end{equation}
where  $ U=U(t+\tau)$, $W=W(t+\tau)$, and  $Z=MW(t)$ which is given.\\
In vector form, \eqref{eq:18} is equivalent to  \begin{equation}\label{F}
F(U,W) = \begin{bmatrix}
F_1(U,W)\\
F_2(U,W)
\end{bmatrix} = \begin{bmatrix}
(M+\tau S(U)) W-\tau RU -Z\\
KU-MW
\end{bmatrix} = \begin{bmatrix}
-\tau R& M+\tau S(U) \\ 
K&-M
\end{bmatrix} \begin{bmatrix}
U \\
W
\end{bmatrix} - \begin{bmatrix}
Z \\
0
\end{bmatrix}
=0
\end{equation} 
which can be solved using Newton-type methods that require the computation of the Jacobian Matrix of $F(U,W)$, given by
\begin{equation}\label{JF0}
 J_F(U,W)=\begin{bmatrix} \tau B(W) -\tau R & M+\tau S(U)\\ K & -M \end{bmatrix}, 
 \end{equation} 
as derived in section \ref{sec:Jac}. Note that $K, M,$ and $R$ are fixed matrices that are computed once, whereas $S(U)$ and $B(W)$ have to be computed for each $U$ and $W$.

\noindent At each time step, it is assumed that $(U(t),W(t))$  was already computed/approximated, then $( U(t+\tau),W(t+\tau))$, the solution of \eqref{HMC-Disc-Comp}, is approximated using Full Newton's method by solving system \eqref{eq:3.6} iteratively till convergence up to some given tolerance
\begin{align}
J_F(U^{(k)},W^{(k)}) \begin{bmatrix}U^{(k+1)}-U^{(k)}\\ W^{(k+1)}-W^{(k)}\end{bmatrix}&=-F(U^{(k)},W^{(k)}) \label{eq:3.6} 
\end{align}
which is equivalent to solving 
\begin{equation}\label{eq:3.7}
\begin{bmatrix}
 \tau B(W^{(k)}) -\tau R & M+\tau S(U^{(k)})\\ 
 K  & -M 
\end{bmatrix} 
\begin{bmatrix}U^{(k+1)}\\ W^{(k+1)}
\end{bmatrix}=\begin{bmatrix} \tau S(U^{(k)})W^{(k)}+MW(t)\\0
\end{bmatrix}
\end{equation}
where $(U^{(0)},W^{(0)})=(U(t),W(t))$. \\

\noindent
 This paper is divided as follows.\\\\
  In section \ref{sec:2}, we prove an apriori error estimate for the solution $\{u_N(t), w_N(t)\}$ to \eqref{HMC-Comp}, specifically Theorem \ref{thrm:apriori}. \\
   In section \ref{sec:3}, we derive Newton's method and prove the existence of a unique solution to \eqref{eq:3.6} as a consequence of  Theorem \ref{thrm:bddR}, for $\tau = O(h^{2.5})$. Moreover, we prove the local convergence of Newton's method in Theorem \ref{thrm:ConvNM3}.\vspace{3mm}\\
However, the  Jacobian of  $F(U,W)$, $J_F(U^{(k)},W^{(k)})$, has to be recomputed at every Newton iteration, which is computationally intense. 
Thus, in section \ref{sec:Var} we discuss two variants of Newton's method:\begin{itemize}
\item Chord's method (section \ref{sec:CM}) which differs from Newton's method in the sense that the Jacobian matrix, $J_F(U^{(0)},W^{(0)})$, is fixed throughout all the iterations within one time step. Thus, the existence of a solution is also a corollary of Theorem \ref{thrm:bddR}, for $\tau = O(h^{2.5})$. Hence, we prove the local convergence of Chord's method in Theorem \ref{thrm:ConvC}.
\item Modified Newton's method  (section \ref{sec:MN}) which avoids computing $B(W^{(k)})$, leading to a modified Jacobian matrix, $\tilde{J}_F(U,W)$. Thus, we prove the existence of a unique solution to the system solved at each iteration of the Modified Newton's method (section \ref{sec:exis2}) for $\tau = O(h^2)$, as a corollary of Theorem \eqref{thrm:bdd22}. Moreover, we prove the global convergence of the method in Theorem  \ref{thrm:ConvMN3} (section \ref{sec:ConvMN}).
\end{itemize}
\noindent In section \ref{Sec:test}, numerical testing on the three Newton-type methods are performed where we compare the number of iterations, runtime and behavior of solution with respect to time. \\
\noindent In section \ref{sec:Conclude} we give concluding remarks.


\section{Apriori Error Estimates for Solutions to \eqref{HMC-Comp}}\label{sec:2}
In this section, we prove an apriori error estimate on $w_N$ and $u_N$ solutions to \eqref{HMC-Comp},  stated as follows.
\begin{theorem}\label{thrm:apriori} For $t = m\tau \leq T$, $m\in \mathbb{N}$, and $\tau < \dfrac{1}{6\loneinftynorm{p}}$, every solution $\{u_N(t), w_N(t)\}$ to \eqref{HMC-Comp}  satisfies:
\begin{eqnarray}
\ltwonorm{u_N(t)}  &\leq & e^{3 T \loneinftynorm{p}} \ltwonorm{w_N(0)} \;\leq \; e^{3 T \loneinftynorm{p}} \ltwonorm{w_0} \\
\ltwonorm{w_N(t)}  &\leq & e^{3 T \loneinftynorm{p}} \ltwonorm{w_N(0)} \;\leq \; e^{3 T \loneinftynorm{p}} \ltwonorm{w_0}\label{app1}
\end{eqnarray}
\end{theorem}
\begin{proof}
Let $v = u_N(s)$ in the second equation of \eqref{HMC-Comp}, then \begin{eqnarray}
\honeinner{u_N(s),u_N(s)} &=& \honenorm{u_N(s)}^2 \;=\; \ltwoinner{w_N(s),u_N(s)} \;\leq \;  \ltwonorm{w_N(s)}\ltwonorm{u_N(s)} \;\leq \;  \ltwonorm{w_N(s)}\honenorm{u_N(s)}\nonumber\\
\implies \ltwonorm{u_N(s)} &\leq & \honenorm{u_N(s)} \;\leq \; \ltwonorm{w_N(s)}, \qquad \qquad \forall s \in \{t, t+\tau \} \label{eqw}
\end{eqnarray}
Let $v = w_N(t+\tau)$ in the first equation of \eqref{HMC-Comp}, and assuming $p \in C^{\infty}$, then 
\begin{eqnarray}
\ltwoinner{w_N(t+\tau)-w_N(t),w_N(t+\tau)} &=& \tau\ltwoinner{\vec{V}(u_N(t+\tau)) \cdot \nabla w_N(t+\tau),w_N(t+\tau)} + {\tau}\ltwoinner{ \vec{V}(p) \cdot \nabla u_N(t+\tau),w_N(t+\tau)} \nonumber\\
\implies \ltwonorm{w_N(t+\tau)}^2 &=& \ltwoinner{w_N(t),w_N(t+\tau)} + {\tau}\ltwoinner{ \vec{V}(p) \cdot \nabla u_N(t+\tau),w_N(t+\tau)} \nonumber\\
&\leq & \ltwonorm{w_N(t)}.\ltwonorm{w_N(t+\tau)} + \tau \ltwonorm{ \vec{V}(p) \cdot \nabla u_N(t+\tau)}.\ltwonorm{w_N(t+\tau)} \nonumber\\
\mbox{and therefore }  \ltwonorm{w_N(t+\tau)}  &\leq & \ltwonorm{w_N(t)} + \tau \ltwonorm{ \vec{V}(p) \cdot \nabla u_N(t+\tau)} \label{eqw1}
\end{eqnarray}
Note that 
\begin{eqnarray}
\vec{V}(p) \cdot \nabla u_N(t+\tau) &=& u_{N,y}(t+\tau) p_{x} - u_{N,x}(t+\tau) p_{y} \;\leq \; \loneinftynorm{p} (u_{N,y}(t+\tau)  - u_{N,x}(t+\tau)) \nonumber\\
\implies \ltwonorm{\vec{V}(p) \cdot \nabla u_N } &\leq & \loneinftynorm{p}( \ltwonorm{u_{N,y}(t+\tau) } + \ltwonorm{u_{N,x}(t+\tau) } )\nonumber\\
&\leq &  2 \loneinftynorm{p} \honenorm{u_{N}(t+\tau) } \;\leq \; 2 \loneinftynorm{p} \ltwonorm{w_N(t+\tau)} \qquad \qquad \mbox{ by \eqref{eqw}} \label{eqV} 
\end{eqnarray}
Replacing \eqref{eqV} in \eqref{eqw1}, and $t+\tau$ by $t$, we get for $\tau < \dfrac{1}{2  \loneinftynorm{p}}$
\begin{eqnarray}
 \ltwonorm{w_N(t)}  &\leq & \ltwonorm{w_N(t-\tau)} + 2 \tau \loneinftynorm{p} \ltwonorm{w_N(t)} \nonumber \\
\implies  \ltwonorm{w_N(t)}  &\leq & \dfrac{1}{1- 2 \tau \loneinftynorm{p} } \ltwonorm{w_N(t-\tau)} 
\end{eqnarray}
Note that for each $\alpha = 2+\delta >2$, there exists $\delta_{\alpha} > 0$, such that if $\tau \loneinftynorm{p} < \delta_{\alpha}$, then $$\dfrac{1}{1- 2 \tau \loneinftynorm{p} } < 1+\alpha \tau \loneinftynorm{p}$$
Specifically, 
\begin{eqnarray}
(1- 2 \tau \loneinftynorm{p})(1+\alpha \tau \loneinftynorm{p}) &=& 1 + (\alpha -2) \tau \loneinftynorm{p} - 2\alpha \tau^2 \loneinftynorm{p}^2 > 1\\
\iff 2\alpha \tau^2 \loneinftynorm{p}^2 &<& (\alpha - 2) \tau \loneinftynorm{p}\\
\mbox{and therefore } \tau &<& \dfrac{ \alpha - 2 }{2\alpha\loneinftynorm{p}} = \dfrac{\delta}{2(2+\delta)\loneinftynorm{p}} = \dfrac{\delta_{\alpha}}{\loneinftynorm{p}}
\end{eqnarray}
Let $\delta = 1$, i.e. $\alpha =3$, if $\tau < \dfrac{1}{6\loneinftynorm{p}}$ then 
\begin{eqnarray}
\dfrac{1}{1- 2 \tau \loneinftynorm{p} }  &\leq & 1+3 \tau \loneinftynorm{p} \nonumber\\
\implies  \ltwonorm{w_N(t)}  &\leq & (1+3\tau \loneinftynorm{p}) \ltwonorm{w_N(t-\tau)}\\
 \mbox{and therefore: }  \ltwonorm{w_N(t)}  &\leq &  (1+3 \tau \loneinftynorm{p})^{m} \ltwonorm{w_N(0)} \;=\; e^{m\ln(1+3 \tau \loneinftynorm{p})} \ltwonorm{w_N(0)} \nonumber\\
 &\leq &  e^{3m \tau \loneinftynorm{p}} \ltwonorm{w_N(0)} =  e^{3 t \loneinftynorm{p}} \ltwonorm{w_N(0)}  \; \leq \;   e^{3 T \loneinftynorm{p}} \ltwonorm{w_N(0)}\label{aperr}
\end{eqnarray}
where $t = m\tau \leq T$ for $m\in \mathbb{N}$. Moreover, since $\ltwonorm{w_N(0)} \leq \ltwonorm{w_0}$, then the result is obtained.  
\end{proof}

\section{Newton' Method}\label{sec:3}
In this section, we discuss in details Newton's method (section \ref{sec:FN}), for solving the full discrete system \eqref{HMC-Disc-Comp}, by first deriving the corresponding Jacobian matrix (section \ref{sec:Jac}). Then, in sections \ref{sec:exis} and  \ref{sec:Conv}, we prove respectively the existence of a unique solution to  \eqref{mat1B}, and the convergence of Newton's method.
\subsection{Jacobian Matrix}\label{sec:Jac}
The Jacobian matrix of  $ F(U,W)$ defined in \eqref{F}, is a $2N \times 2N$  block matrix given by:
$$J_F(U,W)=\begin{bmatrix} F_{1,U} & F_{1,W}\\F_{2,U} & F_{2,W} \end{bmatrix}=\begin{bmatrix} \tau (S(U)W)_U -\tau R & M+\tau S(U)\\ K & -M \end{bmatrix}$$
\linebreak
Let $B(W)= (S(U)W)_U$. Then, one way to obtain the matrix $B(W)$ is based on the observation that the matrix $S(U)$ is linear in $U$, i.e. 
\begin{eqnarray}
U&=&\begin{bmatrix} U_{I_1} \\ U_{I_2}\\\vdots\\U_{I_{N}}\end{bmatrix}=\sum\limits_{j=1}^{N} U_{I_j}e_j \qquad  \implies \qquad S(U)=S\left(\sum\limits_{j=1}^{N} U_{I_j}e_j\right)=\sum\limits_{j=1}^{N} U_{I_j}S(e_j)\nonumber\\
\mbox{and therefore } B(W)&=& (S(U)W)_U=\begin{bmatrix} S(e_1)W & S(e_2)W &\cdots & S(e_{N})W \end{bmatrix}.\label{BWdef}
\end{eqnarray}
Hence, 
\begin{equation}\label{JF}
 J_F(U,W)=\begin{bmatrix} \tau B(W) -\tau R & M+\tau S(U)\\ K & -M \end{bmatrix}, 
 \end{equation} where $K, M,$ and $R$ are fixed matrices that are computed once.\vspace{2mm}\\ However, $S(U)$ and $B(W)$ have to be computed for each $U$ and $W$. $S(e_1), S(e_2), \cdots, S(e_N)$ can be computed once and stored. Yet for large $N$ values storing this set of $N$ matrices might not be feasible. In this case, they can be recomputed once needed.  

Note that $S(U)$ is a block tridiagonal matrix with 2 additional blocks in the upper right  and lower left corner. Moreover, it is a skew-symmetric matrix ($S(U)^T = -S(U)$) that is linear in $U$,  with  6 nonzero entries per row, 6 nonzero entries per column, and zeros on the diagonal assuming the meshing of $\Omega$ is uniform. 
$$ S(U) =\dfrac{1}{6} \left[ \begin{array}{cccccc} S_{1,1} & S_{1,2}& 0&\cdots&0&S_{1,k}\\
S_{2,1} & S_{2,2}& S_{2,3}&0&\cdots&0\\
0&\ddots&\ddots&\ddots&\ddots&\vdots\\
\vdots&\ddots&S_{j,l}&S_{j,j}&S_{j,i}&0\\
0&\cdots&0&S_{i,j}&S_{i,i}&S_{i,k}\\
S_{k,1}&0&\cdots&0&S_{k,i}&S_{k,k}\\
\end{array}\right]  \qquad \qquad with \quad S_{i,j} \equiv S_{i,j}(U)
$$
\noindent where $i = n-2, j = n-3, k= n-1, l=n-4$, and the $3(n-1)$ nonzero block matrices $S_{i,j}$ are of size $(n-1)\times (n-1)$ with $2(n-1)$ nonzero entries each, and the following sparsity patterns:
\begin{itemize}
\item $S_{i,i}$ for $i = 1,.., n-1$ are tridiagonal matrices with zero diagonal entries, and nonzero $S_{i,i}(1,n-1),$ and $ S_{i,i}(n-1,1)$. \vspace{-1mm}
\item $S_{1,n-1}$ and $S_{i+1,i}$ for $i = 1,2,3,.., n-2$ are lower bidiagonal matrices, with nonzero entry in first row and column $n-1$.\vspace{-1mm} 
\item  $S_{n-1,1}$ and $S_{i,i+1}$ for $i = 1,2,..,n-2$ are upper bidiagonal matrices with nonzero entry in first column and row $n-1$.
\end{itemize}
\noindent 
 As for the explicit expressions/values of the entries, refer to appendix A.2 of \cite{FEHMX}.

\noindent Moreover the matrices $B(W)$ and $S(U)$ satisfy the following relation:
\begin{equation}\label{BS}
B(W)U = S(U)W
\end{equation}
since by the linearity of the $S$ matrix we get
\begin{eqnarray}
\label{identity}
B(W)U =\begin{bmatrix} S(e_1)W  &\cdots & S(e_{N})W \end{bmatrix}
\begin{bmatrix} U_{I_1} \\ \vdots\\U_{I_{N}}
\end{bmatrix}=\sum_{j=1}^{N} U_{I_j}S(e_j)W = S\left(\sum_{j=1}^{N} U_{I_j}e_j\right)W=S(U)W\nonumber
\end{eqnarray}

\subsection{Newton's Method}\label{sec:FN}
 At each time step, it is assumed that $(U(t),W(t))$  was already computed/approximated, then $( U(t+\tau),W(t+\tau))$, the solution of \eqref{HMC-Disc-Comp}, is approximated using Newton's method by solving system \eqref{eq:3.62} iteratively till convergence up to some given tolerance 
\begin{align}
J_F(U^{(k)},W^{(k)}) \begin{bmatrix}U^{(k+1)}-U^{(k)}\\ W^{(k+1)}-W^{(k)}\end{bmatrix}&=-F(U^{(k)},W^{(k)}) \label{eq:3.62} 
\end{align}
where $(U^{(0)},W^{(0)})=(U(t),W(t))$, and $J_F(U,W)$   is the Jacobian of  $ F(U,W)$ that has to be recomputed at every Newton iteration. It is possible to solve system \eqref{eq:3.62} directly and obtain the vector $$\begin{bmatrix}
\Delta U\\ \Delta W
\end{bmatrix} = \begin{bmatrix}
U^{(k+1)}-U^{(k)}\\ W^{(k+1)}-W^{(k)}\\
\end{bmatrix}
$$
then $U^{(k+1)}=U^{(k)}+\Delta U$ and $W^{(k+1)}=W^{(k)}+\Delta W$. However, computing $F(U^{(k)},W^{(k)})$ requires 4 matrix-vector multiplications, $MW^{(k)}, S(U^{(k)})W^{(k)}, RU^{(k)}, $ and $
 KU^{(k)}$. 
But by replacing $F(U^{(k)},W^{(k)})$ and $J_F(U^{(k)},W^{(k)})$ by their expressions, \eqref{F} and \eqref{JF}   respectively, and using property \eqref{BS}, system \eqref{eq:3.62} is reduced to the linear system \eqref{New} where the right-hand side vector requires the computation of just one matrix-vector multiplication $S(U^{(k)})W^{(k)}$. 
\begin{eqnarray}
\begin{bmatrix} 
\tau B(W^{(k)}) -\tau R & M+\tau S(U^{(k)})\\ 
K & -M 
\end{bmatrix}
\begin{bmatrix}
 U^{(k+1)}\\
 W^{(k+1)}
 \end{bmatrix}&=& 
 \begin{bmatrix} 
 \tau B(W^{(k)}) -\tau R &M+\tau S(U^{(k)})\\
 K & -M \end{bmatrix}
 \begin{bmatrix}
 U^{(k)}\\
 W^{(k)}
 \end{bmatrix}-F(U^{(k)},W^{(k)})\nonumber\\
\iff \begin{bmatrix} 
 \tau B(W^{(k)}) -\tau R & M+\tau S(U^{(k)})\\
 K & -M 
 \end{bmatrix}
 \begin{bmatrix}
 U^{(k+1)}\\ 
 W^{(k+1)}
 \end{bmatrix}&=& 
 \begin{bmatrix} 
 \tau B(W^{(k)})U^{(k)} -\tau RU^{(k)}+ MW^{(k)}+\tau S(U^{(k)})W^{(k)}\\
 KU^{(k)}  -MW^{(k)} \end{bmatrix}\nonumber\\ 
& & -\begin{bmatrix}
 (M+\tau S(U^{(k)})) W^{(k)}-\tau RU^{(k)}- Z \\
 KU^{(k)}-MW^{(k)} 
\end{bmatrix} \nonumber\\
\begin{bmatrix}
 \tau B(W^{(k)}) -\tau R & M+\tau S(U^{(k)})\\ 
 K  & -M 
\end{bmatrix} 
\begin{bmatrix}U^{(k+1)}\\ W^{(k+1)}
\end{bmatrix}&=&\begin{bmatrix} \tau S(U^{(k)})W^{(k)}+Z\\0
\end{bmatrix} \label{New}
\end{eqnarray}

Thus at the $(k+1)^{th}$ Newton iteration, the Jacobian matrix $J_F(U^{(k)},W^{(k)})$ has to be recomputed by computing $B(W^{(k)})$ and $S(U^{(k)})$. Similarly, the right-hand side vector is computed. Then, system \eqref{New} is solved. At each timestep, the above procedure is repeated  until convergence, i.e. the relative error $\dfrac{||U^{(k+1)} - U^{(k)}||_2}{||U^{(k)}||_2}$  is less than some given tolerence. This procedure is summarized in Algorithm \eqref{alg:HMC-Newton0} where the vectors $U(t), W(t), U^{(k)}(t), W^{(k)}(t)$ are denoted by $U_t, W_t, U_{t,k} , W_{t,k}$ respectively.

At every iteration of Newton's method, there is a need to solve some linear system of form \eqref{mat1B}, where $[\alpha, \beta]^T\in \mathbb{R}^{2N}$. System  \eqref{mat1B} is equivalent to the linear system \eqref{sys1B} using property \eqref{BS}.

 \begin{eqnarray}
 &&J_F(U,W) \begin{bmatrix}
 \alpha\\
 \beta
 \end{bmatrix} 
 = \begin{bmatrix} \tau B(W) -\tau R & M+\tau S(U)\\ K & -M \end{bmatrix} \begin{bmatrix}
 \alpha\\
 \beta
 \end{bmatrix}
 = \begin{bmatrix}
\tau S(U)W+Z\\
0
 \end{bmatrix} \label{mat1B}\\
  &&\iff 
 \begin{cases}
 \tau S(\alpha)W -\tau R \alpha +  M\beta +\tau S(U)\beta =\tau S(U)W+Z&\\
  K\alpha =  M\beta &
 \end{cases}\label{sys1B}
 \end{eqnarray}
Specifically, at iteration $k+1$ of Newton's method $\alpha = U^{(k+1)}, \beta =  W^{(k+1)}, W = W^{(k)},$  $U = U^{(k)}$  and \\$Z = MW^{(0)} = MW(t)$ 
based on \eqref{New}. 

\noindent To prove the convergence of Newton's methods, we prove first the existence of a unique solution of linear system   \eqref{mat1B} in section \ref{sec:exis}, then we 
conclude by the convergence proof in section \ref{sec:Conv}. 
\subsection{Existence of a Unique Solution to \eqref{mat1B}}\label{sec:exis}
Let $M\gamma = \tau S(U)W+Z$. Then,
to prove the existence of a unique solution to linear system \eqref{mat1B}, we show that the Jacobian matrix in invertible. We start by showing that there exists some $C \in \mathbb{R}$ independent of $\tau$ and $h$, such that $$\Mnorm{\alpha}^2 + \Mnorm{\beta}^2 \leq c\Mnorm{\gamma}^2$$  where $\Mnorm{\alpha}^2 =\alpha^TM\alpha$ and $M$ is the Mass matrix. For that purpose, we use variational formulation.

\noindent Let $\phi_N(x,y)  = \sum\limits_{I=1}^N \alpha_I \varphi_I(x,y)$, $\psi_N(x,y)  = \sum\limits_{I=1}^N \beta_I \varphi_I(x,y)$, and $\xi_N(x,y)  = \sum\limits_{I=1}^N \gamma_I \varphi_I(x,y)$, 
then system \eqref{sys1B} can be expressed in variational form elementwise (for $1\leq I\leq N$) as  \eqref{sys1elemB} 
\begin{equation}\label{sys1elemB}
\begin{cases}
- \tau \ltwoinner{\vec{V}(\phi_N) \cdot \nabla w_N,\varphi_I}  + \tau \ltwoinner{\vec{V}(p).\nabla\phi_N, \varphi_{I}} + \ltwoinner{\psi_N,\varphi_I} - \tau\ltwoinner{\vec{V}(u_N) \cdot \nabla\psi_N,\varphi_I} = \ltwoinner{\xi_N,\varphi_I} &\\
\honeinner{\phi_N,\varphi_I} = \ltwoinner{\psi_N,\varphi_I}&\\
\end{cases}
\end{equation}
based on \eqref{elem1}-\eqref{elem5}.
\begin{eqnarray}
(M\beta)_I &=& \sum\limits_{J=1}^N M_{I,J}\beta_J \;=\; \sum\limits_{J=1}^N \ltwoinner{\varphi_I,\varphi_J}\beta_J \;=\; \ltwoinner{\varphi_I,\sum_{J=1}^N \varphi_J \beta_J} \;=\; \ltwoinner{\varphi_I,\psi} \;=\; \ltwoinner{\psi,\varphi_I} \label{elem1}\\
(K\alpha)_I &=& \sum\limits_{J=1}^N K_{I,J}\alpha_J \;=\; \sum\limits_{J=1}^N \honeinner{\varphi_I,\varphi_J}\alpha_J \;=\; \honeinner{\varphi_I,\sum_{J=1}^N \varphi_J \alpha_J} \;=\;  \honeinner{\varphi_I,\phi} \;=\;  \honeinner{\phi,\varphi_I}\label{elem2}\\
(R\alpha)_I &=&\sum\limits_{J=1}^N R_{I,J}\alpha_J \;=\; \sum\limits_{J=1}^N \ltwoinner{\vec{V}(p).\nabla\varphi_{I},\varphi_J} \alpha_J \;=\;   \ltwoinner{\vec{V}(p).\nabla\varphi_{I},\sum_{J=1}^N \varphi_J \alpha_J } \;=\; \ltwoinner{\vec{V}(p).\nabla\varphi_{I},\phi} \nonumber \\
&=& -\ltwoinner{\vec{V}(p).\nabla\phi, \varphi_{I}} \mbox{\qquad \qquad by skew symmetry} \label{elem3}\\
(S(U)\beta)_I  &=&\sum\limits_{J=1}^N S_{I,J}(U)\beta_J \;=\; \sum\limits_{J=1}^N  \ltwoinner{\vec{V}(u_N) \cdot \nabla\varphi_I,\varphi_J}\beta_J \;=\;  \ltwoinner{\vec{V}(u_N) \cdot \nabla\varphi_I,\sum_{J=1}^N \varphi_J\beta_J} \nonumber\\
&=&  \ltwoinner{\vec{V}(u_N) \cdot \nabla\varphi_I,\psi} \;=\; - \ltwoinner{\vec{V}(u_N) \cdot \nabla\psi,\varphi_I}  \label{elem4}\\
(S(\alpha)W)_I &=&  \ltwoinner{\vec{V}(\phi) \cdot \nabla\varphi_I,w_N} \;=\; - \ltwoinner{\vec{V}(\phi) \cdot \nabla w_N,\varphi_I}  \label{elem5}
\end{eqnarray}

\noindent Moreover, $\Mnorm{\beta}^2 = \ltwonorm{\psi_N}^2$ by \eqref{elem1}. Similarly $\Mnorm{\alpha}^2 = \ltwonorm{\phi_N}^2$ and $\Mnorm{\gamma}^2 = \ltwonorm{\xi_N}^2$. Thus, we need to show that 
\begin{equation}\label{eqTP}
 \ltwonorm{\phi_N}^2 +\ltwonorm{\psi_N}^2 \leq  C\ltwonorm{\xi_N}^2
 \end{equation}
For any $v=\sum\limits_{I=1}^N v_I \varphi_I(x,y) \in X_{N,p}$, system \eqref{sys1elemB} can be written as 
\begin{equation}\label{sys1vB}
\begin{cases}
 \ltwoinner{ \psi_N  - \tau\vec{V}(\phi_N) \cdot \nabla w_N    - \tau \vec{V}(u_N) \cdot \nabla\psi_N + \tau \vec{V}(p).\nabla\phi_N,v}   = \ltwoinner{\xi_N,v} &\\
\honeinner{\phi_N,v} = \ltwoinner{\psi_N,v} &\\
\end{cases}
\end{equation}
\begin{theorem}\label{thrm:bddR}
Let $D := D(\Omega,p,T,w_0) =  {c}_{inv}^2 e^{3 T \loneinftynorm{p}} \ltwonorm{w_0} +2\loneinftynorm{p}$ where ${c}_{inv}$ is an inverse inequality constant as provided in Ciarlet (\cite{ciarlet}, Theorem 3.2.6): 
\begin{equation}\label{ineq:ciar}
\forall v\in X_N, \qquad |v|_{1,4} 
\leq  {c}_{inv} h^{-5/4} \ltwonorm{v}
\end{equation}
Then, for  $h< 1$ and 
 $\tau \leq \min\left\{\dfrac{1}{6||p||_{1,\infty}}, \dfrac{h^2}{16{c}_{0,inv}\Mnorm{W(t)}}, \dfrac{h^{5/2}}{2D}\right\} = O(h^{2.5})$, 
\begin{equation}
\ltwonorm{\phi_N}^2 +\ltwonorm{\psi_N}^2 \leq  8\ltwonorm{\xi_N}^2.
\end{equation}

\end{theorem}
\begin{proof} 
By setting $v = \phi_N$ in the second equation of system \eqref{sys1vB} and using Cauchy-Schwarz we get 
\eqref{eq23B}. 
\begin{eqnarray}
\ltwonorm{\phi_N} ^2 &\leq & \honenorm{\phi_N}^2 \;=\;\ltwoinner{\psi_N,\phi_N}  \;\leq \; \ltwonorm{\psi_N}  \ltwonorm{\phi_N}  \nonumber\\
\implies \ltwonorm{\phi_N} &\leq & \ltwonorm{\psi_N}   \label{eq23B}\\
\therefore \ltwonorm{\phi_N}^2 +\ltwonorm{\psi_N}^2 &\leq & 2\ltwonorm{\psi_N}^2 \qquad \qquad \label{eq24B}
\end{eqnarray}
Thus, to obtain our result, we seek an upper bound on $\ltwonorm{\psi_N}^2$ in terms of $\ltwonorm{\xi_N}^2$. \\
Let $v = \psi_N$ in the first equation of system \eqref{sys1vB} we get \eqref{eqq1B}. Then, using Cauchy-Schwarz we get \eqref{eqq2B}.
\begin{eqnarray}
  \ltwoinner{\xi_N,\psi_N} &=&  \ltwoinner{ \psi_N  - \tau\vec{V}(\phi_N) \cdot \nabla w_N    - \tau \vec{V}(u_N) \cdot \nabla\psi_N + \tau \vec{V}(p).\nabla\phi_N,\psi_N}   \label{eqq1B}\\
\ltwonorm{\psi_N}^2 &=& \ltwoinner{\xi_N,\psi_N} + \ltwoinner{  \tau\vec{V}(\phi_N) \cdot \nabla w_N,\psi_N} + \ltwoinner{  \tau \vec{V}(u_N) \cdot \nabla\psi_N,\psi_N} - \ltwoinner{  \tau \vec{V}(p).\nabla\phi_N,\psi_N} \nonumber\\
  &=& \ltwoinner{\xi_N,\psi_N} + \ltwoinner{  \tau\vec{V}(\phi_N) \cdot \nabla w_N,\psi_N} - \ltwoinner{  \tau \vec{V}(p).\nabla\phi_N,\psi_N}\qquad \qquad \mbox{using skew-symmetry} \nonumber\\
&\leq & \ltwonorm{\xi_N}  \ltwonorm{\psi_N} + \tau \ltwonorm{\vec{V}(\phi_N) \cdot \nabla w_N} \ltwonorm{\psi_N} + \tau \ltwonorm{\vec{V}(p).\nabla\phi_N} \ltwonorm{\psi_N}  \nonumber\\
\therefore  \ltwonorm{\psi_N} &\leq & \ltwonorm{\xi_N}  + \tau \ltwonorm{\vec{V}(\phi_N) \cdot \nabla w_N} + \tau \ltwonorm{\vec{V}(p).\nabla\phi_N} \label{eqq2B}
\end{eqnarray}
We upper bound the last two terms of \eqref{eqq2B} in terms of $ \ltwonorm{\psi_N}$. 
\begin{eqnarray}
\vec{V}(\phi_N) \cdot \nabla w_N &=&- \phi_{N,y} w_{N,x} + \phi_{N,x} w_{N,y} = - \vec{V}(w_N) \cdot \nabla \phi_N  \nonumber\\
 \ltwonorm{\vec{V}(\phi_N) \cdot \nabla w_N } &\leq & \ltwonorm{\phi_{N,y} w_{N,x}} + \ltwonorm{\phi_{N,x} w_{N,y}} \nonumber\\
&\leq & \lfonorm{\phi_{N,y}} \lfonorm{w_{N,x}} + \lfonorm{\phi_{N,x}} \lfonorm{w_{N,y}} \qquad \quad \mbox{By Holder's inequality}\nonumber\\
&\leq &  |w_N|_{1,4} \;.\; |\phi_N|_{1,4} \nonumber\\
&\leq & (c_{inv} h^{-5/4} \ltwonorm{w_N}) \;.\; (c_{inv} h^{-5/4} \ltwonorm{\phi_N}) \qquad \quad \mbox{By}  \eqref{ineq:ciar}\qquad \label{eqCia}\\
&\leq & c_{inv}^2 h^{-5/2} \ltwonorm{w_N}  \;.\; \ltwonorm{\psi_N} \qquad \qquad \qquad \qquad  \mbox{By \eqref{eq23B}} \label{eq3A}
\end{eqnarray}
Assuming $p\in C^{\infty}$, then
\begin{eqnarray}
\vec{V}(p) \cdot \nabla \phi_N &=& \phi_{N,y} p_{x} - \phi_{N,x} p_{y} \;\leq \; \loneinftynorm{p} (\phi_{N,y}  - \phi_{N,x}) \nonumber\\
\implies \ltwonorm{\vec{V}(p) \cdot \nabla \phi_N } &\leq & \loneinftynorm{p}( \ltwonorm{\phi_{N,y} } + \ltwonorm{\phi_{N,x} } )
\;\leq\;   2 \loneinftynorm{p} \honenorm{\phi_{N} } \label{eqP}\\
& \leq &   2 \loneinftynorm{p} \ltwonorm{\psi_N} \qquad \qquad \qquad \qquad  \mbox{By \eqref{eq23B}}\label{eq3B}
\end{eqnarray}
Replacing \eqref{eq3A} and \eqref{eq3B} in \eqref{eqq2B} we get 
\begin{eqnarray}
 \ltwonorm{\psi_N} &\leq & \ltwonorm{\xi_N}  + \tau c_{inv}^2 h^{-5/2} \ltwonorm{w_N} \ltwonorm{\psi_N} + 2\tau   \loneinftynorm{p} \ltwonorm{\psi_N} \nonumber\\
 \implies \ltwonorm{\psi_N} &\leq & \dfrac{1}{\tilde{c}}  \ltwonorm{\xi_N} \nonumber\\
 \mbox{and therefore }  \ltwonorm{\phi_N}^2 +\ltwonorm{\psi_N}^2 &\leq & 2\ltwonorm{\psi_N}^2 \;\leq \;  \dfrac{2}{\tilde{c}^2}  \ltwonorm{\xi_N}^2 
\end{eqnarray}
where  $c =  \dfrac{2}{\tilde{c}^2}$ in \eqref{eqTP} and $\tilde{c} = 1-  \tau c_{inv}^2 h^{-5/2} \ltwonorm{w_N} - 2\tau  \loneinftynorm{p}$.\\
Note that, assuming $h<1$ and using the apriori error estimates \eqref{app1} on $\ltwonorm{w_N}$ for $\tau < \dfrac{1}{6\loneinftynorm{p}}$ we get
\begin{eqnarray}
\tilde{c} = \tilde{c}(\tau, h)  &=& 1-  \tau c_{inv}^2 h^{-5/2} \ltwonorm{w_N(t)} - 2\tau  \loneinftynorm{p} \;=\; 1 -  \tau h^{-5/2}(c_{inv}^2 \ltwonorm{w_N(t)} +2h^{5/2}\loneinftynorm{p} ) \nonumber\\
&\geq & 1 -  \tau h^{-5/2}\left(c_{inv}^2 e^{3 T \loneinftynorm{p}} \ltwonorm{w_0} +2h_0^{5/2}\loneinftynorm{p} \right) \nonumber \\
 &\geq & 1 -  \tau h^{-5/2}\left( c_{inv}^2 e^{3 T \loneinftynorm{p}} \ltwonorm{w_0} +2\loneinftynorm{p} \right)  = 1-\tau h^{-5/2} D \label{taucond} 
\end{eqnarray}
where $D = c_{inv}^2 e^{3 T \loneinftynorm{p}} \ltwonorm{w_0} +2\loneinftynorm{p}$ is a constant independent of $\tau$ and $h$.

\noindent Thus, if $\tau h^{-5/2}\leq \dfrac{1}{2D}$, and $h<1$, then   $\tilde{c} \geq  \dfrac{1}{2}$ for 
 $\tau \leq \dfrac{h^{5/2}}{2D}$ , and  therefore
$
c \;=\;  \dfrac{2}{\tilde{c}^2} \;\leq\;8
$.
\end{proof}
\noindent A consequence of Theorem \eqref{thrm:bddR} is the existence of a unique solution to the linear system \eqref{mat1B}.
\begin{theorem}\label{thrm:NewtonConv}
Let $D := D(\Omega,p,T,w_0) = c_{inv}^2 e^{3 T \loneinftynorm{p}} \ltwonorm{w_0} +2\loneinftynorm{p}$
then system \eqref{mat1B} has a unique solution for  $h<1$ and 
 $\tau \leq \min\left\{\dfrac{1}{6||p||_{1,\infty}}, \dfrac{h^2}{16{c}_{0,inv}\Mnorm{W(t)}}, \dfrac{h^{5/2}}{2D}\right\} = O(h^{2.5})$.
\end{theorem}
\begin{proof}
Let $\gamma = 0$, then $\xi_N = 0$ and by theorem \ref{thrm:bddR} $$\ltwonorm{\phi_N}^2 +\ltwonorm{\psi_N}^2 \leq  0$$ 
for  $h<1$ and 
 $\tau \leq \dfrac{h^{5/2}}{2D}$.

\noindent Thus, $\ltwonorm{\phi_N}^2 = \Mnorm{\alpha}^2 = 0$ and   $\ltwonorm{\psi_N}^2 = \Mnorm{\beta}^2 = 0$, implying that $\alpha = \beta = 0$. Thus, $Null\{J_F(U,W)\} = \{0\}$ implying that $J_F(U,W)$ is invertible and the linear system \eqref{mat1B} has a unique solution.
\end{proof}
\subsection{Convergence}\label{sec:Conv}
We are approximating the solution of the nonlinear system \eqref{HMC-Disc-Comp} by using Newton's method \eqref{New}. 
Let

\begin{minipage}{0.4\textwidth}
\begin{equation}\label{def1}
\begin{cases}
\alpha = U(t+\tau)&\\
\alpha^{(k)} = U^{(k)}(t)& \\
\alpha^{(0)} = U^{(0)}(t) = U(t)&\\
\phi_N(x,y)  = \sum\limits_{I=1}^N \alpha_I \varphi_I(x,y) &\\
\phi_N^{(k)}(x,y) = \sum\limits_{I=1}^N \alpha_I^{(k)} \varphi_I(x,y)& 
\end{cases}
\end{equation}
\end{minipage}
\begin{minipage}{0.15\textwidth}
\qquad 
\end{minipage}
\begin{minipage}{0.4\textwidth}
\begin{equation}\label{def2}
\begin{cases}
\beta = W(t+\tau)&\\
\beta^{(k)} = W^{(k)}(t)& \\
\beta^{(0)} = W^{(0)}(t) = W(t) = \gamma&\\
\psi_N(x,y)  = \sum\limits_{I=1}^N \beta_I \varphi_I(x,y)&\\
\psi_N^{(k)}(x,y) = \sum\limits_{I=1}^N \beta_I^{(k)} \varphi_I(x,y)& 
\end{cases}
\end{equation}
\end{minipage}

 
\noindent Then system \eqref{HMC-Disc-Comp} can be expressed in variational form for any $v=\sum\limits_{I=1}^N v_I \varphi_I(x,y) \in X_{N,p}$ using \eqref{elem1}-\eqref{elem5} as  
\begin{equation}\label{sys6V}
\begin{cases}
 \ltwoinner{ \psi_N   - \tau \vec{V}(\phi_N) \cdot \nabla\psi_N + \tau \vec{V}(p).\nabla\phi_N,v}   = \ltwoinner{\psi^{(0)}_N,v} &\\
\honeinner{\phi_N,v} = \ltwoinner{\psi_N,v} &
\end{cases}
\end{equation}

\noindent Similarly, the iterative Newton's method, \eqref{New} or equivalently \eqref{sys1B}, can be expressed in variational form as 
\begin{equation}\label{sysNV}
\begin{cases}
 \ltwoinner{ \psi_N^{(k+1)}  - \tau\vec{V}(\phi_N^{(k+1)}) \cdot \nabla \psi_N^{(k)}    - \tau \vec{V}(\phi_N^{(k)}) \cdot \nabla\psi_N^{(k+1)} + \tau \vec{V}(p).\nabla\phi_N^{(k+1)},v}   = \ltwoinner{\psi^{(0)}_N,v} - \tau  \ltwoinner{\vec{V}(\phi_N^{(k)}) \cdot \nabla \psi_N^{(k)}, v} \\
\honeinner{\phi_N^{(i)},v} = \ltwoinner{\psi_N^{(i)},v} \qquad  \qquad \mbox{for } i=\{k,k+1\} 
\end{cases}
\end{equation}
Let $e^{(k)} = \psi_N - \psi_N^{(k)}$ and $g^{(k)} = \phi_N - \phi_N^{(k)}$, then to prove the convergence of Newton's method to the unique solution of \eqref{HMC-Disc-Comp}, we prove that there exists some constant $c<1$ such 
\begin{equation}
\ltwonorm{e^{(k+1)}}^2 + \ltwonorm{g^{(k+1)}}^2 \leq \ltwonorm{e^{(0)}}^2\; c^{2(k+1)}  
\end{equation}
\begin{theorem}\label{thrm:convNM2}  Assume that $\phi_N^{(0)}$ is chosen such that $\forall k \geq k_0 \geq 0 , \ltwonorm{\phi_N^{(k+1)} - \phi_N^{(k)}} < \epsilon_{tol}$.
\\ Then, for 
 $\tau \leq \min\left\{\dfrac{1}{6||p||_{1,\infty}}, \dfrac{h^2}{16{c}_{0,inv}\Mnorm{W(t)}}, \dfrac{h^{5/2}}{2D_1 }\right\} = O(h^{2.5})$ there exists a constant $c<1$ such that $$\ltwonorm{e^{(k+1)}}^2 \; \leq \; c^2 \ltwonorm{e^{(k)}}^2$$
 where $D_1 := D_1(\Omega,T,p,w_0) =c_{inv}^2 (\epsilon_{tol}+ e^{3 T \loneinftynorm{p}} \ltwonorm{w_0}) + 2\loneinftynorm{p}$ and $h<1$.
\end{theorem}
\begin{proof}
By theorem \ref{thrm:NewtonConv}, \eqref{sysNV} has a unique solution $\{\phi_N^{(k+1)}, \psi_N^{(k+1)}\}$ for 
 $\tau \leq \min\left\{\dfrac{1}{6||p||_{1,\infty}}, \dfrac{h^2}{16{c}_{0,inv}\Mnorm{W(t)}}, \dfrac{h^{5/2}}{2D}\right\} $ and $h<1$, where  $D := D(\Omega,p,T,w_0) = c_{inv}^2 e^{3 T \loneinftynorm{p}} \ltwonorm{w_0} +2\loneinftynorm{p}$.\\
Then, by subtracting the second equation of \eqref{sysNV} from that of \eqref{sys6V}, we get \eqref{eq2Conv1} for $i = \{k,k+1\}$. Letting $v = g^{(i)}$ we get  \eqref{eq2Conv2}
\begin{eqnarray}
\honeinner{g^{(i)},v} &=& \ltwoinner{e^{(i)},v} \label{eq2Conv1}\\
  \honenorm{g^{(i)}}^2 &=& \ltwoinner{e^{(i)},g^{(i)}} \;\leq \;  \ltwonorm{e^{(i)}} \ltwonorm{g^{(i)}} \;\leq \;  \ltwonorm{e^{(i)}} \honenorm{g^{(i)}} \nonumber\\
\therefore \ltwonorm{g^{(i)}} &\leq & \honenorm{g^{(i)}} \;\leq\;  \ltwonorm{e^{(i)}}\label{eq2Conv2}
\end{eqnarray}
By subtracting the first equations of \eqref{sysNV} from that of \eqref{sys6V}, we get \eqref{eqConv1} by linearity of $\vec{V}$ operator.
\begin{eqnarray}
\ltwoinner{ e^{(k+1)},v} &=& \tau\ltwoinner{\vec{V}(\phi_N) \cdot \nabla\psi_N - \vec{V}(\phi_N^{(k+1)}) \cdot \nabla \psi_N^{(k)} - \vec{V}(\phi_N^{(k)}) \cdot \nabla\psi_N^{(k+1)} +\vec{V}(\phi_N^{(k)}) \cdot \nabla \psi_N^{(k)} - \vec{V}(p).\nabla(g^{(k+1)}),v}\nonumber \\
&=& \tau\ltwoinner{\vec{V}(g^{(k+1)}) \cdot \nabla\psi_N + \vec{V}(\phi_N^{(k+1)}) \cdot \nabla (e^{(k)}) +\vec{V}(\phi_N^{(k)}) \cdot \nabla (e^{(k+1)} -e^{(k)})  -  \vec{V}(p).\nabla(g^{(k+1)}),v} \label{eqConv1}\\
&=& \tau\ltwoinner{\vec{V}(g^{(k+1)}) \cdot \nabla\psi_N + \vec{V}(g^{(k)} - g^{(k+1)}) \cdot \nabla (e^{(k)}) +\vec{V}(\phi_N^{(k)}) \cdot \nabla (e^{(k+1)})  -  \vec{V}(p).\nabla(g^{(k+1)}),v} \label{eqConvv1}
\end{eqnarray}
Let $v = e^{(k+1)}$ in  \eqref{eqConvv1}, then
\begin{eqnarray}
\ltwonorm{e^{(k+1)}}^2 
&=&\tau\ltwoinner{\vec{V}(g^{(k+1)}) \cdot \nabla\psi_N + \vec{V}(g^{(k)}-g^{(k+1)}) \cdot \nabla e^{(k)} , e^{(k+1)}}  - \tau \ltwoinner{ \vec{V}(p).\nabla(g^{(k+1)}),e^{(k+1)}}\nonumber\\
&\leq & \tau\ltwonorm{\vec{V}(g^{(k+1)}) \cdot \nabla\psi_N + \vec{V}(g^{(k)}-g^{(k+1)}) \cdot \nabla e^{(k)} - \vec{V}(p).\nabla(g^{(k+1)})} \ltwonorm{e^{(k+1)}} \nonumber\\
\therefore \ltwonorm{e^{(k+1)}} &\leq & \tau\ltwonorm{\vec{V}(g^{(k+1)}) \cdot \nabla\psi_N} + \tau\ltwonorm{\vec{V}(g^{(k)}-g^{(k+1)}) \cdot \nabla e^{(k)} }+ \tau\ltwonorm{\vec{V}(p).\nabla(g^{(k+1)})} \label{eqConv3}
\end{eqnarray}
Similarly to \eqref{eqCia}, by Ciarlet (\cite{ciarlet}, Theorem 3.2.6), we have \eqref{eq:C1}-\eqref{eq:C33}
\begin{eqnarray}
\ltwonorm{\vec{V}(g^{(k+1}) \cdot \nabla \psi_N} &\leq &c_{inv}^2 h^{-5/2} \ltwonorm{\psi_N} \;.\; \ltwonorm{g^{(k+1)}} \;\leq \; c_{inv}^2 h^{-5/2} \ltwonorm{\psi_N} \;.\; \ltwonorm{e^{(k+1)}} \label{eq:C1}\\
\ltwonorm{\vec{V}(g^{(k)} - g^{(k+1}) \cdot \nabla e^{(k)}}  &\leq &c_{inv}^2 h^{-5/2} \ltwonorm{e^{(k)}} \;.\; \ltwonorm{g^{(k)}-g^{(k+1}}  =c_{inv}^2 h^{-5/2} \ltwonorm{e^{(k)}} \;.\; \ltwonorm{\phi_N^{(k+1)}-\phi_N^{(k)}} \nonumber\\
&\leq &c_{inv}^2 \epsilon_{tol} h^{-5/2} \ltwonorm{e^{(k)}}
\label{eq:C33}\\
 \ltwonorm{\vec{V}(p) \cdot \nabla  g^{(k+1)} } &\leq &  2 \loneinftynorm{p} \;.\;\honenorm{ g^{(k+1)} } \quad \mbox{Similarly to \eqref{eqP}\qquad} \nonumber \\
 &\leq &  2 \loneinftynorm{p} \;.\;\ltwonorm{ e^{(k+1)} }\qquad \mbox{By } \eqref{eq2Conv2}\quad \label{eq:C4}
\end{eqnarray} 
Replace equations \eqref{eq:C1}, \eqref{eq:C33} and \eqref{eq:C4} in \eqref{eqConv3} and using the apriori estimate on  $\ltwonorm{\psi_N}$ (Theorem \ref{thrm:apriori}), we get
\begin{eqnarray}
\ltwonorm{e^{(k+1)}} &\leq & \tau\ltwonorm{\vec{V}(g^{(k+1)}) \cdot \nabla\psi_N} + \tau\ltwonorm{\vec{V}(g^{(k)}-g^{(k+1)}) \cdot \nabla e^{(k)} }+ \tau\ltwonorm{\vec{V}(p).\nabla(g^{(k+1)})} \nonumber\\
&\leq & \tau c_{inv}^2 h^{-5/2} a_0 \,.\, \ltwonorm{e^{(k+1)}} + \tau c_{inv}^2 \epsilon_{tol} h^{-5/2} \ltwonorm{e^{(k)}}  + 2\tau   \loneinftynorm{p} \,.\,\ltwonorm{ e^{(k+1)} }\qquad
\nonumber\\
\therefore \ltwonorm{e^{(k+1)}} &\leq & \dfrac{ \tau c_{inv}^2 \epsilon_{tol} h^{-5/2} }{1 -   \tau c_{inv}^2 h^{-5/2} a_0  -  2\tau   \loneinftynorm{p}}\ltwonorm{e^{(k)}} \;=\; {c} \, \ltwonorm{e^{(k)}} \label{eqConv4}
\end{eqnarray}
where $a_0 = e^{3 T \loneinftynorm{p}} \ltwonorm{w_0}$, and $\tau \leq   \dfrac{ 1 }{ c_{inv}^2 h^{-5/2} a_0  + 2   \loneinftynorm{p}} = \dfrac{ h^{5/2} }{ c_{inv}^2 a_0  + 2  h^{5/2}  \loneinftynorm{p}}$. In addition, \\
${c}= \dfrac{ \tau c_{inv}^2 \epsilon_{tol} h^{-5/2} }{1 -   \tau c_{inv}^2 h^{-5/2} a_0  -  2\tau   \loneinftynorm{p}} <1$ if and only if $$ \tau < \dfrac{1}{c_{inv}^2 h^{-5/2}( \epsilon_{tol} + a_0) + 2\loneinftynorm{p}} = \dfrac{h^{5/2}}{c_{inv}^2 ( \epsilon_{tol} + a_0)  + 2h^{5/2}\loneinftynorm{p}}< \dfrac{ h^{5/2} }{ c_{inv}^2 a_0  + 2  h^{5/2}  \loneinftynorm{p}}.$$
Moreover, since $h<1$, then $\dfrac{h^{5/2}}{2D_1} < \dfrac{h^{5/2}}{D_1} < \dfrac{h^{5/2}}{c_{inv}^2 (\epsilon_{tol}  +  a_0 )+ 2h^{5/2}\loneinftynorm{p}}$. Thus, let $\tau < \dfrac{h^{5/2}}{2D_1} < \dfrac{h^{5/2}}{2D}$
 where $D_1 = c_{inv}^2 (\epsilon_{tol}  +  a_0) + 2\loneinftynorm{p} =  c_{inv}^2 \epsilon_{tol} +D$,
  which ends the proof.
\end{proof}
\noindent A consequence of Theorem \eqref{thrm:convNM2} is the local convergence of Newton's method.
\begin{theorem}\label{thrm:ConvNM3}
 Assume that $\phi_N^{(0)}$ is chosen such that $\forall k \geq k_0 \geq 0, \ltwonorm{\phi_N^{(k+1)} - \phi_N^{(k)}} < \epsilon_{tol}$. 
 Then, for \\$\tau \leq \min\left\{\dfrac{1}{6||p||_{1,\infty}}, \dfrac{h^2}{16{c}_{0,inv}\Mnorm{W(t)}}, \dfrac{h^{5/2}}{2D_1 }\right\}$,  Newton's method converges to the unique solution of \eqref{HMC-Disc-Comp}, 
$$\lim\limits_{k\rightarrow \infty} \phi_N^{(k)} = \phi_N \qquad and \qquad \lim\limits_{k\rightarrow \infty} \psi_N^{(k)} = \psi_N$$
where $D_1 := D_1(\Omega,T,p,w_0) =c_{inv}^2 (\epsilon_{tol}+ e^{3 T \loneinftynorm{p}} \ltwonorm{w_0}) + 2\loneinftynorm{p}$ and $h<1$.
 
\end{theorem}
\begin{remark} This additional assumption ($\forall k \geq k_0 \geq 0, \ltwonorm{\phi_N^{(k+1)} - \phi_N^{(k)}} < \epsilon_{tol}$ ) is computational in nature, since it could be used as a stopping criteria for any iterative method solving a nonlinear problem. 
\end{remark}
\begin{proof}
By theorem \ref{thrm:convNM2}, there exists $c<1$ such that
\begin{eqnarray}
\ltwonorm{e^{(k+1)}}^2 & \leq & c^2 \ltwonorm{e^{(k)}}^2 \qquad \qquad \qquad \qquad \qquad 
\implies \ltwonorm{e^{(k+1)}}^2 \; \leq \; c^{2(k+1)} \ltwonorm{e^{(0)}}^2 \nonumber\\
\lim\limits_{k\rightarrow \infty} \ltwonorm{e^{(k+1)}}^2  & \leq &\lim\limits_{k\rightarrow \infty}  c^{2(k+1)} \ltwonorm{e^{(0)}}^2 \;=\; 0  \quad \qquad \;\iff \quad \lim\limits_{k\rightarrow \infty} \ltwonorm{e^{(k+1)}}^2 = 0\nonumber\\
\ltwonorm{e^{(k+1)}}^2 + \ltwonorm{g^{(k+1)}}^2& \leq & 2\ltwonorm{e^{(k+1)}}^2 \;\leq \;  2c^2 \ltwonorm{e^{(k)}}^2 \qquad \quad 
\implies  \ltwonorm{e^{(k+1)}}^2 + \ltwonorm{g^{(k+1)}}^2 \; \leq \; 2 c^{2(k+1)} \ltwonorm{e^{(0)}}^2\nonumber\\
\lim\limits_{k\rightarrow \infty}  \ltwonorm{e^{(k+1)}}^2 + \ltwonorm{g^{(k+1)}}^2& \leq & \lim\limits_{k\rightarrow \infty} 2 c^{2(k+1)} \ltwonorm{e^{(0)}}^2 \;=\;0\qquad \quad \iff  \quad  \lim\limits_{k\rightarrow \infty}  \ltwonorm{g^{(k+1)}}^2 = 0\nonumber
\end{eqnarray}
$\therefore \lim\limits_{k\rightarrow \infty} \psi_N^{(k+1)} \;=\;  \psi_N \qquad and \qquad   \lim\limits_{k\rightarrow \infty} \phi_N^{(k+1)} \;=\;  \phi_N 
$
\end{proof}
\section{Variants of Newton's method}\label{sec:Var}
In this section we discuss two variants of Newton's method, Chord's method (section \ref{sec:CM}) and Modified Newton's method (section \ref{sec:MN}), that are less computationally  intensive than Newton's method .  
\subsection{Chord's Method}\label{sec:CM}
To avoid recomputing the Jacobian matrix at each Newton iteration, Chord's method approximates the solution of \eqref{HMC-Disc-Comp},  by solving system \eqref{chord} iteratively till convergence up to some given tolerance
\begin{align}
J_F(U^{(0)},W^{(0)}) \begin{bmatrix}U^{(k+1)}-U^{(k)}\\ W^{(k+1)}-W^{(k)}\end{bmatrix}&=-F(U^{(k)},W^{(k)}) \label{chord} 
\end{align}
where $(U^{(0)},W^{(0)})=(U(t),W(t))$, and $J_F(U^{(0)},W^{(0)})$   is the Jacobian matrix that is computed once per time iteration. 

By replacing $F(U^{(k)},W^{(k)})$ and $J_F(U^{(0)},W^{(0)})$ by their expressions, \eqref{JF} and \eqref{F} respectively, system \eqref{chord} is reduced to the linear system \eqref{eq:chord}, where computing the right-hand side vector requires 3 matrix-vector multiplications. However, using property \eqref{BS}, this can be reduced to just 2 matrix-vector multiplications, $S(U^{(k)})(W^{(0)} - W^{(k)})$ and $S(U^{(0)})W^{(k)}$, as shown in \eqref{eq:chord2}. Algorithm \eqref{alg:HMC-Chord} summarizes the procedure.
\begin{eqnarray}
\begin{bmatrix}
 \tau B(W^{(0)}) -\tau R & M+\tau S(U^{(0)})\\ 
 K & -M 
 \end{bmatrix}
  \begin{bmatrix}
  U^{(k+1)}\\ 
  W^{(k+1)}
  \end{bmatrix}&=&
   \begin{bmatrix} 
   \tau B(W^{(0)}) -\tau R &M+\tau S(U^{(0)})\\ 
   K & -M
    \end{bmatrix}
   \begin{bmatrix}
   U^{(k)}\\ 
   W^{(k)}
   \end{bmatrix}  -F(U^{(k)},W^{(k)}) \nonumber\\
  \iff
 \begin{bmatrix}
  \tau B(W^{(0)}) -\tau R & M+\tau S(U^{(0)})\\
   K  & -M  
   \end{bmatrix} 
   \begin{bmatrix}
   U^{(k+1)}\\
    W^{(k+1)}
    \end{bmatrix}
    &=&\begin{bmatrix}
     \tau B(W^{(0)})U^{(k)} -\tau R U^{(k)}+ (M+\tau S(U^{(0)}))W^{(k)}\\
      K U^{(k)} -MW^{(k)} 
      \end{bmatrix}\nonumber\\
      &&-\begin{bmatrix}
      (M+\tau S(U^{(k)})) W^{(k)}-\tau RU^{(k)}- Z \\
       KU^{(k)}-MW^{(k)}
        \end{bmatrix}\nonumber\\
        \label{eq:chord}
\begin{bmatrix} 
\tau B(W^{(0)}) -\tau R & M+\tau S(U^{(0)})\\
 K  & -M  
 \end{bmatrix} 
 \begin{bmatrix}
 U^{(k+1)}\\ 
 W^{(k+1)}
 \end{bmatrix}
 &=&
 \begin{bmatrix}
  \tau B(W^{(0)})U^{(k)} +\tau S(U^{(0)})W^{(k)}-\tau S(U^{(k)}) W^{(k)}+ Z \\
  0
  \end{bmatrix} \\
  \label{eq:chord2}
  &=&
 \begin{bmatrix}
  \tau S(U^{(k)})W^{(0)} +\tau S(U^{(0)})W^{(k)}-\tau S(U^{(k)}) W^{(k)}+ Z \\
  0
  \end{bmatrix} 
\end{eqnarray}

Note that \eqref{eq:chord2} is equivalent to system \eqref{mat1B} and \eqref{sys1B} where $\alpha = U^{(k+1)}, \beta =  W^{(k+1)}, W = W^{(0)}, U = U^{(0)},$ and the right hand side is replaced by  $M\gamma = \tau S(U^{(k)})W^{(0)}+\tau S(U^{(0)})W^{(k)}-\tau S(U^{(k)})W^{(k)}+Z$. Thus, the existence of a unique solution to \eqref{eq:chord2} is a corollary of theorem \ref{thrm:bddR} proven in section \ref{sec:exis}.

\subsubsection{Convergence}\label{sec:ConvC}
We are approximating the solution of the nonlinear system \eqref{HMC-Disc-Comp} by using Chord's method \eqref{chord} or equivalently \eqref{eq:chord2}, 
%
which can be expressed in variational form for any $v=\sum\limits_{I=1}^N v_I \varphi_I(x,y) \in X_{N,p}$ using \eqref{elem1}-\eqref{elem5}  and \eqref{def1}-\eqref{def2} as  

\begin{equation}\label{sysCV}
\begin{cases}
 \ltwoinner{ \psi_N^{(k+1)}  - \tau\vec{V}(\phi_N^{(k+1)}) \cdot \nabla \psi_N^{(0)}    - \tau \vec{V}(\phi_N^{(0)}) \cdot \nabla\psi_N^{(k+1)} + \tau \vec{V}(p).\nabla\phi_N^{(k+1)},v}   = \ltwoinner{\psi^{(0)}_N,v} + \tau  \ltwoinner{\vec{V}(\phi_N^{(k)}) \cdot \nabla \psi_N^{(k)}, v} \\
  \qquad \qquad \qquad \qquad \qquad \qquad \qquad \qquad \qquad \qquad \qquad \qquad  \qquad \qquad \qquad- \tau  \ltwoinner{\vec{V}(\phi_N^{(0)}) \cdot \nabla \psi_N^{(k)}, v} - \tau  \ltwoinner{\vec{V}(\phi_N^{(k)}) \cdot \nabla \psi_N^{(0)}, v}\\
\honeinner{\phi_N^{(i)},v} = \ltwoinner{\psi_N^{(i)},v} \qquad  \qquad \mbox{for } i=\{k,k+1\} 
\end{cases}
\end{equation}
Let $e^{(k)} = \psi_N - \psi_N^{(k)}$ and $g^{(k)} = \phi_N - \phi_N^{(k)}$, then we prove the local convergence of Chord's method.
\begin{theorem}\label{thrm:convC22} Assume that $\psi_N^{(0)}$ and $\phi_N^{(0)}$ are chosen such that $ \ltwonorm{\phi_N - \phi_N^{(0)}} < c_1$ and  $\forall k\geq k_0 \geq 0$,$$ \ltwonorm{\phi_N^{(k+1)} - \phi_N^{(k)}} < \epsilon_{tol1} \qquad \mbox{and} \qquad\ltwonorm{\psi_N^{(k+1)} - \psi_N^{(k)}} < \epsilon_{tol2}.$$ Then, for 
 $\tau \leq \min\left\{\dfrac{1}{6||p||_{1,\infty}}, \dfrac{h^2}{16{c}_{0,inv}\Mnorm{W(t)}}, \dfrac{h^{5/2}}{D_2 }\right\} = O(h^{2.5})$ there exists a constant $c<1$ such that $$\ltwonorm{e^{(k+1)}}^2 \; \leq \; c^2 \ltwonorm{e^{(k)}}^2$$
 where $D_2 := D_2(\Omega,T,p,w_0) =c_{inv}^2( 4c_1+ \epsilon_{tol1} + \epsilon_{tol2}+ e^{3 T \loneinftynorm{p}} \ltwonorm{w_0})+ 2\loneinftynorm{p}$, and $h<1$.
\end{theorem}
\begin{proof}
Similarly to theorem \ref{thrm:NewtonConv}, and as a consequence of theorem \ref{thrm:bddR}, \eqref{sysCV} has a unique solution $\{\phi_N^{(k+1)}, \psi_N^{(k+1)}\}$ for 
 $\tau \leq \min\left\{\dfrac{1}{6||p||_{1,\infty}}, \dfrac{h^2}{16{c}_{0,inv}\Mnorm{W(t)}}, \dfrac{h^{5/2}}{2D}\right\} $, where  $D  = c_{inv}^2 e^{3 T \loneinftynorm{p}} \ltwonorm{w_0} +2\loneinftynorm{p}$ and $h<1$.
Then, by subtracting the second equation of \eqref{sysCV} from that of \eqref{sys6V}, we get \eqref{eq2ConvC1} for $i = \{k,k+1\}$. Letting $v = g^{(i)}$ we get  \eqref{eq2ConvC2}
\begin{eqnarray}
\honeinner{g^{(i)},v} &=& \ltwoinner{e^{(i)},v} \label{eq2ConvC1}\\
  \honenorm{g^{(i)}}^2 &=& \ltwoinner{e^{(i)},g^{(i)}} \;\leq \;  \ltwonorm{e^{(i)}} \ltwonorm{g^{(i)}} \;\leq \;  \ltwonorm{e^{(i)}} \honenorm{g^{(i)}} \nonumber\\
\therefore \ltwonorm{g^{(i)}} &\leq & \honenorm{g^{(i)}} \;\leq\;  \ltwonorm{e^{(i)}}\label{eq2ConvC2}
\end{eqnarray}
By subtracting the first equations of \eqref{sysCV} from that of \eqref{sys6V}, we get \eqref{eqConvC1} by linearity of $\vec{V}$ operator.
\begin{eqnarray}
\ltwoinner{ e^{(k+1)},v} &=& \tau\ltwoinner{\vec{V}(\phi_N) \cdot \nabla\psi_N - \vec{V}(\phi_N^{(k+1)}) \cdot \nabla \psi_N^{(0)} - \vec{V}(\phi_N^{(0)}) \cdot \nabla\psi_N^{(k+1)} -\vec{V}(\phi_N^{(k)}) \cdot \nabla \psi_N^{(k)} - \vec{V}(p).\nabla(g^{(k+1)}),v}\nonumber \\
&& + \tau \ltwoinner{\vec{V}(\phi_N^{(0)}) \cdot \nabla \psi_N^{(k)}+\vec{V}(\phi_N^{(k)}) \cdot \nabla \psi_N^{(0)}, v}\label{eqConvC1} \\
&=& \tau\ltwoinner{\vec{V}(g^{(k+1)}) \cdot \nabla\psi_N + \vec{V}(\phi_N^{(k+1)}) \cdot \nabla (e^{(0)}) +\vec{V}(\phi_N^{(0)}) \cdot \nabla (e^{(k+1)}-e^{(k)})   -  \vec{V}(p).\nabla(g^{(k+1)}),v} \nonumber\\
&& + \tau \ltwoinner{
\vec{V}(\phi_N^{(k)}) \cdot \nabla (e^{(k)} - e^{(0)})
, v} \nonumber\\
&=&\tau\ltwoinner{\vec{V}(g^{(k+1)}) \cdot \nabla\psi_N + \vec{V}(g^{(k)}-g^{(k+1)}) \cdot \nabla e^{(0)} +\vec{V}(\phi_N^{(0)}) \cdot \nabla (e^{(k+1)}-e^{(k)})   -  \vec{V}(p).\nabla(g^{(k+1)}),v} \nonumber\\
\nonumber\\
&&+\tau \ltwoinner{ \vec{V}(\phi_N) \cdot \nabla e^{(k)} -
\vec{V}(g^{(k)}) \cdot \nabla e^{(k)} , v} 
\label{eqConvvC1}
\end{eqnarray}
Note that
\begin{eqnarray}
- \ltwoinner{ \vec{V}(g^{(k)}) \cdot \nabla e^{(k)} , v}  &=&  \ltwoinner{ \vec{V}(g^{(k+1)} - g^{(k)}) \cdot \nabla e^{(k)} , v}  +  \ltwoinner{ \vec{V}(g^{(k+1)}) \cdot \nabla ( e^{(k+1)} -  e^{(k)}) , v}\nonumber\\
&& -  \ltwoinner{ \vec{V}(g^{(k+1)}) \cdot \nabla e^{(k+1)} , v} \label{eq:V1}\\
\ltwoinner{ \vec{V}(\phi_N^{(0)}) \cdot \nabla (e^{(k+1)}-e^{(k)}) , v} &=& \ltwoinner{ \vec{V}(g^{(0)}) \cdot \nabla (e^{(k)}-e^{(k+1)}) ,v} + \ltwoinner{ \vec{V}(\phi_N) \cdot \nabla (e^{(k+1)}-e^{(k)}) , v} \qquad \qquad \label{eq:V2}
\end{eqnarray}
Replacing \eqref{eq:V1} and  \eqref{eq:V2} in \eqref{eqConvvC1} we get
\begin{eqnarray}
\ltwoinner{ e^{(k+1)},v} 
 &=&\tau\ltwoinner{\vec{V}(g^{(k+1)}) \cdot \nabla\psi_N + \vec{V}(g^{(k)}-g^{(k+1)}) \cdot \nabla e^{(0)} +\vec{V}(g^{(0)}) \cdot \nabla (e^{(k)}-e^{(k+1)})   -  \vec{V}(p).\nabla(g^{(k+1)}),v} \nonumber\\
&&+\tau \ltwoinner{ \vec{V}(g^{(k+1)} - g^{(k)}) \cdot \nabla e^{(k)}+ \vec{V}(g^{(k+1)}) \cdot \nabla ( e^{(k+1)} -  e^{(k)}) , v}\nonumber\\
&& - \tau \ltwoinner{ \vec{V}(g^{(k+1)}) \cdot \nabla e^{(k+1)} , v} +\tau \ltwoinner{ \vec{V}(\phi_N) \cdot \nabla e^{(k+1)} ,v}
\label{eqConvvC2}
\end{eqnarray}

\noindent Let $v = e^{(k+1)}$ in  \eqref{eqConvvC2}, then
\begin{eqnarray}
\ltwonorm{e^{(k+1)}} 
&\leq &\tau\ltwonorm{\vec{V}(g^{(k+1)}) \cdot \nabla\psi_N + \vec{V}(g^{(k)}-g^{(k+1)}) \cdot \nabla e^{(0)} +\vec{V}(g^{(0)}) \cdot \nabla (e^{(k)}-e^{(k+1)})   -  \vec{V}(p).\nabla(g^{(k+1)})} \nonumber\\
\nonumber\\
&&+\tau \ltwonorm{ \vec{V}(g^{(k+1)} - g^{(k)}) \cdot \nabla e^{(k)}+\vec{V}(g^{(k+1)}) \cdot \nabla ( e^{(k+1)} -  e^{(k)}) } \nonumber\\
&\leq &\tau\ltwonorm{\vec{V}(g^{(k+1)}) \cdot \nabla\psi_N} + \tau\ltwonorm{\vec{V}(g^{(k)}-g^{(k+1)}) \cdot \nabla e^{(0)}} +\tau\ltwonorm{\vec{V}(g^{(0)}) \cdot \nabla (e^{(k)}-e^{(k+1)})}    
\nonumber\\
&&+ \tau \ltwonorm{\vec{V}(g^{(k+1)} - g^{(k)}) \cdot \nabla e^{(k)}}+\tau \ltwonorm{\vec{V}(g^{(k+1)}) \cdot \nabla ( e^{(k+1)} -  e^{(k)}) }  + \tau\ltwonorm{\vec{V}(p).\nabla(g^{(k+1)})}\label{eqConvC3}
\end{eqnarray}
Similarly to \eqref{eqCia}, by  \eqref{ineq:ciar} (\cite{ciarlet}, Theorem 3.2.6), we have \eqref{eq:CC1}-\eqref{eq:CC44}
\begin{eqnarray}
\ltwonorm{\vec{V}(g^{(k+1}) \cdot \nabla \psi_N} &\leq &c_{inv}^2 h^{-5/2} \ltwonorm{\psi_N} \;.\; \ltwonorm{g^{(k+1)}} \;\leq \; c_{inv}^2 h^{-5/2} \ltwonorm{\psi_N} \;.\; \ltwonorm{e^{(k+1)}} \label{eq:CC1}\\
\ltwonorm{\vec{V}(g^{(k)} - g^{(k+1}) \cdot \nabla e^{(0)}}  &\leq &c_{inv}^2 h^{-5/2} \ltwonorm{e^{(0)}} \ltwonorm{g^{(k)} - g^{(k+1)}} \;\leq \;c_{inv}^2  h^{-5/2}  c_1 \left( \ltwonorm{e^{(k)}}+\ltwonorm{e^{(k+1)}}\right) \qquad\\
\ltwonorm{\vec{V}(g^{(0)}) \cdot \nabla (e^{(k)}-e^{(k+1)})}    &\leq &c_{inv}^2  h^{-5/2}  \ltwonorm{g^{(0)}} \ltwonorm{e^{(k)} - e^{(k+1)}}\;\leq \; c_{inv}^2  h^{-5/2} c_1\left( \ltwonorm{e^{(k)}}+\ltwonorm{e^{(k+1)}}\right)\\
 \ltwonorm{\vec{V}(g^{(k+1)} - g^{(k)}) \cdot \nabla e^{(k)}}&\leq & c_{inv}^2 h^{-5/2} \ltwonorm{g^{(k+1)} - g^{(k)}} \;.\; \ltwonorm{e^{(k)}} \;\leq \;  c_{inv}^2 h^{-5/2}  \epsilon_{tol1} \ltwonorm{e^{(k)}} \\
  \ltwonorm{\vec{V}(g^{(k+1)}) \cdot \nabla ( e^{(k+1)} -  e^{(k)}) }&\leq & c_{inv}^2 h^{-5/2} \ltwonorm{g^{(k+1)}} \;.\; \ltwonorm{ e^{(k+1)} -e^{(k)}}\;\leq \;  c_{inv}^2 h^{-5/2}  \epsilon_{tol2} \ltwonorm{g^{(k+1)}} \qquad \label{eq:CC44}\\ 
  \ltwonorm{\vec{V}(p) \cdot \nabla  g^{(k+1)} } &\leq &  2 \loneinftynorm{p} \;.\;\honenorm{ g^{(k+1)} } \quad \mbox{Similarly to \eqref{eqP}\qquad} \nonumber \\
 &\leq &  2 \loneinftynorm{p} \;.\;\ltwonorm{ e^{(k+1)} }\qquad \mbox{By } \eqref{eq2ConvC2}\quad \label{eq:CC4}
\end{eqnarray} 
Replace equations \eqref{eq:CC1} - \eqref{eq:CC4} in \eqref{eqConvC3}, then 
\begin{eqnarray}
\ltwonorm{e^{(k+1)}} &\leq & 2\tau \loneinftynorm{p} \;.\;\ltwonorm{ e^{(k+1)} }+\tau c_{inv}^2 h^{-5/2} \left[ a_0 \ltwonorm{e^{(k+1)}} +2c_1\left( \ltwonorm{e^{(k)}}+\ltwonorm{e^{(k+1)}}\right) \right.\nonumber \\
&&\left. +\epsilon_{tol1} \ltwonorm{e^{(k)}} +\epsilon_{tol2} \ltwonorm{e^{(k+1)}}\right]\nonumber \\
\therefore\ltwonorm{e^{(k+1)}} &\leq & \dfrac{\tau c_{inv}^2 h^{-5/2} (2c_1+\epsilon_{tol1})}{1 - 2\tau \loneinftynorm{p}-\tau c_{inv}^2 h^{-5/2}(a_0+2c_1+\epsilon_{tol2})}  \ltwonorm{e^{(k)}} \;=\; c \ltwonorm{e^{(k)}}
 \label{eqConvC4}
\end{eqnarray}
where $a_0 = e^{3 T \loneinftynorm{p}} \ltwonorm{w_0}$, 
$\tau \leq \dfrac{ 1 }{  c_{inv}^2 h^{-5/2} (a_0 +2c_1+\epsilon_{tol2}) + 2   \loneinftynorm{p}} = \dfrac{  h^{5/2}  }{  c_{inv}^2(a_0 +2c_1+\epsilon_{tol2}) + 2    h^{5/2} \loneinftynorm{p}}$. In addition,  
${c}= \dfrac{\tau c_{inv}^2 h^{-5/2} (2c_1+\epsilon_{tol1})}{1 - 2\tau \loneinftynorm{p}-\tau c_{inv}^2 h^{-5/2}(a_0+2c_1+\epsilon_{tol2}) } \;<\;1, $ if and only if $$ \tau < \dfrac{1}{ h^{-5/2} \tilde{D}_2  + 2\loneinftynorm{p}} = \dfrac{h^{5/2}}{\tilde{D}_2  + 2h^{5/2}\loneinftynorm{p}} < \dfrac{  h^{5/2}  }{  c_{inv}^2(a_0 +2c_1+\epsilon_{tol2}) + 2    h^{5/2} \loneinftynorm{p}} $$
 where $\tilde{D}_2 = c_{inv}^2 (a_0 + 4c_1+ \epsilon_{tol1} + \epsilon_{tol2})$. Moreover, since $h<1$, then $\dfrac{h^{5/2}}{2D_2} < \dfrac{h^{5/2}}{D_2} < \dfrac{h^{5/2}}{\tilde{D}_2+ 2h^{5/2}\loneinftynorm{p}}$ where $D_2 = \tilde{D}_2 + 2\loneinftynorm{p}$. Thus, let $\tau < \dfrac{h^{5/2}}{2D_2} < \dfrac{h^{5/2}}{2D}$
 where $D = c_{inv}^2  a_0 + 2\loneinftynorm{p}$,
  which ends the proof.
\end{proof}
\noindent A corollary of Theorem \eqref{thrm:convC22} is the local convergence of Chord's method.
\begin{theorem}\label{thrm:ConvC}
Assume that $\psi_N^{(0)}$ and $\phi_N^{(0)}$ are chosen such that $ \ltwonorm{\phi_N - \phi_N^{(0)}} < c_1$ and  $\forall k\geq k_0 \geq 0$,$$ \ltwonorm{\phi_N^{(k+1)} - \phi_N^{(k)}} < \epsilon_{tol1} \qquad \mbox{and} \qquad\ltwonorm{\psi_N^{(k+1)} - \psi_N^{(k)}} < \epsilon_{tol2}.$$ Then, for 
 $\tau \leq \min\left\{\dfrac{1}{6||p||_{1,\infty}}, \dfrac{h^2}{16{c}_{0,inv}\Mnorm{W(t)}}, \dfrac{h^{5/2}}{D_2 }\right\} = O(h^{2.5})$ 
  Chord's method converges to the unique solution of \eqref{HMC-Disc-Comp},
$$\lim\limits_{k\rightarrow \infty} \phi_N^{(k)} = \phi_N \qquad and \qquad \lim\limits_{k\rightarrow \infty} \psi_N^{(k)} = \psi_N$$
where $D_2 := D_2(\Omega,T,p,w_0) =c_{inv}^2( 4c_1+ \epsilon_{tol1} + \epsilon_{tol2}+ e^{3 T \loneinftynorm{p}} \ltwonorm{w_0})+ 2\loneinftynorm{p}$, and $h<1$.
\end{theorem}
\begin{proof}
By theorem \ref{thrm:convC22}, there exists $c<1$ such that
\begin{eqnarray}
\ltwonorm{e^{(k+1)}}^2 & \leq & c^2 \ltwonorm{e^{(k)}}^2 \qquad \qquad \qquad \qquad \qquad 
\implies \ltwonorm{e^{(k+1)}}^2 \; \leq \; c^{2(k+1)} \ltwonorm{e^{(0)}}^2 \nonumber\\
\lim\limits_{k\rightarrow \infty} \ltwonorm{e^{(k+1)}}^2  & \leq &\lim\limits_{k\rightarrow \infty}  c^{2(k+1)} \ltwonorm{e^{(0)}}^2 \;=\; 0  \quad \qquad \;\iff \quad \lim\limits_{k\rightarrow \infty} \ltwonorm{e^{(k+1)}}^2 = 0\nonumber\\
\ltwonorm{e^{(k+1)}}^2 + \ltwonorm{g^{(k+1)}}^2& \leq & 2\ltwonorm{e^{(k+1)}}^2 \;\leq \;  2c^2 \ltwonorm{e^{(k)}}^2 \qquad \quad 
\implies  \ltwonorm{e^{(k+1)}}^2 + \ltwonorm{g^{(k+1)}}^2 \; \leq \; 2 c^{2(k+1)} \ltwonorm{e^{(0)}}^2\nonumber\\
\lim\limits_{k\rightarrow \infty}  \ltwonorm{e^{(k+1)}}^2 + \ltwonorm{g^{(k+1)}}^2& \leq & \lim\limits_{k\rightarrow \infty} 2 c^{2(k+1)} \ltwonorm{e^{(0)}}^2 \;=\;0\qquad \quad \iff  \quad  \lim\limits_{k\rightarrow \infty}  \ltwonorm{g^{(k+1)}}^2 = 0\nonumber
\end{eqnarray}
$\therefore \lim\limits_{k\rightarrow \infty} \psi_N^{(k+1)} \;=\;  \psi_N \qquad and \qquad   \lim\limits_{k\rightarrow \infty} \phi_N^{(k+1)} \;=\;  \phi_N 
$
\end{proof}

\subsection{Modified Newton's Method}\label{sec:MN}
In both Newton and Chord's methods, the computation of the Jacobian matrix requires the computation of the $B(W)$ matrix that consists of $N$ matrix-vector multiplications. To avoid the computation of the $B(W)$ matrix, we introduce a Modified Newton's Method. Starting with Newton's equation \eqref{New}  or equivalently \eqref{eq:3.1}
 \begin{equation}
\label{eq:3.1}
\iff
\left\{\begin{array}{ll} \tau B(W^{(k)})U^{(k+1)} -\tau RU^{(k+1)}+MW^{(k+1)}+\tau S(U^{(k)})W^{(k+1)}=\tau S(U^{(k)})W^{(k)}+Z\\
KU^{(k+1)}-MW^{(k+1)}=0
\end{array}
\right. 
\end{equation}
and using property \eqref{BS}, we approximate $B(W^{(k)})U^{(k+1)}$ by $S(U^{(k)})W^{(k)}$,
 $$B(W^{(k)})U^{(k+1)}=S(U^{(k+1)})W^{(k)}\approx S(U^{(k)})W^{(k)}$$
 leading to the Modified Newton system \eqref{ModNew} whose right-hand side is fixed throughout the Modified Newton iterations. Moreover, the modified Jacobian matrix $\tilde{J_F}(U,W)$ \eqref{ModNew} requires updating one of its blocks by computing $S(U^{(k)})$ at each Modified Newton iteration. The procedure is summarized in Algorithm \eqref{alg:HMC-Modified Newton}.
 \begin{eqnarray}
&&\left\{\begin{array}{ll} \tau S(U^{(k)}) W^k  -\tau RU^{(k+1)}+MW^{(k+1)}+\tau S(U^{(k)})W^{(k+1)})=\tau S(U^{(k)})W^{(k)}+Z\\
KU^{(k+1)}-MW^{(k+1)}=0
\end{array}
\right.\nonumber \\
\iff &&
\left\{\begin{array}{ll}  -\tau RU^{(k+1)}+MW^{(k+1)}+\tau S(U^{(k)})W^{(k+1)}=Z\\
KU^{(k+1)}-MW^{(k+1)}=0
\end{array}
\right. \nonumber\\
 \iff && \label{ModNew}
 \begin{bmatrix}  -\tau R & M+\tau S(U^{(k)})\\ K  & -M  \end{bmatrix} \begin{bmatrix}U^{(k+1)}\\ W^{(k+1)}\end{bmatrix}=\begin{bmatrix} Z\\0\end{bmatrix}
 \end{eqnarray}
Note that the Modified Newton method, defined by \eqref{ModNew}, is equivalent to the iterative solution ($x^{(k+1)} = G(x^{(k)})$) of the fixed point problem of \eqref{HMC-Disc-Comp} or equivalently \eqref{eq:18}, i.e.
$$ \begin{bmatrix}U\\ W\end{bmatrix}=\begin{bmatrix}  -\tau R & M+\tau S(U)\\ K  & -M  \end{bmatrix}^{-1}\begin{bmatrix} Z\\0\end{bmatrix} = \tilde{J_F}(U,W)^{-1}\begin{bmatrix} Z\\0\end{bmatrix} = G([U,\, W]^T)$$
where we prove in section \eqref{sec:exis2} that the Modified Jacobian matrix $\tilde{J_F}(U,W)$ is invertible.

\noindent At every iteration of the Modified Newton's method, there is a need to solve a system of form \eqref{mat1B22}, where $[\alpha, \beta]^T\in \mathbb{R}^{2N}$. System  \eqref{mat1B22} is equivalent to system \eqref{sys1B22}.
 \begin{eqnarray}
 &&\tilde{J}_F(U,W) \begin{bmatrix}
 \alpha\\
 \beta
 \end{bmatrix} 
 = \begin{bmatrix}  -\tau R & M+\tau S(U)\\ K & -M \end{bmatrix} \begin{bmatrix}
 \alpha\\
 \beta
 \end{bmatrix}
 = \begin{bmatrix}
 M\gamma\\
0
 \end{bmatrix} \label{mat1B22}\\
  &&\iff 
 \begin{cases}
 -\tau R \alpha +  M\beta +\tau S(U)\beta = M\gamma&\\
  K\alpha =  M\beta &
 \end{cases}\label{sys1B22}
 \end{eqnarray}
where at the $k+1^{th}$ iteration $\alpha = U^{(k+1)}, \beta =  W^{(k+1)}, \gamma = W^{(0)}, W = W^{(k)},$ and $U = U^{(k)}$ based on \eqref{ModNew}. 

\noindent To prove the convergence of this method (section \ref{sec:ConvMN}), we prove first the existence of a unique solution of system   \eqref{mat1B22} in section \ref{sec:exis2}.

\subsubsection{Existence of a Unique Solution to \eqref{mat1B22}}\label{sec:exis2}
To prove the existence of a unique solution to \eqref{mat1B22}, we start by showing that there exists some $C \in \mathbb{R}$ independent of $\tau$ and $h$, such that $\Mnorm{\alpha}^2 + \Mnorm{\beta}^2 \leq C\Mnorm{\gamma}^2$  
using variational formulation.\\
\noindent Let $\phi_N(x,y)  = \sum\limits_{I=1}^N \alpha_I \varphi_I(x,y)$, $\psi_N(x,y)  = \sum\limits_{I=1}^N \beta_I \varphi_I(x,y)$, and $\xi_N(x,y)  = \sum\limits_{I=1}^N \gamma_I \varphi_I(x,y)$, 
then system \eqref{sys1B22} can be expressed in variational form elementwise (for $1\leq I\leq N$) as  \eqref{sys1elemB22} based on \eqref{elem1}-\eqref{elem5}.
\begin{equation}\label{sys1elemB22}
\begin{cases}
  \tau \ltwoinner{\vec{V}(p).\nabla\phi_N, \varphi_{I}} + \ltwoinner{\psi_N,\varphi_I} - \tau\ltwoinner{\vec{V}(u_N) \cdot \nabla\psi_N,\varphi_I} = \ltwoinner{\xi_N,\varphi_I} &\\
\honeinner{\phi_N,\varphi_I} = \ltwoinner{\psi_N,\varphi_I}&\\
\end{cases}
\end{equation}

\noindent Moreover, $\Mnorm{\beta}^2 = \ltwonorm{\psi_N}^2$ by \eqref{elem1}. Similarly $\Mnorm{\alpha}^2 = \ltwonorm{\phi_N}^2$ and $\Mnorm{\gamma}^2 = \ltwonorm{\xi_N}^2$. Thus, we need to show that 
\begin{equation}\label{eqTP22}
 \ltwonorm{\phi_N}^2 +\ltwonorm{\psi_N}^2 \leq  C\ltwonorm{\xi_N}^2
 \end{equation}
For any $v=\sum\limits_{I=1}^N v_I \varphi_I(x,y) \in X_{N,p}$, system \eqref{sys1elemB22} can be written as 
\begin{equation}\label{sys1vB22}
\begin{cases}
 \ltwoinner{ \psi_N   - \tau \vec{V}(u_N) \cdot \nabla\psi_N + \tau \vec{V}(p).\nabla\phi_N,v}   = \ltwoinner{\xi_N,v} &\\
\honeinner{\phi_N,v} = \ltwoinner{\psi_N,v} &\\
\end{cases}
\end{equation}
\begin{theorem}\label{thrm:bdd22}
Let  $ \tau \leq \min\left\{\dfrac{1}{4\loneinftynorm{p}}, \dfrac{h^2}{16{c}_{0,inv}\Mnorm{W(t)}}\right\} = O(h^2)
 $ 
then \vspace{-2mm} 
\begin{equation}
\ltwonorm{\phi_N}^2 +\ltwonorm{\psi_N}^2 \leq  8\ltwonorm{\xi_N}^2.
\end{equation}
\end{theorem}
\begin{proof} 
Similarly to \eqref{eq23B}, By setting $v = \phi_N$ in the second equation of system \eqref{sys1vB22} and using Cauchy-Schwarz we get 
\begin{eqnarray}
\implies \ltwonorm{\phi_N} &\leq & \ltwonorm{\psi_N}   \label{eq23B22}\\
 \therefore  \ltwonorm{\phi_N}^2 +\ltwonorm{\psi_N}^2 &\leq & 2\ltwonorm{\psi_N}^2 \qquad \qquad \label{eq24B22}
\end{eqnarray}
Thus, to obtain our result we must upper bound $\ltwonorm{\psi_N}^2$ in terms of $\ltwonorm{\xi_N}^2$. \\
Let $v = \psi_N$ in the first equation of system \eqref{sys1vB22} we get \eqref{eqq1B22}. Then, using Cauchy-Schwarz we get \eqref{eqq2B22}.
\begin{eqnarray}
  \ltwoinner{\xi_N,\psi_N} &=&  \ltwoinner{ \psi_N    - \tau \vec{V}(u_N) \cdot \nabla\psi_N + \tau \vec{V}(p).\nabla\phi_N,\psi_N}   \label{eqq1B22}\\
\ltwonorm{\psi_N}^2 &=& \ltwoinner{\xi_N,\psi_N}  + \ltwoinner{  \tau \vec{V}(u_N) \cdot \nabla\psi_N,\psi_N} - \ltwoinner{  \tau \vec{V}(p).\nabla\phi_N,\psi_N} \nonumber\\
  &=& \ltwoinner{\xi_N,\psi_N} - \ltwoinner{  \tau \vec{V}(p).\nabla\phi_N,\psi_N}\qquad \qquad \mbox{using skew-symmetry} \nonumber\\
&\leq & \ltwonorm{\xi_N}  \ltwonorm{\psi_N}  + \tau \ltwonorm{\vec{V}(p).\nabla\phi_N} \ltwonorm{\psi_N}  \nonumber\\
\therefore  \ltwonorm{\psi_N} &\leq & \ltwonorm{\xi_N}  + \tau \ltwonorm{\vec{V}(p).\nabla\phi_N} \label{eqq2B22}
\end{eqnarray}
Assuming $p\in C^{\infty}$, we upper bound the last terms of \eqref{eqq2B22} in terms of $ \ltwonorm{\psi_N}$ similarly to \eqref{eq3B}.
\begin{eqnarray}
\implies \ltwonorm{\vec{V}(p) \cdot \nabla \phi_N } 
&\leq &  2 \loneinftynorm{p} \honenorm{\phi_{N} } 
\; \leq \;   2 \loneinftynorm{p} \ltwonorm{\psi_N} \label{eq3B22}
\end{eqnarray}
Replacing \eqref{eq3B22} in \eqref{eqq2B22} we get 
\begin{eqnarray}
 \ltwonorm{\psi_N} &\leq & \ltwonorm{\xi_N}  + 2\tau   \loneinftynorm{p} \ltwonorm{\psi_N} \nonumber\\
 \implies \ltwonorm{\psi_N} &\leq & \dfrac{1}{\tilde{C}}  \ltwonorm{\xi_N} \nonumber\\
 \therefore   \ltwonorm{\phi_N}^2 +\ltwonorm{\psi_N}^2 &\leq & 2\ltwonorm{\psi_N}^2 \;\leq \;  \dfrac{2}{\tilde{C}^2}  \ltwonorm{\xi_N}^2 
\end{eqnarray}
where  $C =  \dfrac{2}{\tilde{C}^2}$ in \eqref{eqTP22}, $\tilde{C} = 1 - 2\tau  \loneinftynorm{p}$.
Thus, if $\tau \leq \dfrac{1}{4 \loneinftynorm{p}}$, then   $\tilde{C} \geq  \dfrac{1}{2}$, and  therefore
$
C \;=\;  \dfrac{2}{\tilde{C}^2} \;\leq\;8
$.
\end{proof}
\noindent A corollary of Theorem \eqref{thrm:bdd22} is the existence of a unique solution to system \eqref{mat1B22}.
\begin{theorem} Let  $ \tau \leq \min\left\{\dfrac{1}{4\loneinftynorm{p}}, \dfrac{h^2}{16{c}_{0,inv}\Mnorm{W(t)}}\right\} = O(h^2) $  then system \eqref{mat1B22} has a unique solution.
\end{theorem}
\begin{proof}
Let $\gamma = 0$, then $\xi_N = 0$ and by theorem \ref{thrm:bdd22} $$\ltwonorm{\phi_N}^2 +\ltwonorm{\psi_N}^2 \leq  0$$ 
for  $\tau \leq \dfrac{1}{4\loneinftynorm{p}}$.
\noindent Thus, $\ltwonorm{\phi_N}^2 = \Mnorm{\alpha}^2 = 0$ and   $\ltwonorm{\psi_N}^2 = \Mnorm{\beta}^2 = 0$, implying that $\alpha = \beta = 0$. \\Thus, $Null\{\tilde{J}_F(U,W)\} = \{0\}$, implying that $\tilde{J}_F(U,W)$ is invertible and system \eqref{mat1B22} has a unique solution.
\end{proof}
\subsubsection{Convergence}\label{sec:ConvMN}
We are approximating the solution of the nonlinear system \eqref{HMC-Disc-Comp} by using Modified Newton's method \eqref{ModNew}, 
%
which can be expressed in variational form for any $v=\sum\limits_{I=1}^N v_I \varphi_I(x,y) \in X_{N,p}$ using \eqref{elem1}-\eqref{elem5} and \eqref{def1}-\eqref{def2} as  
%
\begin{equation}\label{sysMNV}
\begin{cases}
 \ltwoinner{ \psi_N^{(k+1)}  - \tau\vec{V}(\phi_N^{(k)}) \cdot \nabla \psi_N^{(k+1)}   + \tau \vec{V}(p).\nabla\phi_N^{(k+1)},v}   = \ltwoinner{\psi_N^{(0)},v} \\
\honeinner{\phi_N^{(i)},v} = \ltwoinner{\psi_N^{(i)},v} \qquad  \qquad \mbox{for } i=\{k,k+1\} \\
\end{cases}
\end{equation}
Let $e^{(k)} = \psi_N - \psi_N^{(k)}$ and $g^{(k)} = \phi_N - \phi_N^{(k)}$, then we prove the convergence of Modified Newton's method.
\begin{theorem}\label{thrm:ConvMN}
Let $\tau \leq \min\left\{\dfrac{1}{6||p||_{1,\infty}}, \dfrac{h^2}{16{c}_{0,inv}\Mnorm{W(t)}}, \dfrac{h^{5/2}}{D_3}\right\} = O(h^{2.5})$, 
then there exists a constant $c<1$ such that  
\begin{equation}
\ltwonorm{e^{(k+1)}}^2 + \ltwonorm{g^{(k+1)}}^2 \leq c^2 \left( \ltwonorm{e^{(k)}}^2 + \ltwonorm{g^{(k)}}^2\right)
\end{equation} where $D_3 := D_3(\Omega,p,T,w_0) = c_{inv}^2 e^{3 T \loneinftynorm{p}} \ltwonorm{w_0} +2 \loneinftynorm{p}$, and $h<1$.
\end{theorem}
\begin{proof}By theorem \ref{thrm:bdd22}, \eqref{sysMNV} has a unique solution $\{\phi_N^{(k+1)}, \psi_N^{(k+1)}\}$ for 
 $\tau \leq \min\left\{\dfrac{1}{4||p||_{1,\infty}}, \dfrac{h^2}{16{c}_{0,inv}\Mnorm{W(t)}}\right\}$.\\
Then, by subtracting the second equation of \eqref{sysMNV} from that of \eqref{sys6V}, we get \eqref{eq2Conv1M} for $i = \{k,k+1\}$. Letting $v = g^{(i)}$ we get  \eqref{eq2Conv2M}
\begin{eqnarray}
\honeinner{g^{(i)},v} &=& \ltwoinner{e^{(i)},v} \label{eq2Conv1M}\\
  \honenorm{g^{(i)}}^2 &=& \ltwoinner{e^{(i)},g^{(i)}} \;\leq \;  \ltwonorm{e^{(i)}} \ltwonorm{g^{(i)}} \;\leq \;  \ltwonorm{e^{(i)}} \honenorm{g^{(i)}} \nonumber\\
\therefore \ltwonorm{g^{(i)}} &\leq & \honenorm{g^{(i)}} \;\leq\;  \ltwonorm{e^{(i)}}\label{eq2Conv2M}
\end{eqnarray}
By subtracting the first equations of \eqref{sysMNV} from that of \eqref{sys6V}, we get \eqref{eqConv1M} by linearity of $\vec{V}(.)$ operator.
\begin{eqnarray}
\ltwoinner{ e^{(k+1)},v} &=& \tau\ltwoinner{\vec{V}(\phi_N) \cdot \nabla\psi_N  - \vec{V}(\phi_N^{(k)}) \cdot \nabla\psi_N^{(k+1)}  - \vec{V}(p).\nabla(g^{(k+1)}),v}\nonumber \\
&=& \tau\ltwoinner{\vec{V}(g^{(k)}) \cdot \nabla\psi_N + \vec{V}(\phi_N^{(k)}) \cdot \nabla (e^{(k+1)})   -  \vec{V}(p).\nabla(g^{(k+1)}),v} \label{eqConv1M}
\end{eqnarray}
Let $v = e^{(k+1)}$ in  \eqref{eqConv1M}, then
\begin{eqnarray}
\ltwonorm{e^{(k+1)}}^2 &=&\tau\ltwoinner{\vec{V}(g^{(k)}) \cdot \nabla\psi_N,e^{(k+1)}} + \tau\ltwoinner{\vec{V}(\phi_N^{(k)}) \cdot \nabla (e^{(k+1)}),e^{(k+1)}}   -  \tau\ltwoinner{\vec{V}(p).\nabla(g^{(k+1)}),e^{(k+1)}}\nonumber\\
&=&\tau\ltwoinner{\vec{V}(g^{(k)}) \cdot \nabla\psi_N,e^{(k+1)}}  -  \tau\ltwoinner{\vec{V}(p).\nabla(g^{(k+1)}),e^{(k+1)}}\label{eqConv3M}\\
&\leq & \tau\ltwonorm{\vec{V}(g^{(k)}) \cdot \nabla\psi_N }\ltwonorm{ e^{(k+1)}}  + \tau \ltwonorm{ \vec{V}(p).\nabla(g^{(k+1)})} \ltwonorm{e^{(k+1)}} \nonumber\\
\therefore \ltwonorm{e^{(k+1)}} &\leq & \tau\ltwonorm{\vec{V}(g^{(k)}) \cdot \nabla\psi_N}  + \tau \ltwonorm{ \vec{V}(p).\nabla g^{(k+1)}} \label{eqConv4M}
\end{eqnarray}
Similarly to \eqref{eqCia}, by Ciarlet (\cite{ciarlet}, Theorem 3.2.6), we have \eqref{eq:C1M}, then for $\tau \leq \dfrac{1}{6\loneinftynorm{p} }$ we get \eqref{eq:C2M}
\begin{eqnarray}
\ltwonorm{\vec{V}(g^{(k)}) \cdot \nabla \psi_N} &\leq &c_{inv}^2 h^{-5/2} \ltwonorm{\psi_N} \;.\; \ltwonorm{g^{(k)}} \label{eq:C1M}\\
 &\leq &c_{inv}^2 h^{-5/2} e^{3 T \loneinftynorm{p}} \ltwonorm{w_0} \;.\; \ltwonorm{g^{(k)}}  =h^{-5/2} \tilde{D}_3 \ltwonorm{g^{(k)}}\qquad \mbox{By } \eqref{app1}\quad \label{eq:C2M}\\
 \ltwonorm{\vec{V}(p) \cdot \nabla  g^{(k+1)} } &\leq &  2 \loneinftynorm{p} \;.\;\honenorm{ g^{(k+1)} } \quad \mbox{Similarly to \eqref{eqP}\qquad} \nonumber \\
 &\leq &  2 \loneinftynorm{p} \;.\;\ltwonorm{ e^{(k+1)} }\qquad \mbox{By } \eqref{eq2Conv2}\quad \label{eq:C3M}
\end{eqnarray} 
where $\tilde{D}_3:= \tilde{D}_3(\Omega,p,T,w_0) =c_{inv}^2  e^{3 T \loneinftynorm{p}} \ltwonorm{w_0}$.
Replacing \eqref{eq:C2M} and \eqref{eq:C3M} in \eqref{eqConv4M}, we get \eqref{eq:ConvM}
\begin{eqnarray}
 \ltwonorm{e^{(k+1)}} &\leq &  \dfrac{\tau h^{-5/2} \tilde{D}_3}{1-2 \tau \loneinftynorm{p}} \ltwonorm{g^{(k)}}
 = {c}\ltwonorm{g^{(k)}} \nonumber\\
\implies  \ltwonorm{e^{(k+1)}}^2 &\leq & c^2 \ltwonorm{g^{(k)}}^2 \label{eq:ConvM}\\
\therefore \ltwonorm{e^{(k+1)}}^2 + \ltwonorm{g^{(k+1)}}^2 &\leq & 2\ltwonorm{e^{(k+1)}}^2 \; \leq \;  2c^2 \ltwonorm{g^{(k)}}^2 \;\leq \; c^2\left( \ltwonorm{e^{(k)}}^2 +\ltwonorm{g^{(k)}}^2  \right)\quad \mbox{By}\; \eqref{eq2Conv2M} \qquad \label{eq:ConvMF}
\end{eqnarray}
where ${c} =  \dfrac{\tau h^{-5/2} \tilde{D}_3}{1-2 \tau \loneinftynorm{p}}$. If $\tau < \dfrac{1}{h^{-5/2} \tilde{D}_3 +2  \loneinftynorm{p}} = \dfrac{h^{5/2}}{ \tilde{D}_3 +2 h^{5/2} \loneinftynorm{p}}$, then  $c<1$.
Moreover, assuming $h < 1$, then $\tau< \dfrac{h^{5/2}}{ {D}_3 } <  \dfrac{h^{5/2}}{ \tilde{D}_3 +2 h^{5/2} \loneinftynorm{p}}$ where $ D_3 = \tilde{D}_3 +2 \loneinftynorm{p}$, which ends the proof.
\end{proof}
\noindent A corollary of Theorem \eqref{thrm:ConvMN} is the global convergence of Modified Newton's method.
\begin{theorem}\label{thrm:ConvMN3}
Let $\tau \leq \min\left\{\dfrac{1}{6||p||_{1,\infty}}, \dfrac{h^2}{16{c}_{0,inv}\Mnorm{W(t)}}, \dfrac{h^{5/2}}{D_3}\right\} = O(h^{2.5})$,
 then, for any choice of initial guesses $\{\phi_N^{(0)}, \psi_N^{(0)} \}$, Modified Newton's method converges to the unique solution of \eqref{HMC-Disc-Comp},
$$\lim\limits_{k\rightarrow \infty} \phi_N^{(k)} = \phi_N \qquad and \qquad \lim\limits_{k\rightarrow \infty} \psi_N^{(k)} = \psi_N,$$
where $D_3 := D_3(\Omega,p,T,w_0) = c_{inv}^2 e^{3 T \loneinftynorm{p}} \ltwonorm{w_0} +2 \loneinftynorm{p}$, and $h<1$.
\end{theorem}
\begin{proof}
By theorem \ref{thrm:ConvMN}, there exists $c<1$ such that
\begin{eqnarray}
\ltwonorm{e^{(k+1)}}^2 + \ltwonorm{g^{(k+1)}}^2 &\leq & c^2 \left( \ltwonorm{e^{(k)}}^2 + \ltwonorm{g^{(k)}}^2\right)
\; \leq \; c^{2(k+1)} \left( \ltwonorm{e^{(0)}}^2 + \ltwonorm{g^{(0)}}^2\right)\nonumber \\
\implies \lim\limits_{k\rightarrow \infty} \left(\ltwonorm{e^{(k+1)}}^2 + \ltwonorm{g^{(k+1)}}^2\right) &\leq & \lim\limits_{k\rightarrow \infty}  c^{2(k+1)}  \left(\ltwonorm{e^{(0)}}^2 + \ltwonorm{g^{(0)}}^2\right) \;=\; 0 \nonumber\\
\therefore \lim\limits_{k\rightarrow \infty} \ltwonorm{e^{(k+1)}}^2 &=& 0 \qquad and \qquad  \lim\limits_{k\rightarrow \infty} \ltwonorm{g^{(k+1)}}^2 \;=\; 0\nonumber
\end{eqnarray}
$\therefore \lim\limits_{k\rightarrow \infty} \psi_N^{(k+1)} \;=\;  \psi_N \qquad and \qquad   \lim\limits_{k\rightarrow \infty} \phi_N^{(k+1)} \;=\;  \phi_N$.
\end{proof}
\section{Computer Simulations and Testings} \label{Sec:test}
We implement the three discussed methods, Newton (Algorithm \ref{alg:HMC-Newton0}), Chord (Algorithm \ref{alg:HMC-Chord}), Modified Newton (Algorithm \ref{alg:HMC-Modified Newton}) for $k_{max} = 20$ and $tol = 10^{-6}$ using Freefem++ \cite{MR3043640}, a programming language and software focused on solving partial differential equations using the finite element method. 
\begin{algorithm}[h!]
\centering 
\caption{  Solving HM using Newton's Method }
{\renewcommand{\arraystretch}{1.3}
\begin{algorithmic}[1]
\Statex{\textbf{Input:} \;\;\;$A$: stiffness  matrix ; $M$: Mass matrix; $K=M+A$;  $S(U)$; $R$: As defined in \cite{FEHMX}; $B(W)$: Matrix \eqref{BWdef};}
\Statex{ \qquad \qquad $U_0;W_0$: the discrete initial condition vectors; $T$: end time; $\tau$: time step; $N$:total number of mesh nodes;}
 \Statex{\qquad\quad \;\;\;$k_{max}$: maximum Newton iterations; $tol$: Newton's relative error stopping  tolerance.  }
 \Statex{\textbf{Output:} $U$: $N\times (T/\tau +1)$ matrix containing the computed vectors $U_t$ for $t=0,\tau, 2\tau,\cdots,T$ \vspace{3mm}}
\For {$t= 0:\tau:T$\vspace{2mm}} 
\State $U_{t,0}=U_{t}$;\;\; $W_{t,0}=W_{t}$;\;\;$error = 1$; \;\; $k=0$\vspace{1mm}
\State Let $Z=M*W_{t};$ \vspace{1mm}
\While {( $error > tol$ and $k < k_{max}$ )\vspace{2mm}}
\State Let $ g=\tau S(U_{t,k})*W_{t,k}+Z$;\vspace{1mm}
\State Let $r(0:N-1)= g$ and $  r(n:2*N-1)=0$;\vspace{1mm}
\State Let $J_1=\tau B(W_{t,k})-\tau R$, and $J2= \tau S(U_{t,k})+M$; \vspace{1mm}
\State Construct $J: \quad  J=[[J1,J2],[K,-M]]$;\vspace{1mm}
\State Solve for $V$:\quad $J*V = r$;\vspace{1mm}
\State Let $U_{t,k+1}=V(0:N-1)$ and $W_{t,k+1}=V(N:2*N-1)$;\vspace{1mm}
\State $ error = \dfrac{||U_{t,k+1}-U_{t,k}||}{||U_{t,k}||}$;\;\; $k=k+1$;\vspace{1mm}
\EndWhile \vspace{1mm}
\State $U_{t+1} = U_{t,k}$;
\State  $W_{t+1} = W_{t,k}$;
\EndFor
\end{algorithmic}}
\label{alg:HMC-Newton0}
\end{algorithm}

\begin{algorithm}[h!]
\centering
\caption{  Solving HM using Chord's Method }
{\renewcommand{\arraystretch}{1.3}
\begin{algorithmic}[1]
\Statex{\textbf{Input:} \;\;\;$A$: stiffness  matrix ; $M$: Mass matrix; $K=M+A$;  $S(U)$; $R$: As defined in \cite{FEHMX}; $B(W)$: Matrix \eqref{BWdef};}
\Statex{ \qquad \qquad $U_0;W_0$: the discrete initial condition vectors; $T$: end time; $\tau$: time step; $N$:total number of mesh nodes;}
 \Statex{\qquad\quad \;\;\;$k_{max}$: maximum Chord iterations; $tol$: Chord's  relative error stopping tolerance. }
 \Statex{\textbf{Output:} $U$: $N\times (T/\tau +1)$ matrix containing the computed vectors $U_t$ for $t=0,\tau, 2\tau,\cdots,T$ \vspace{3mm}}
\For {$t= 0:\tau:T$\vspace{2mm}} 
\State $U_{t,0}=U_{t}$;\;\;$W_{t,0}=W_{t}$;\;\; $error = 1$; \;\; $k=0$\vspace{1mm}
\State Let $Z=M*W_{t};$ Let $J_1=\tau B(W_{t,0})-\tau R;$  Let $J2= \tau S(U_{t,0})+M$; \vspace{1mm}
\State Construct $J: \quad  J=[[J1,J2],[K,-M]]$;\vspace{1mm}
\While {( $error > tol$ and $k < k_{max}$ )\vspace{2mm}}
\State Let $ g= \tau S(U_{t,k})*(W_{t,0} - W_{t,k}) + \tau S(U_{t,0})*W_{t,k} + Z$;\vspace{1mm}
\State Let $r(0:N-1)= g$ and $  r(n:2*N-1)=0$;\vspace{1mm}
\State Solve for $V$:\quad $J*V = r$;\vspace{1mm}
\State let $U_{t,k+1}=V(0:N-1)$ and $W_{t,k+1}=V(n:2*N-1)$;\vspace{1mm}
\State $ error = \dfrac{||U_{t,k+1}-U_{t,k}||}{||U_{t,k}||}$;\;\; $k=k+1$;\vspace{1mm}
\EndWhile \vspace{1mm}
\State $U_{t+1} = U_{t,k}$;\;\;\;  $W_{t+1} = W_{t,k}$;
\EndFor
\end{algorithmic}}
\label{alg:HMC-Chord}
\end{algorithm}

\begin{algorithm}[H]
\centering
\caption{Solving HM using Modified Newton Method }
{\renewcommand{\arraystretch}{1.3}
\begin{algorithmic}[1]
\Statex{\textbf{Input:} \;\;\;$A$: stiffness  matrix ; $M$: Mass matrix; $K=M+A$;  $S(U)$; $R$: As defined in \cite{FEHMX}; }
\Statex{ \qquad \qquad $U_0;W_0$: the discrete initial condition vectors; $T$: end time; $\tau$: time step; $N$:total number of mesh nodes;}
 \Statex{\qquad\quad \;\;\;$k_{max}$: maximum Modified Newton iterations; $tol$:  Modified Newton's relative error stopping tolerance.   }
 \Statex{\textbf{Output:} $U$: $N\times (T/\tau +1)$ matrix containing the computed vectors $U_t$ for $t=0,\tau, 2\tau,\cdots,T$ \vspace{3mm}}
\For {$t= 0:\tau:T$\vspace{2mm}} 
\State $U_{t,0}=U_{t}$;\; $W_{t,0}=W_{t}$;\; $error = 1$; \; $k=0$\vspace{1mm}
\State Let $Z=M*W_{t};$ Let $r(0:N-1)= Z;$ and $  r(N:2*N-1)=0$;\vspace{1mm}
\While {( $error > tol$ and $k < k_{max}$ )\vspace{2mm}}
\State let $J_2=M+\tau S(U_{t,k})$;\vspace{1mm}
\State Construct $J: \quad  J=[[-\tau R,J2],[K,-M]]$;\vspace{1mm}
\State Solve for $V$:\quad $J*V = r$;\vspace{1mm}
\State let $U_{t,k+1}=V(0:N-1)$ and $W_{t,k+1}=V(n:2*N-1)$;\vspace{1mm}
\State $ error = \dfrac{||U_{t,k+1}-U_{t,k}||}{||U_{t,k}||}$;\;\; $k=k+1$;\vspace{1mm}
\EndWhile \vspace{1mm}
\State $U_{t+1} = U_{t,k}$;\;\;\;  $W_{t+1} = W_{t,k}$;
\EndFor
\end{algorithmic}}
\label{alg:HMC-Modified Newton}
\end{algorithm}

We consider the same initial conditions $u_0$ as in \cite{FEHM} (Table \ref{tab:test}), and compare the obtained solutions of the three method, the required runtime, and the number of iterations per time step (Table \ref{tab:time}). Moreover, we compare them with the semi-linear algorithm introduced in \cite{FEHM}.

\begin{table}[H]
\centering
\setlength{\tabcolsep}{5pt}
\renewcommand{\arraystretch}{2}
\begin{tabular}{|c|c|c|c|}
\hline
&$\Omega$&$u_0(x,y)$&$p(x,y) = \ln\dfrac{n_0}{w_{ci}}$\\
\hline
Test 1& $[0,1]\times[0,1]$&$10^{-5}\sin(10\pi y) $&$12x$\\ \hline
Test 2& $[0,\pi]\times[0,\pi]$& $10^{-5}\sin(3y) $&$12x$ \\ \hline
Test 3& $[0,\pi]\times[0,\pi]$& $10^{-5}\sin(3x) $&$12x$ \\ \hline
Test 4& $[0,\pi]\times[0,\pi]$&$10^{-10}xy(x-2)\sin(x)$&$12x$\\  \hline
Test 5&$[0,20]\times[0,20]$& $-10^{-5}(x -10)e^{-0.5(x-10)^2-0.5(y-10)^2} $&$\ln(10^{13}e^{-(x-10)^2/64-(y-10)^2/64})$\\
\hline
\end{tabular}
\caption{\small Considered Test Cases} \label{tab:test}
\end{table}

 Similarly to \cite{FEHM}, we consider a square domain $[x_0,x_n]\times [y_0,y_n]$ with a uniform mesh in the $x$ and $y$ direction ($x_i - x_{i-1} = \dfrac{x_n - x_0}{n} = y_i - y_{i-1} = \dfrac{y_n - y_0}{n}$ for $i = 1,2,..,n$ and $n$ intervals in the $x$ and $y$ directions respectively) and the finite element $\mathbb{P}1$ space with periodic boundary conditions, using appropriate Freefem++ functions. Even though theoretically $\tau = O(h^{2.5})$, we use $\tau = O(h)$ in our testings, specifically $\tau = 0.1$. Moreover, the simulation is stopped once the maximum value of $u(t)$ at one of the mesh nodes is $0.3$, which corresponds to the maximum value attained physically.

The function $p$ from the initial Hasegawa-Mima PDE and the initial condition $u_0$ are  given as input.
 As for the initial condition $w_0 = u_0 - \Delta u_0$ it could be given as input if $u_0$ is a simple function. However, for any function $u_0$, we compute the vector $W_0 = W(0)$ by solving the linear system $${M}* W_0 = {K}*U_0$$ where the vector $U_0 = U(0)$.
    The matrices ${M}$, ${K}$, ${R}$ and ${S}({U}^k)$ are generated in Freefem++ using their corresponding variational formulations. As for $B(W^k)$, it requires the generation of $N = (n-1)^2$ matrices ${S}({e}^j)$ of size $N\times N$, where each is multiplied by $W^k$. The $N$ matrices ${S}({e}^j)$ could be generated once and stored. However, this would require $N^3$ words, which is memory-bound  (for $n=65$, $N^3 = 6.87\times10^{10}$). Thus, the matrices ${S}({e}^j)$ are regenerated at every iteration to compute $B(W^k)$. 
    
The algorithmic difference between the three methods is that in Newton $B(W^k)$ is computed at every iteration, in Chord it is computed once per time-step, and in Modified Newton it is not computed at all.   
The effect of the $B(W^k)$ computation on the runtimes of the three methods is evident in Table \ref{tab:time}.
 We compare the three methods for $16$ partition intervals in the $x$ and $y$ directions, since for finer meshes the matrices will be larger, particularly the $B(W^k)$ matrix. We also consider end time $T = 10$ with time-step $\tau = 0.1$. 
 
 \begin{table}[H]
\centering
\renewcommand{\arraystretch}{2}
\begin{tabular}{|c|c|c|c|c|c|c|c|c|c|}
\hline
&\multicolumn{3}{|c|}{Newton (Algorithm \ref{alg:HMC-Newton0})}& \multicolumn{3}{|c|}{Chord (Algorithm \ref{alg:HMC-Chord})}&\multicolumn{3}{|c|}{Modified Newton (Algorithm \ref{alg:HMC-Modified Newton})}\\\hline
&Iter&RelErr&Time(s)&Iter&RelErr&Time(s)&Iter&RelErr&Time(s)\\\hline
Test 1&2&	$10^{-14}$	&304.577&2	&$10^{-10}$	&146.755&2	&$10^{-14}$	&1.35253\\\hline
Test 2&2	&$10^{-16}$&	281.969&2&$10^{-12}$&138.803&2	&$10^{-16}$	&1.37822\\\hline
Test 3&1&	$10^{-15}	$&141.301&1&	$10^{-15}	$&	140.674&1&	$10^{-16}$	& 0.763779\\\hline
Test 4&2	&$10^{-14}$&	274.895&2&	$10^{-12}$	&139&2&	$10^{-12}$	&1.32606\\\hline
Test 5&2&	$10^{-14}$	&277.537&2	&$10^{-10}$&	138.5&2&	$10^{-10}$	&1.38517\\\hline
\end{tabular}
\caption{\small Comparison of the 3 methods for the 5 test cases, with $T = 10$, $n=17$, $\tau = 0.1$ showing the number of method iterations per time step (Iter), the last relative error per time step (RelErr) and the total runtime of the algorithm (Time). }
 \label{tab:time}
\end{table}
For all the tests except the third, two method iterations are performed per time-step, where Chord's method is $2$ times faster than Newton's method, since half the  $B(W^k)$ matrices are computed. Moreover, Modified Newton's method is $200$ times faster than Newton's method since it avoids computing $200$  $B(W^k)$ matrices. In the second iteration, the relative error $\dfrac{||U_{t,k+1}-U_{t,k}||}{||U_{t,k}||}$ varies between $10^{-10}$ and $10^{-16}$, which is much smaller than the tolerance $tol = 10^{-6}$. This implies that the change in the solution in the second iteration is relatively negligible. 
 
 For Test 3, one method iteration is performed per time-step. Thus Chord's method and Newton's method require the same runtime, since the same number of $B(W^k)$ matrices are computed. Moreover, Modified Newton's method is $100$ times faster than Newton's method since it avoids computing $100$  $B(W^k)$ matrices.

\begin{figure}[H]
\centering
\includegraphics[scale=0.35]{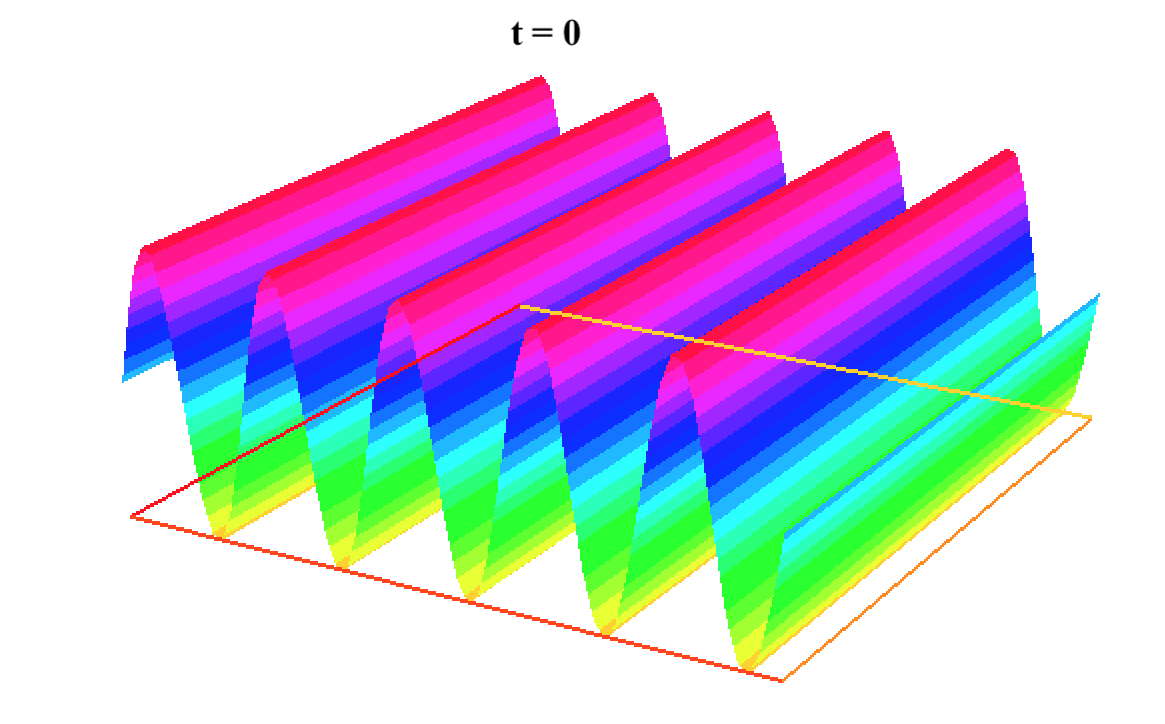}\qquad
\includegraphics[scale=0.35]{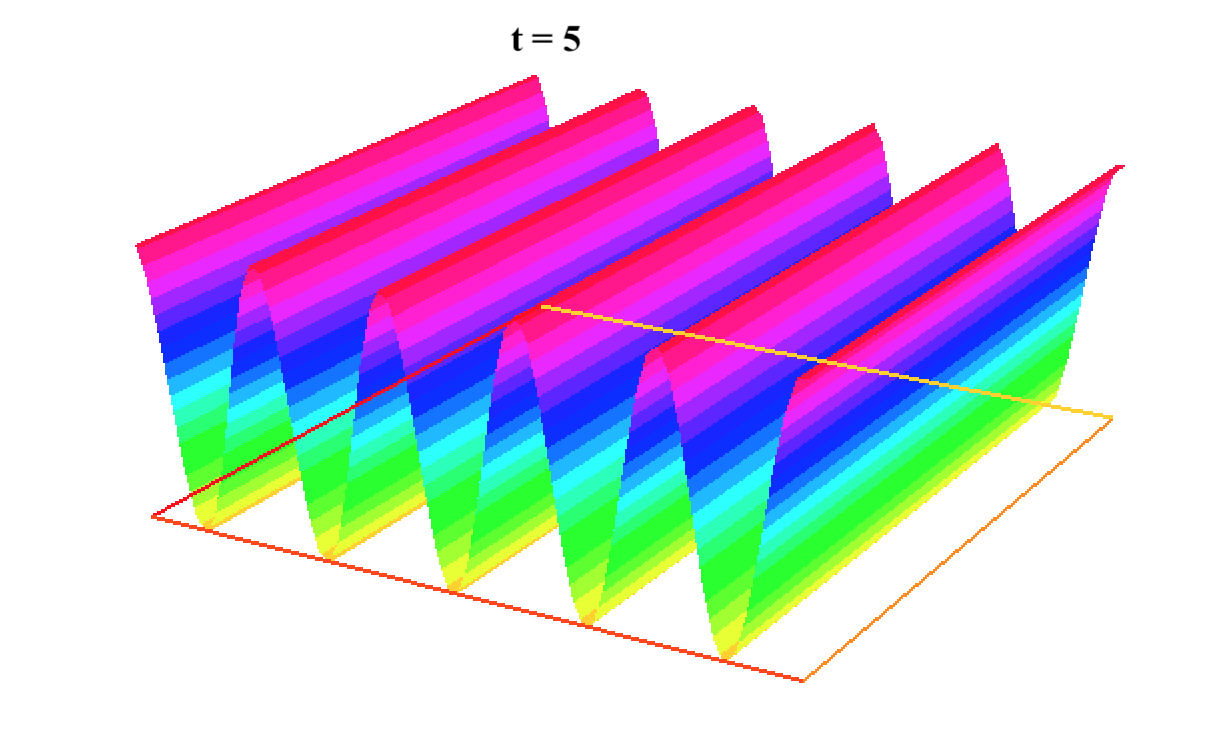}
\includegraphics[scale=0.35]{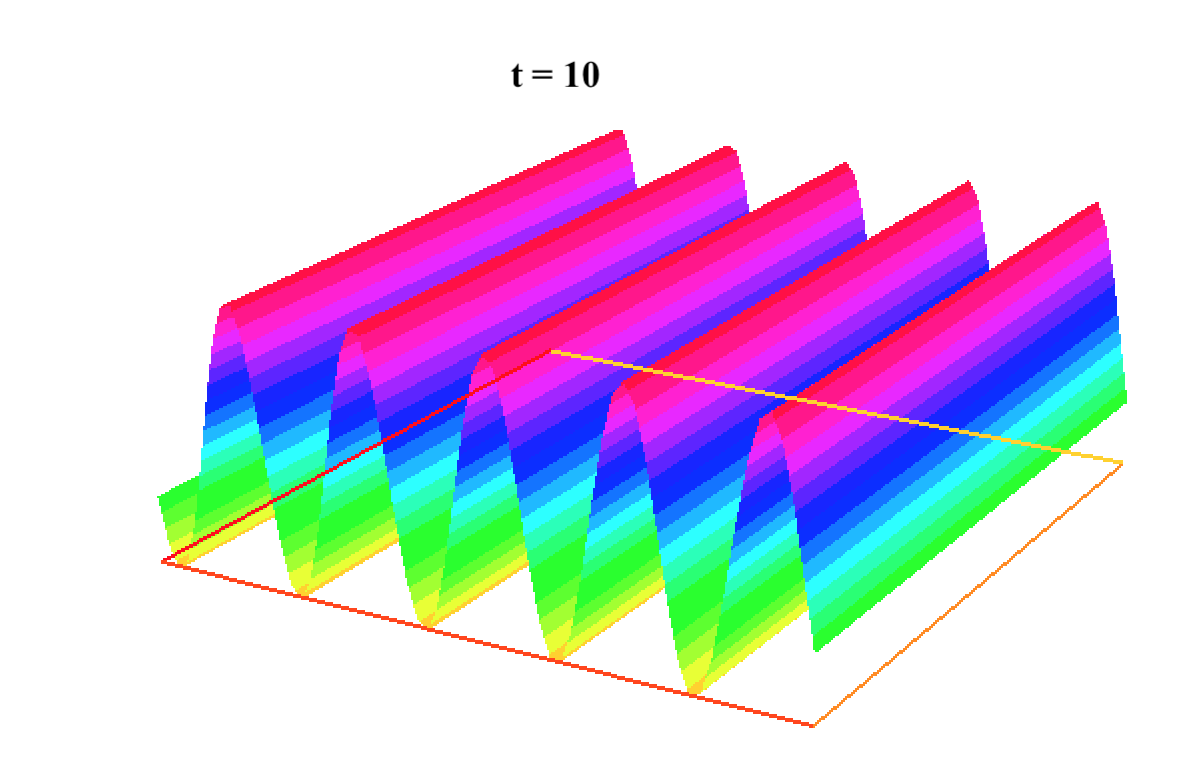}\qquad
\hspace{-2mm} \includegraphics[scale=0.35]{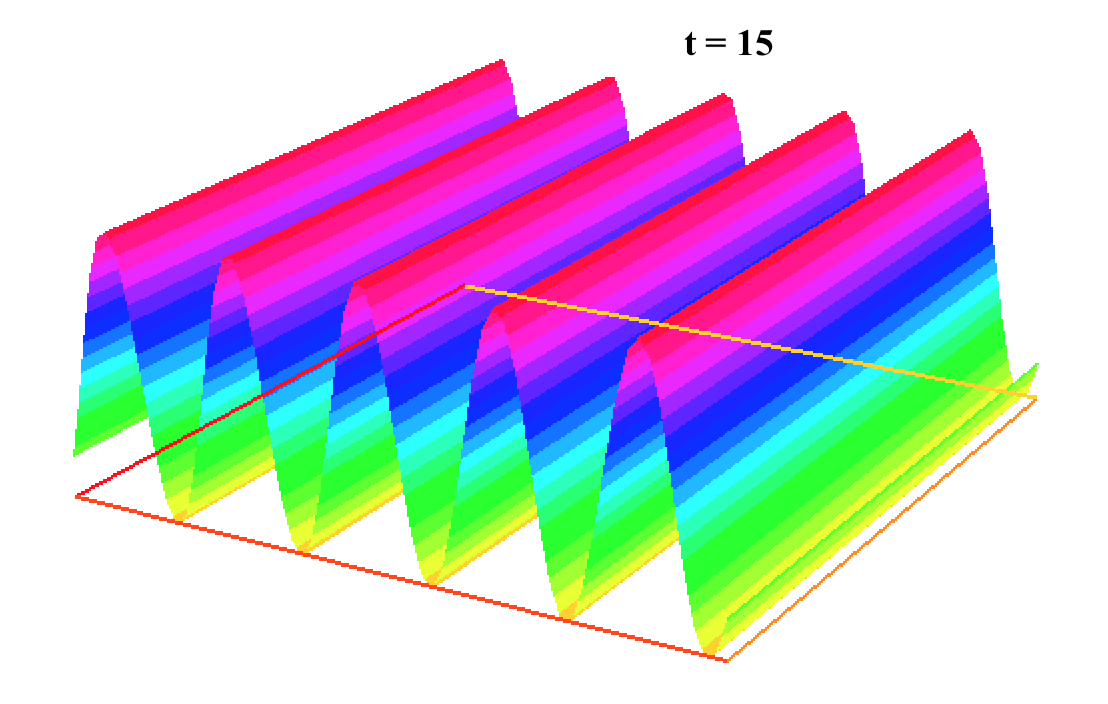}
\includegraphics[scale=0.35]{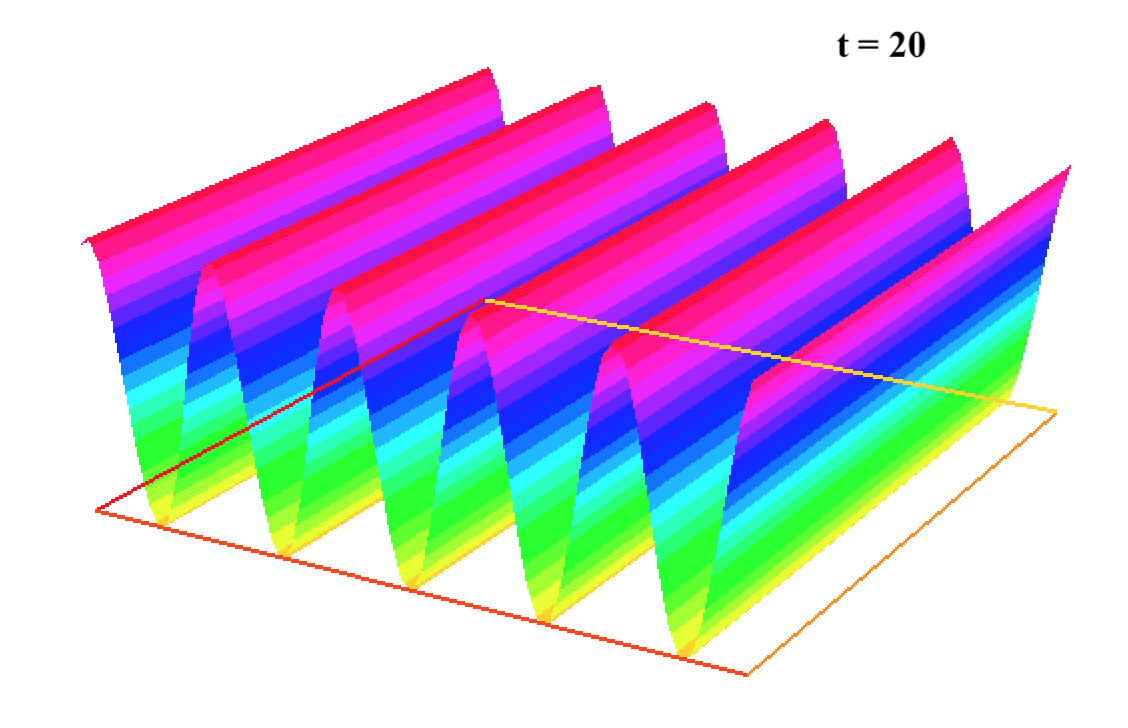} \qquad \qquad\qquad   \qquad \qquad
\includegraphics[scale=0.4]{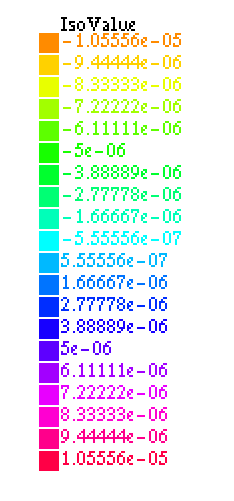}
\caption{\it \small Time evolution of solution $u$ of $\eqref{HMC-Disc-Comp}$ for Test1 using Algorithm \ref{alg:HMC-Modified Newton}, with $\tau = 0.1$, and a $65 \times 65$ grid on  $\Omega = [0,1] \times [0,1]$.}\label{fig:sin9piy}
\end{figure}

\begin{figure}[H]
\centering
\includegraphics[scale=0.25]{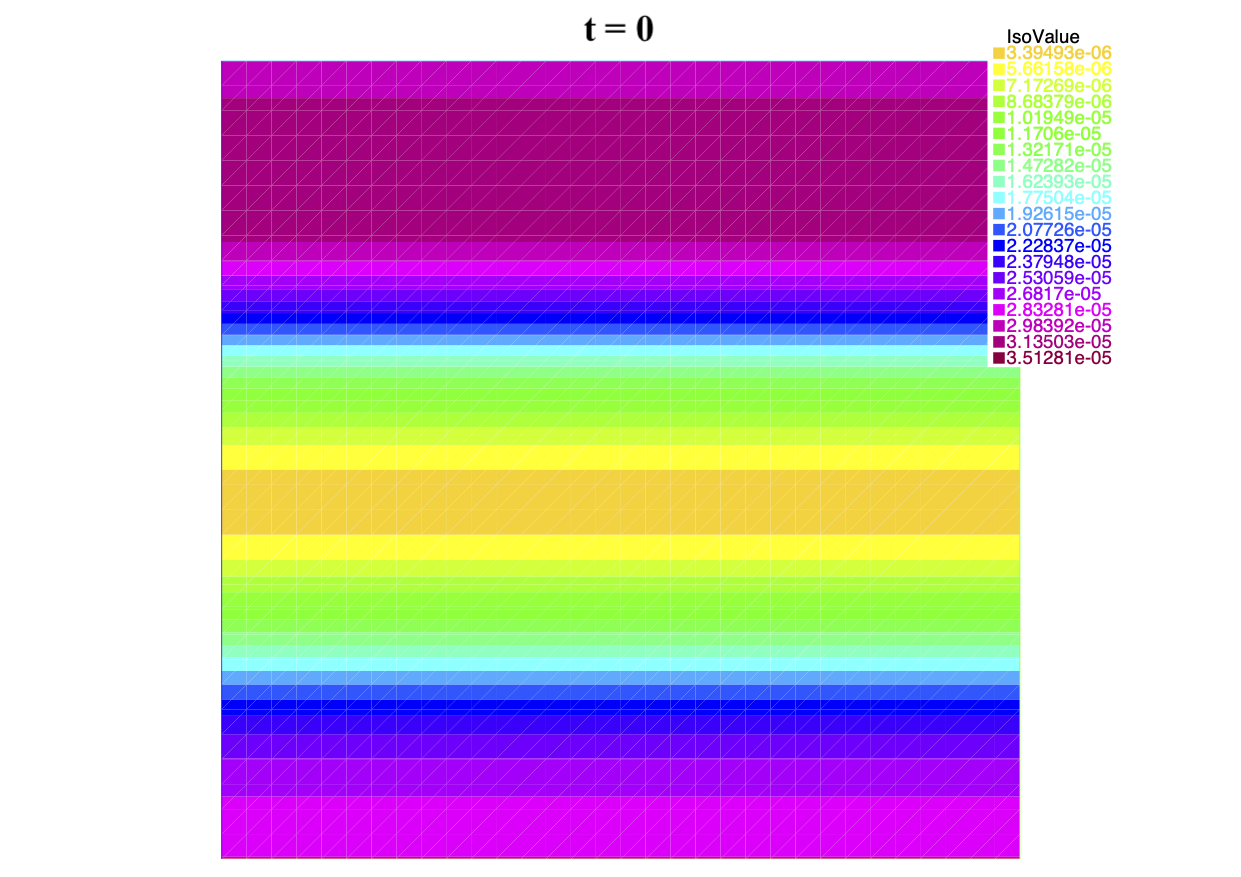}
\includegraphics[scale=0.25]{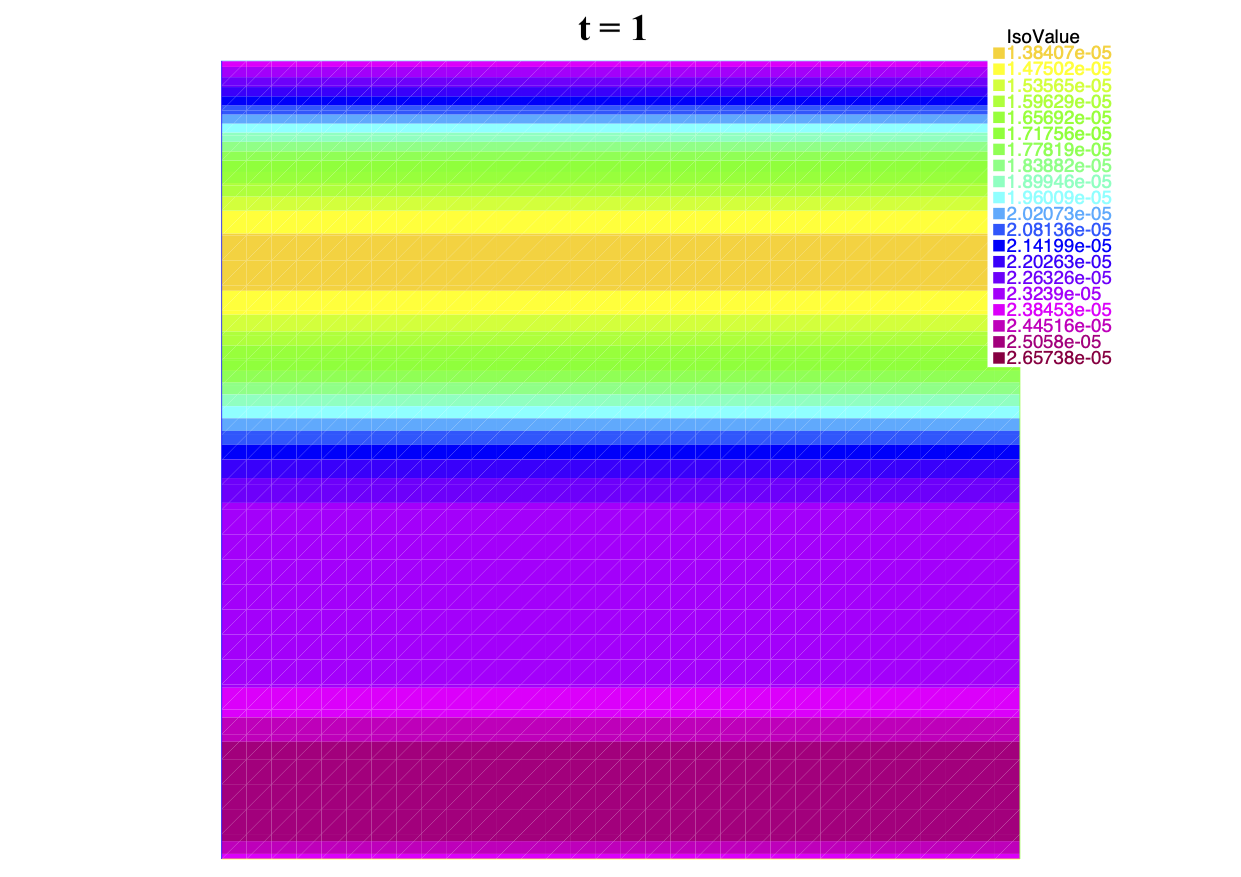}
\includegraphics[scale=0.25]{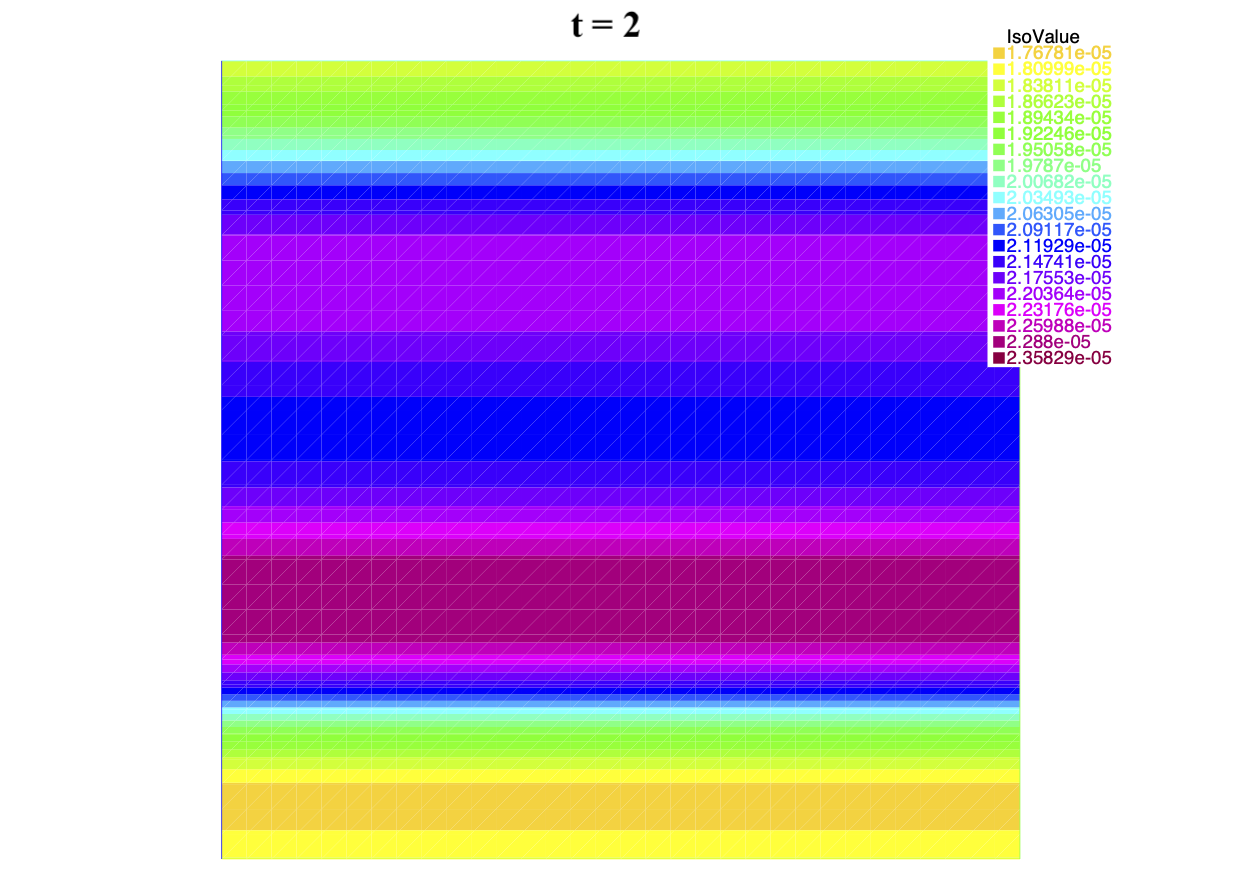}\\ \vspace{2mm}
\includegraphics[scale=0.25]{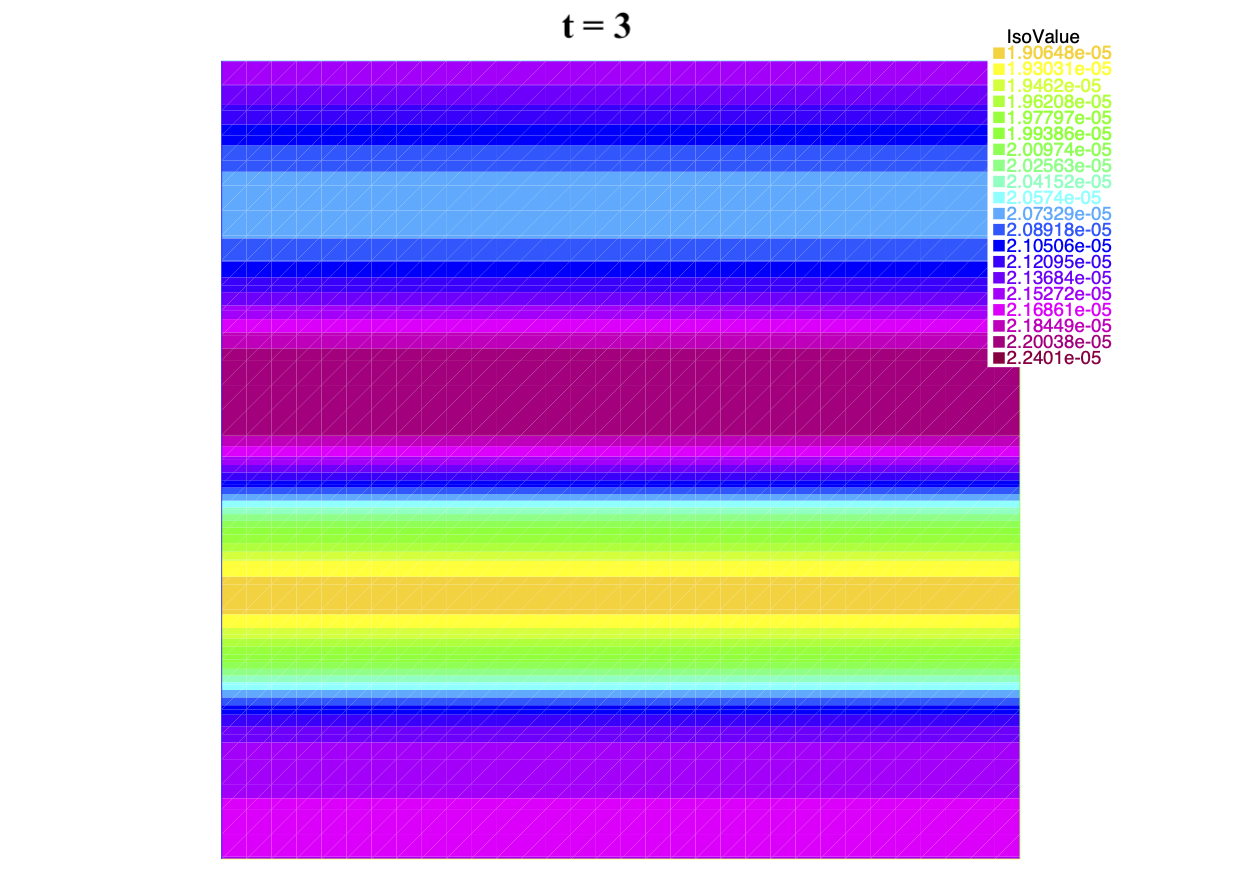}
\includegraphics[scale=0.25]{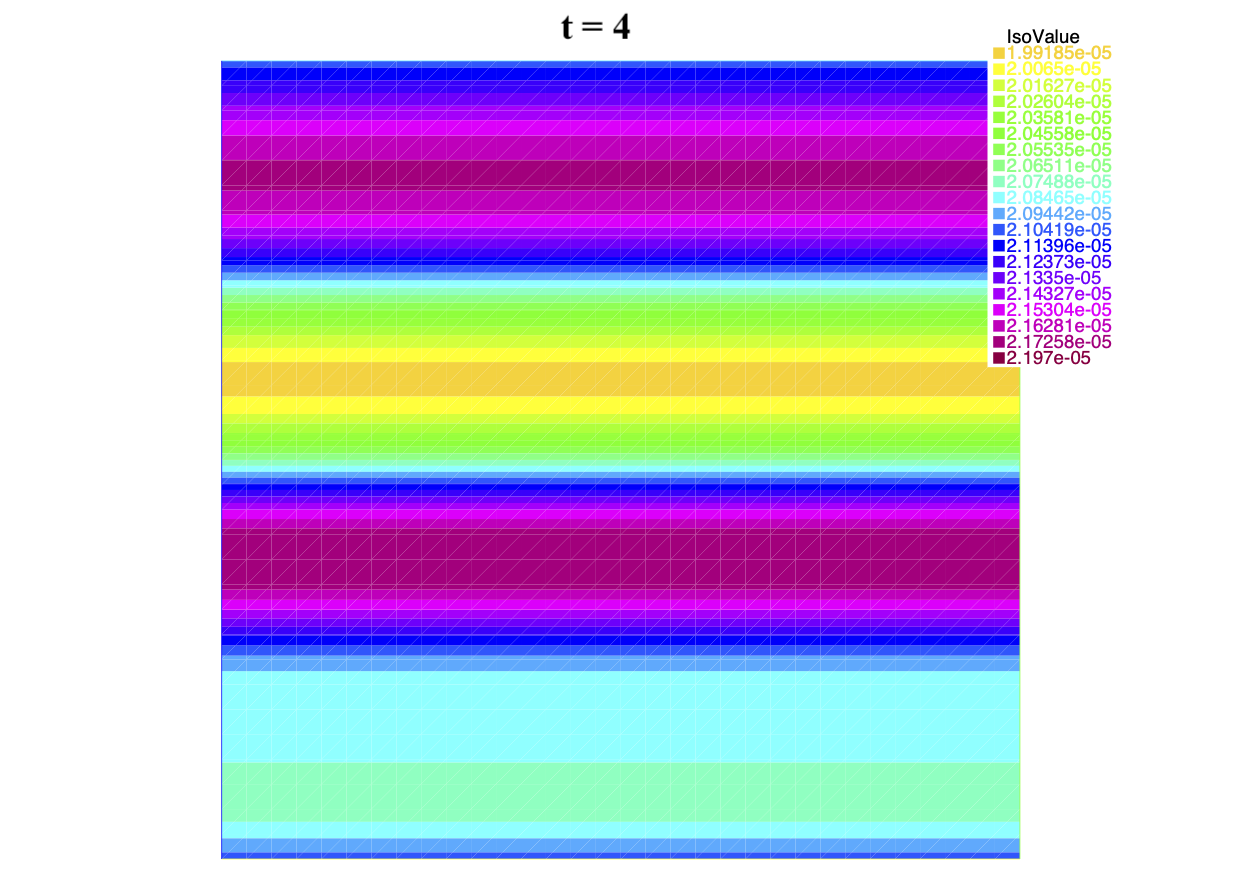}
\includegraphics[scale=0.25]{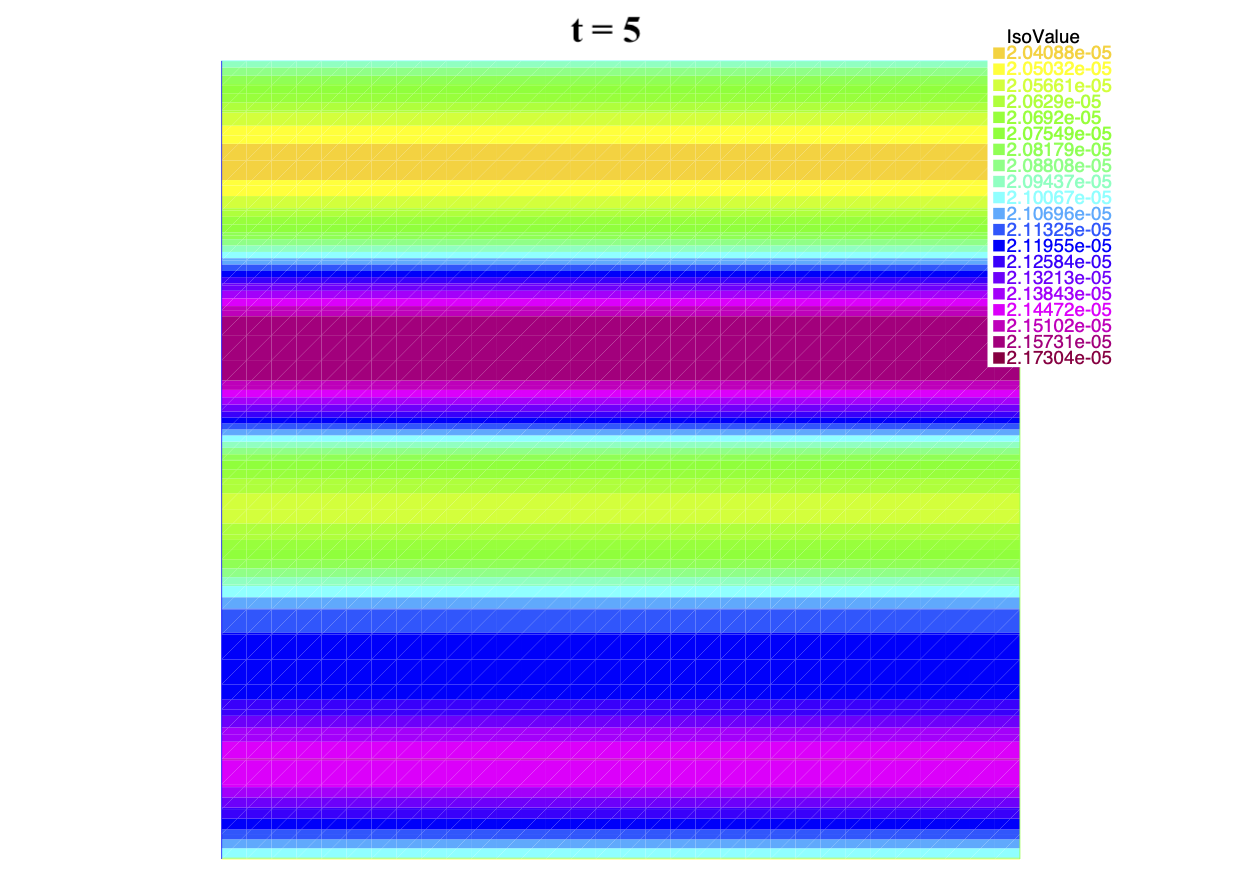}\\ \vspace{2mm}
\includegraphics[scale=0.25]{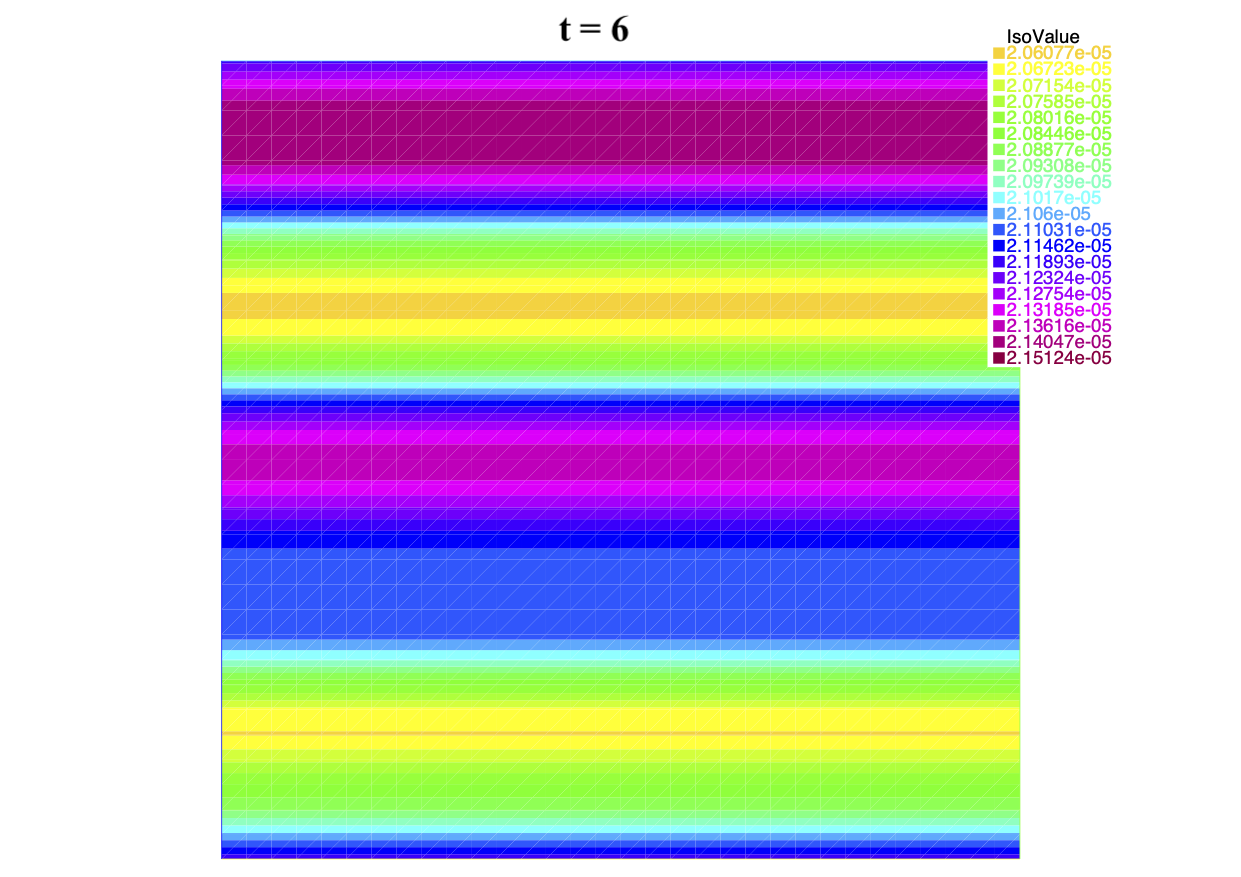}
\includegraphics[scale=0.25]{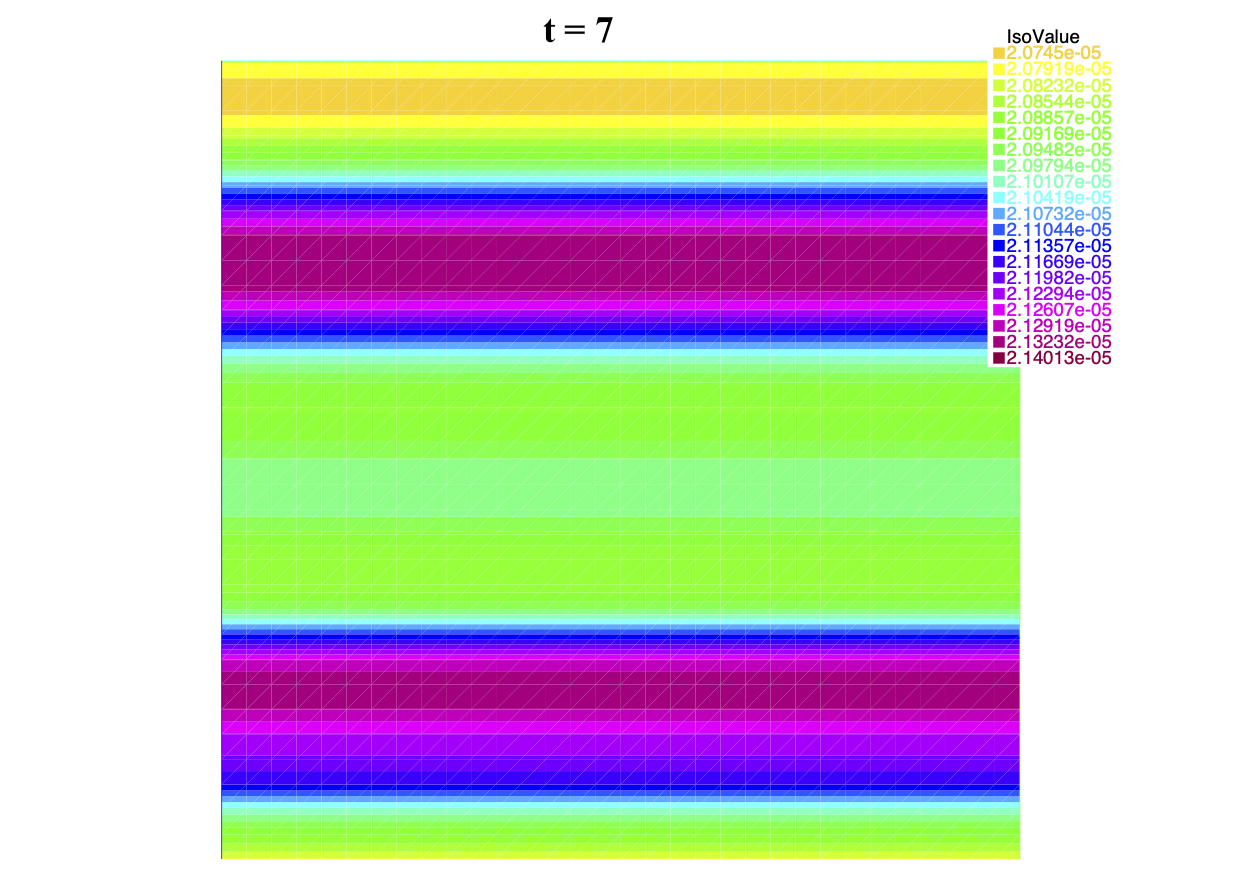}
\includegraphics[scale=0.25]{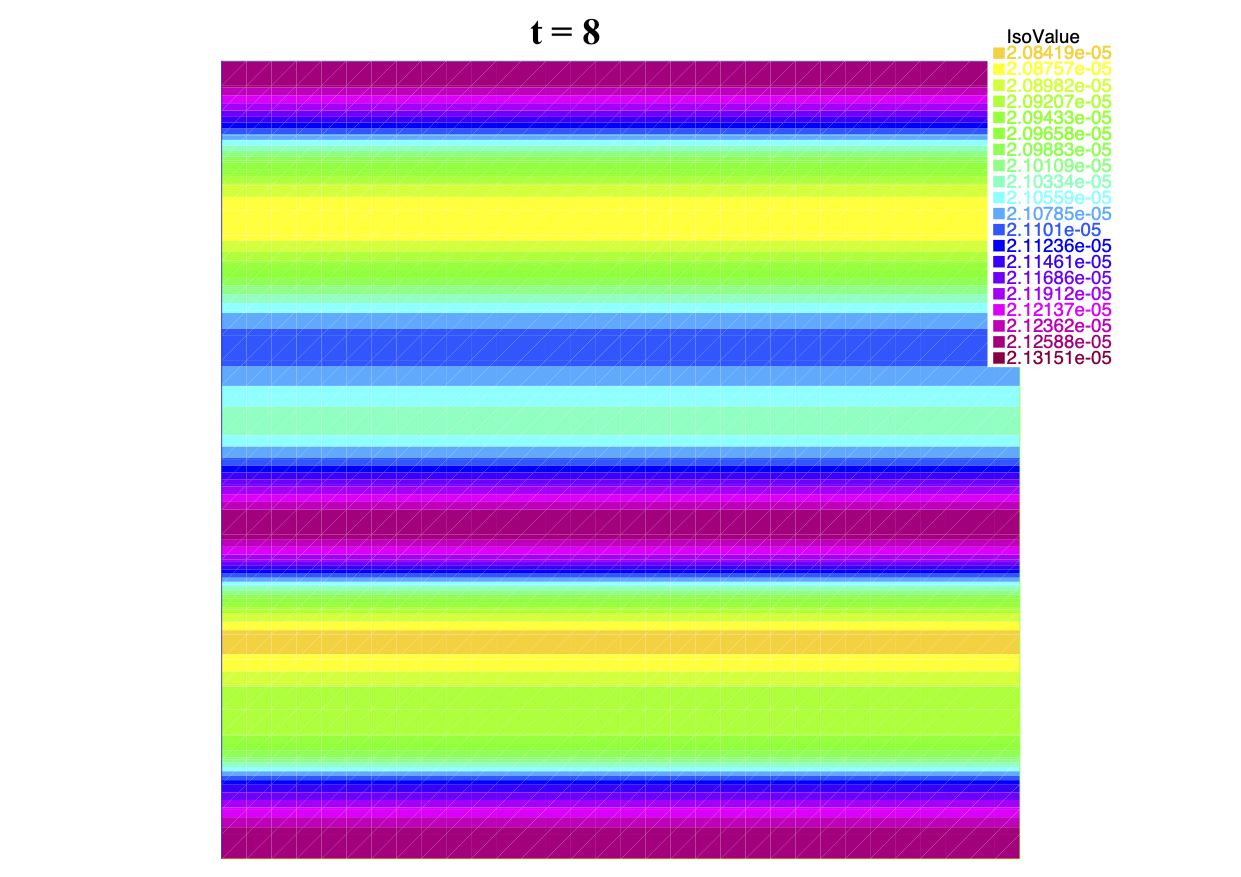}\\ \vspace{2mm}
\caption{\it \small Time evolution of solution $u$ of $\eqref{HMC-Disc-Comp}$ for Test 2 using Algorithm \ref{alg:HMC-Modified Newton}, with $\tau = 0.1$, and a $33 \times 33$ grid on  $\Omega = [0,\pi] \times [0,\pi]$.}\label{fig:sin3y32}
\end{figure}

\begin{figure}[H]
\centering
\includegraphics[scale=0.25]{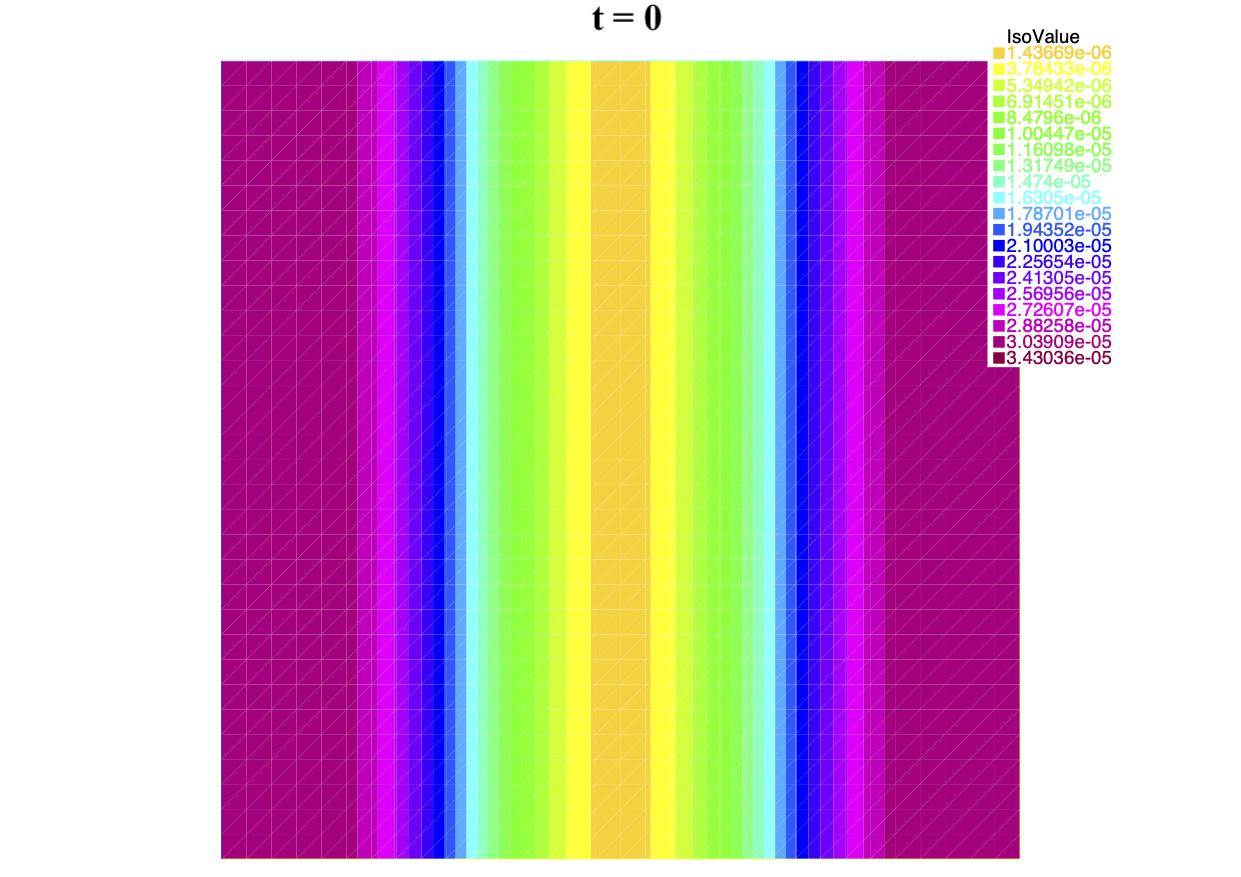}
\includegraphics[scale=0.25]{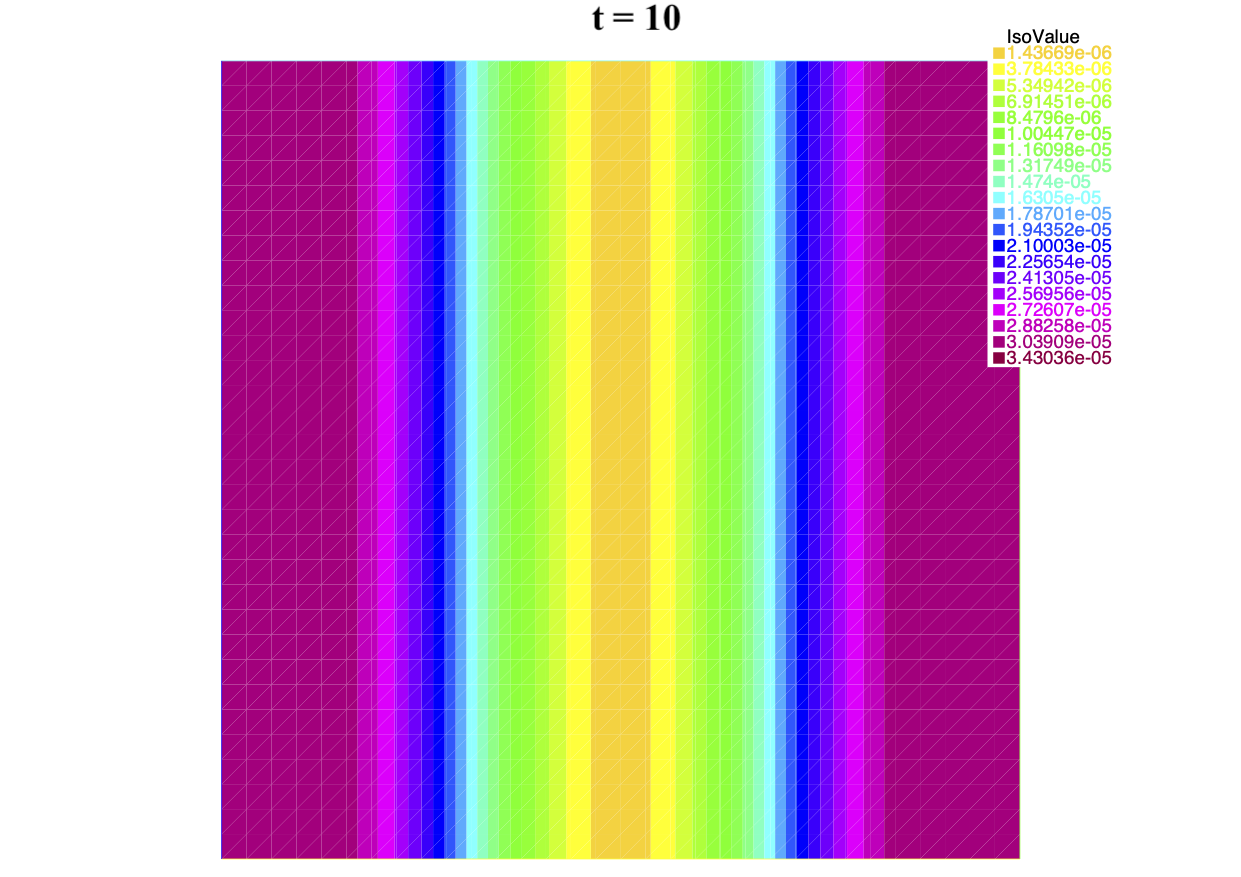}
\includegraphics[scale=0.25]{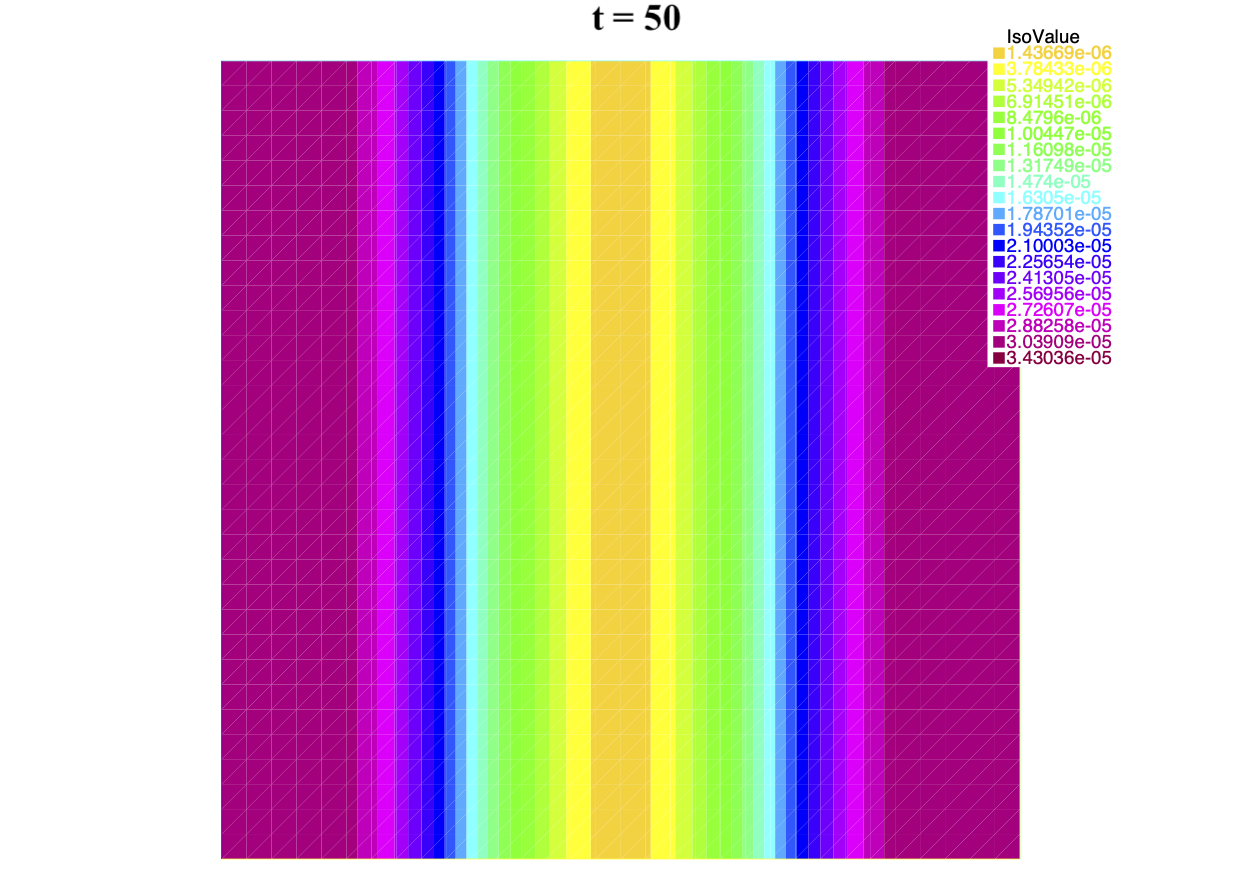}
\caption{\it \small Time evolution of solution $u$ of $\eqref{HMC-Disc-Comp}$ for Test 3 using Algorithm \ref{alg:HMC-Modified Newton}, with $\tau = 0.1$, and a $33 \times 33$ grid on  $\Omega = [0,\pi] \times [0,\pi]$.}\label{fig:sin3x-32}
\end{figure}

\begin{figure}[H]
\centering
\includegraphics[scale=0.25]{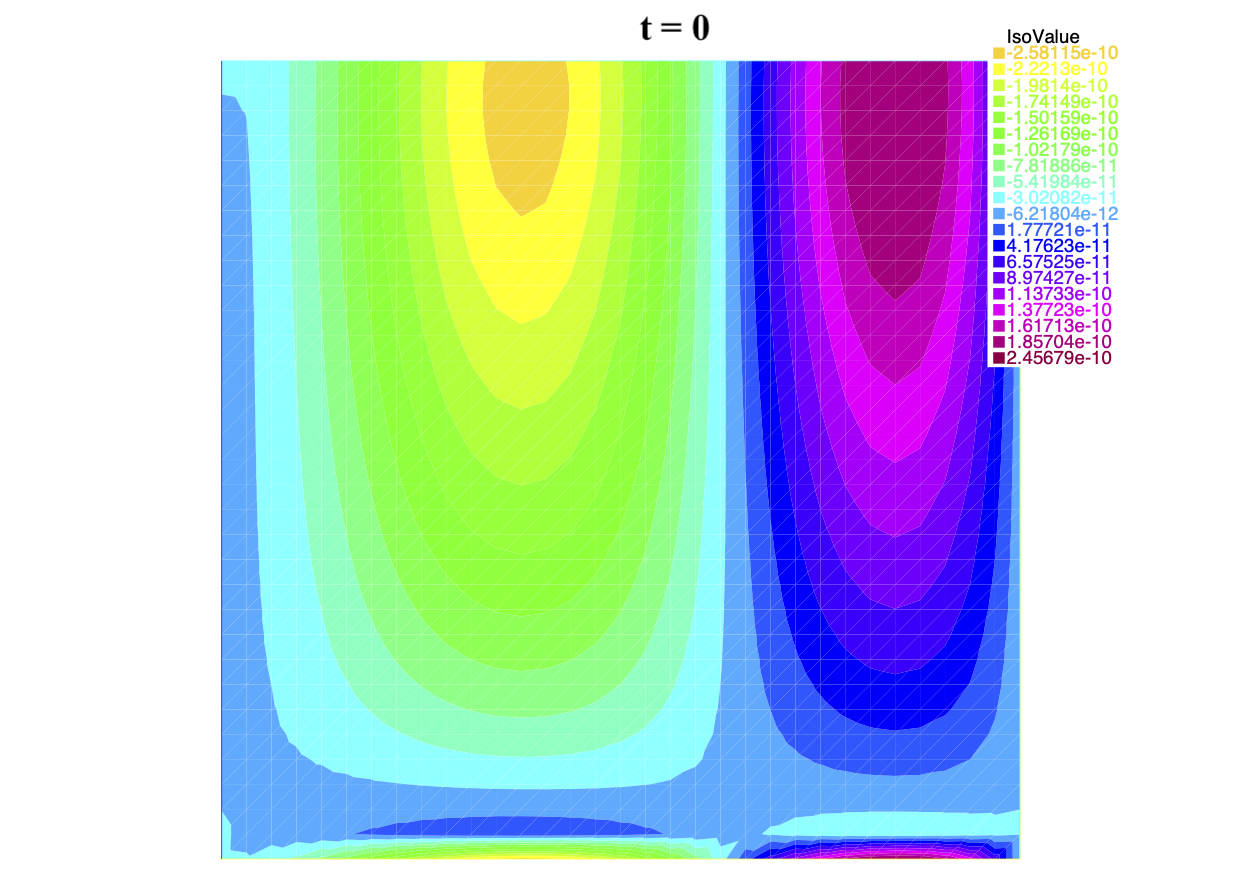}
\includegraphics[scale=0.25]{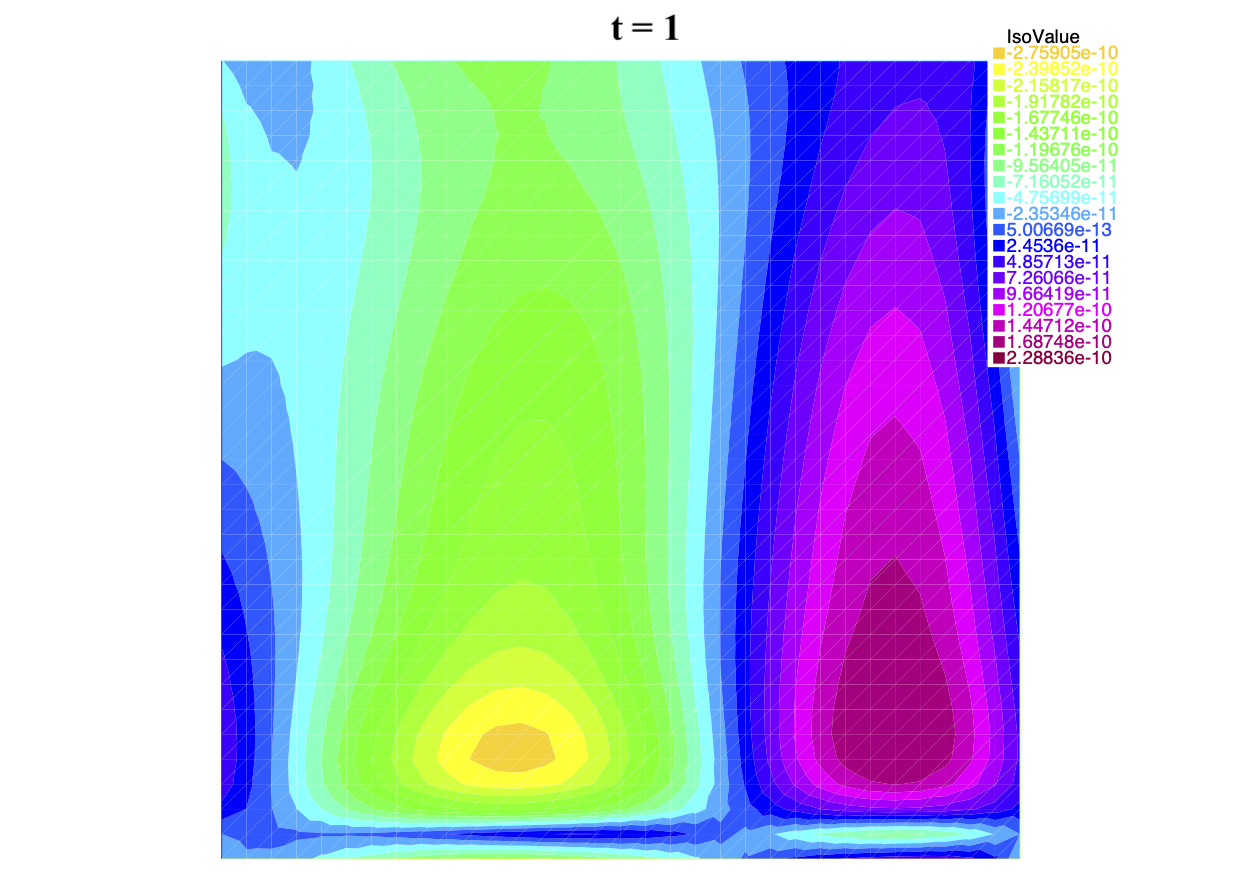}
\includegraphics[scale=0.25]{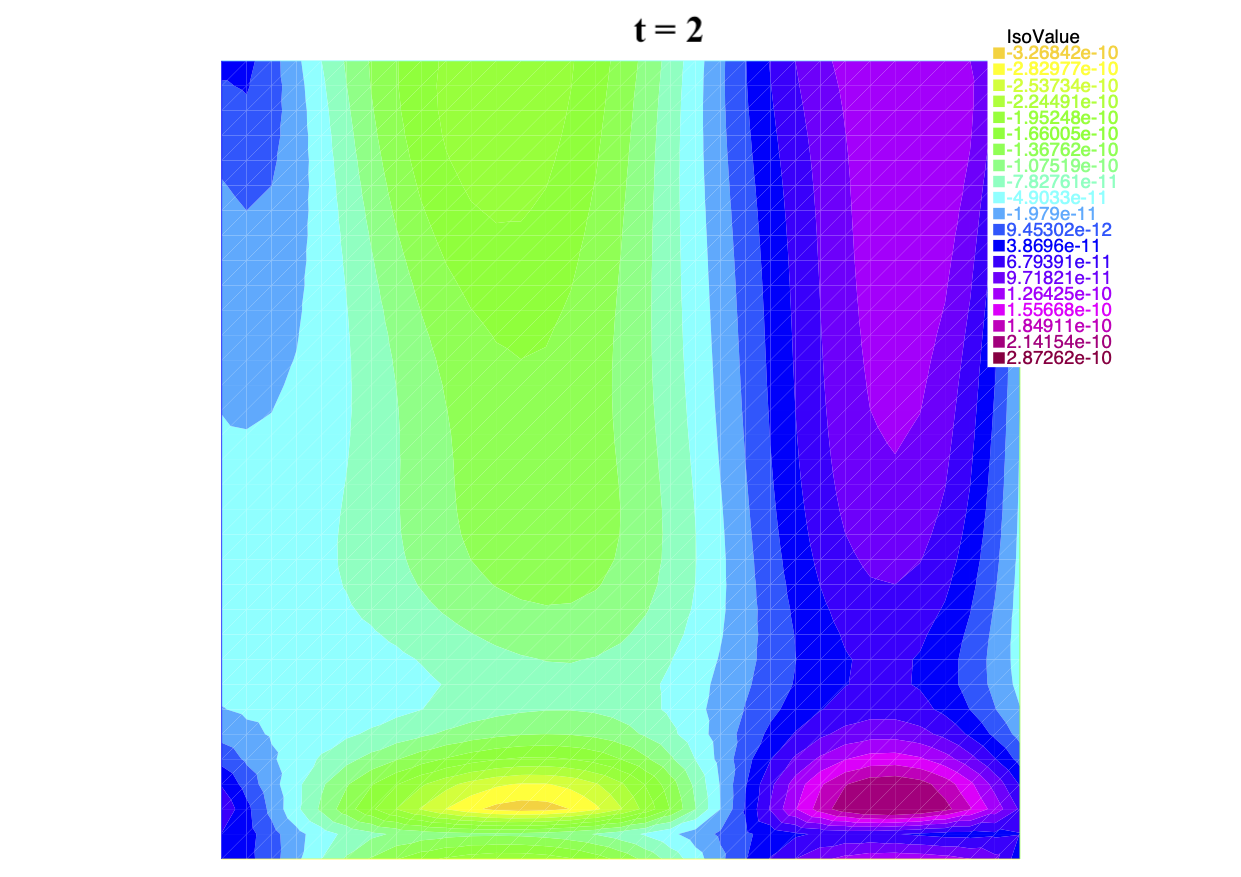}
\includegraphics[scale=0.25]{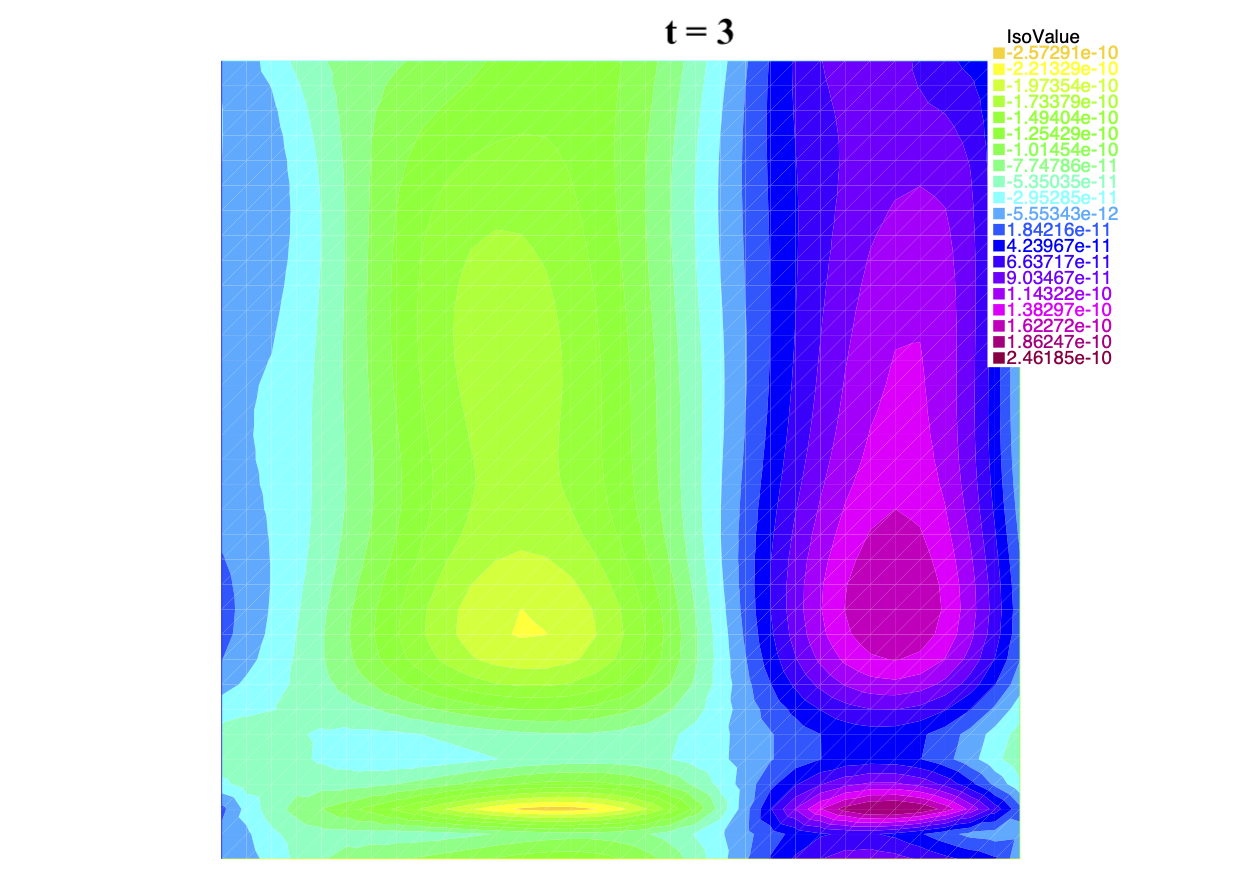}
\includegraphics[scale=0.25]{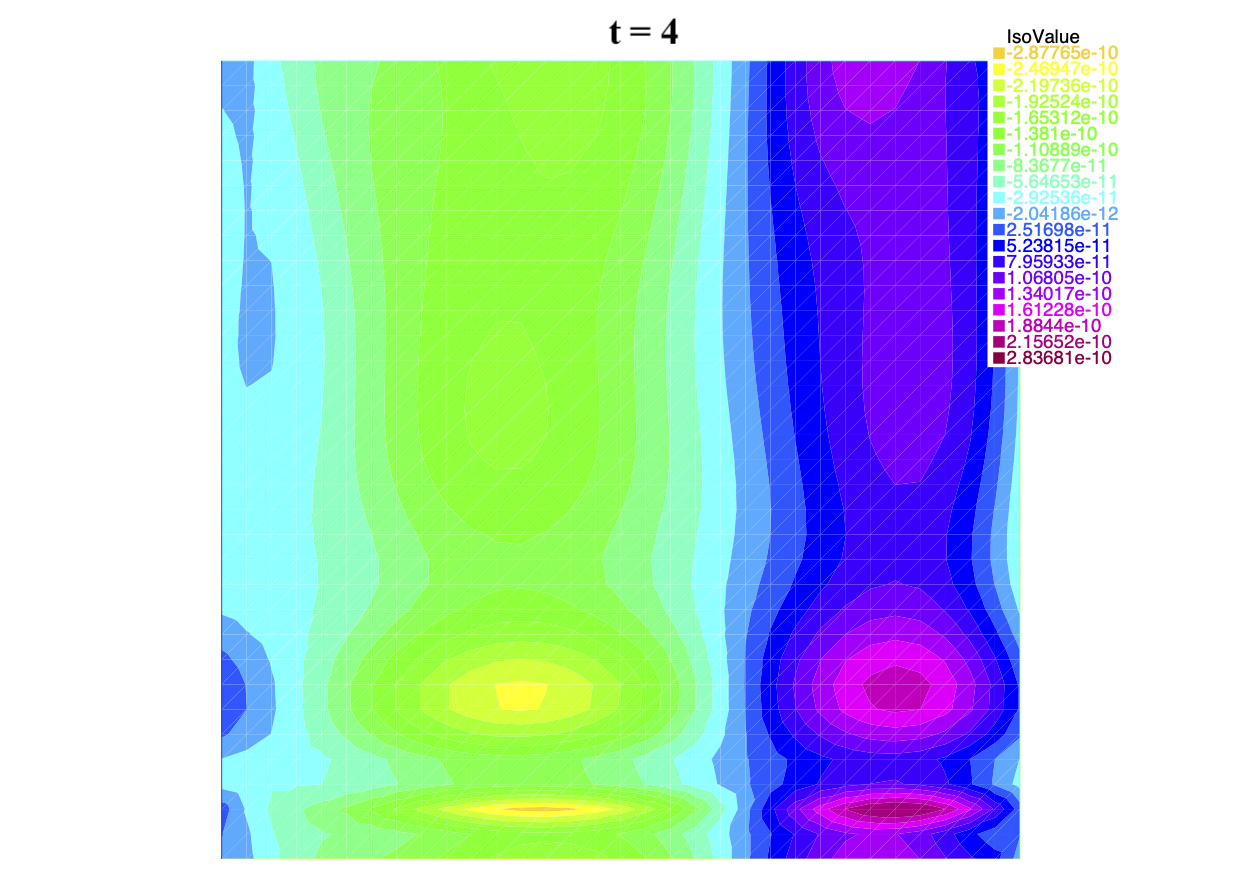}
\includegraphics[scale=0.25]{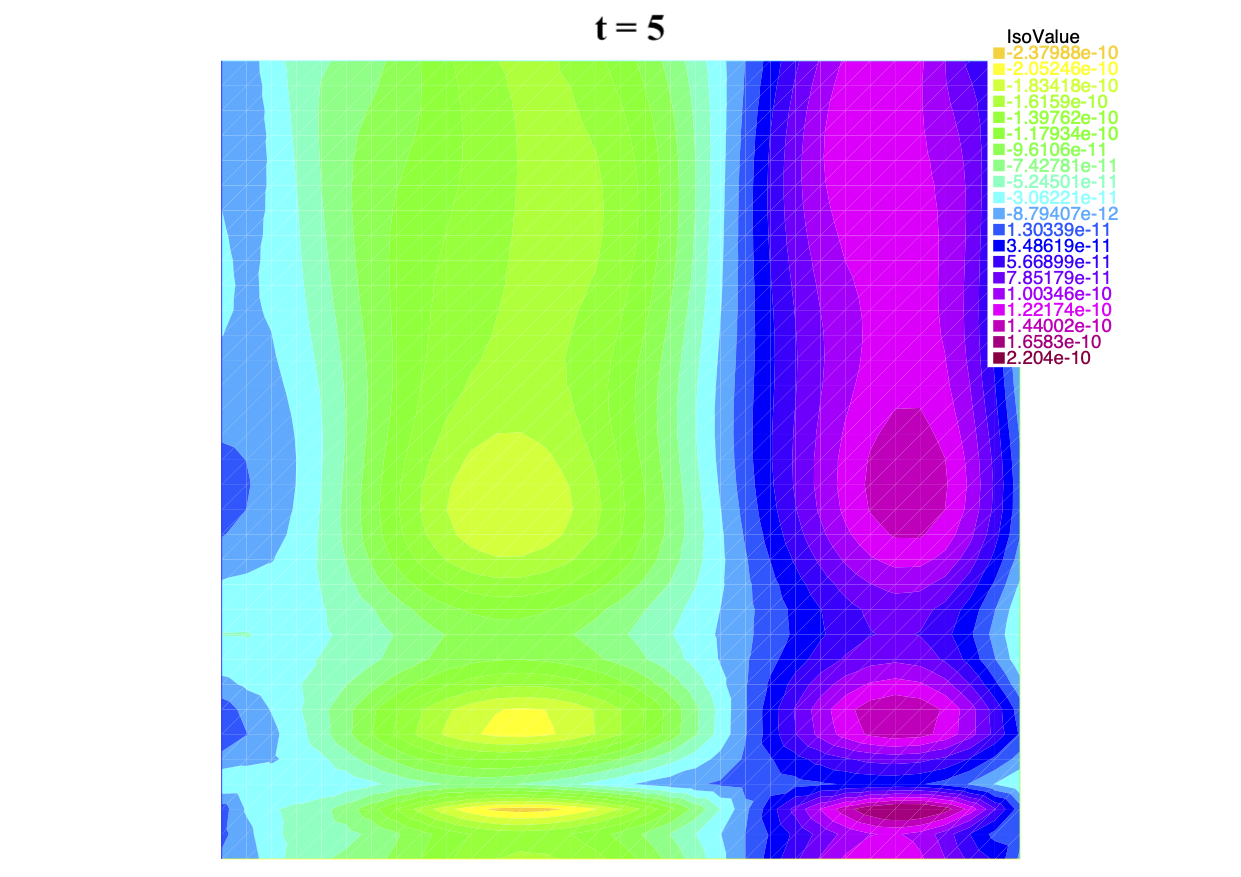}
\includegraphics[scale=0.25]{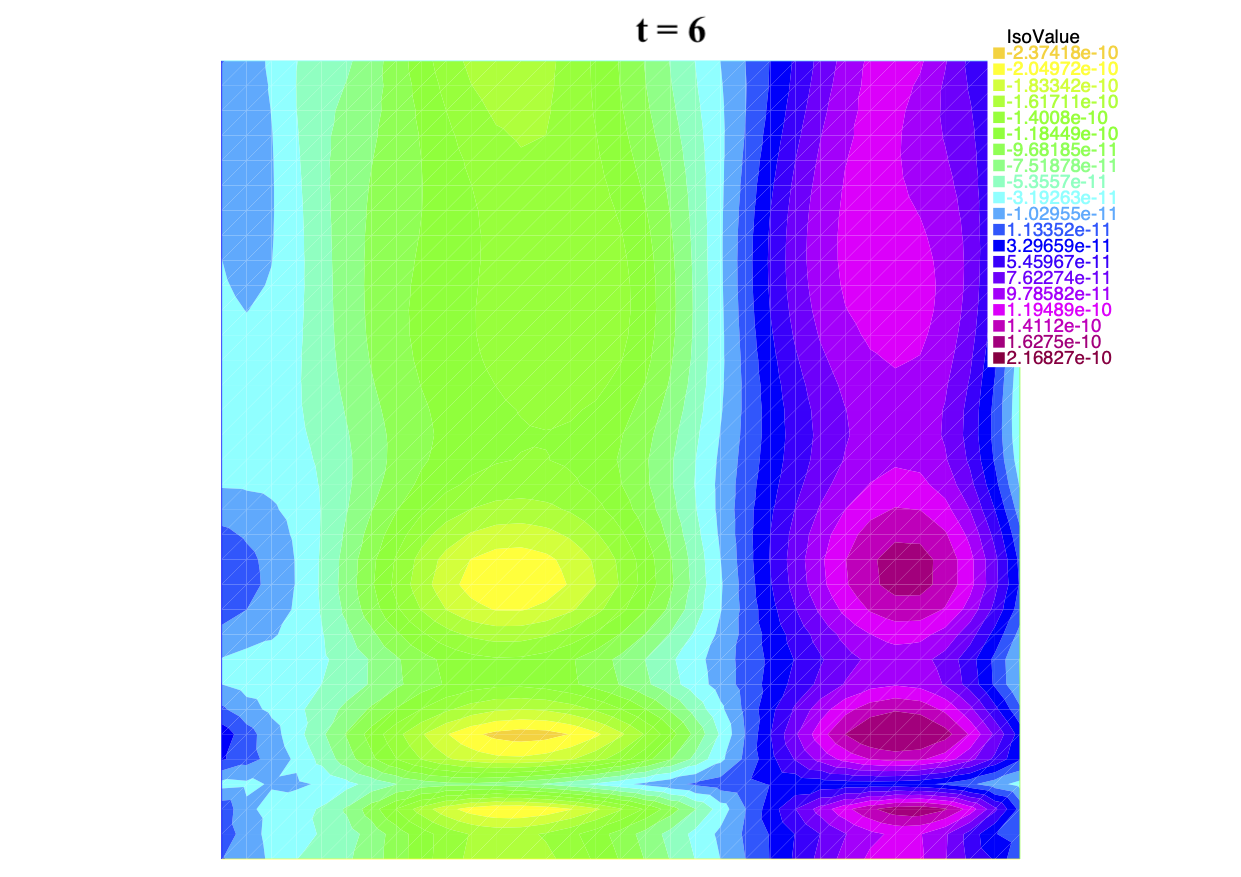}
\includegraphics[scale=0.25]{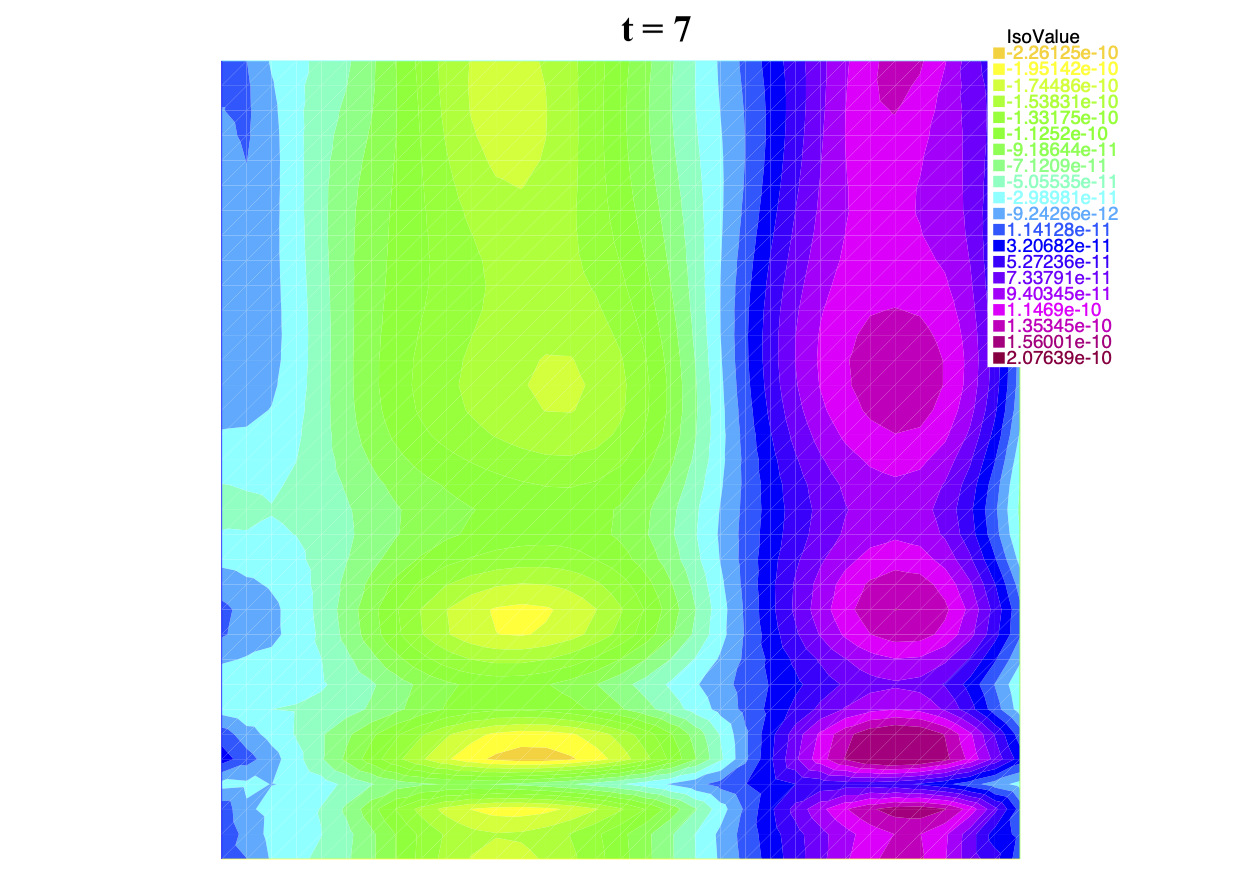}
\includegraphics[scale=0.25]{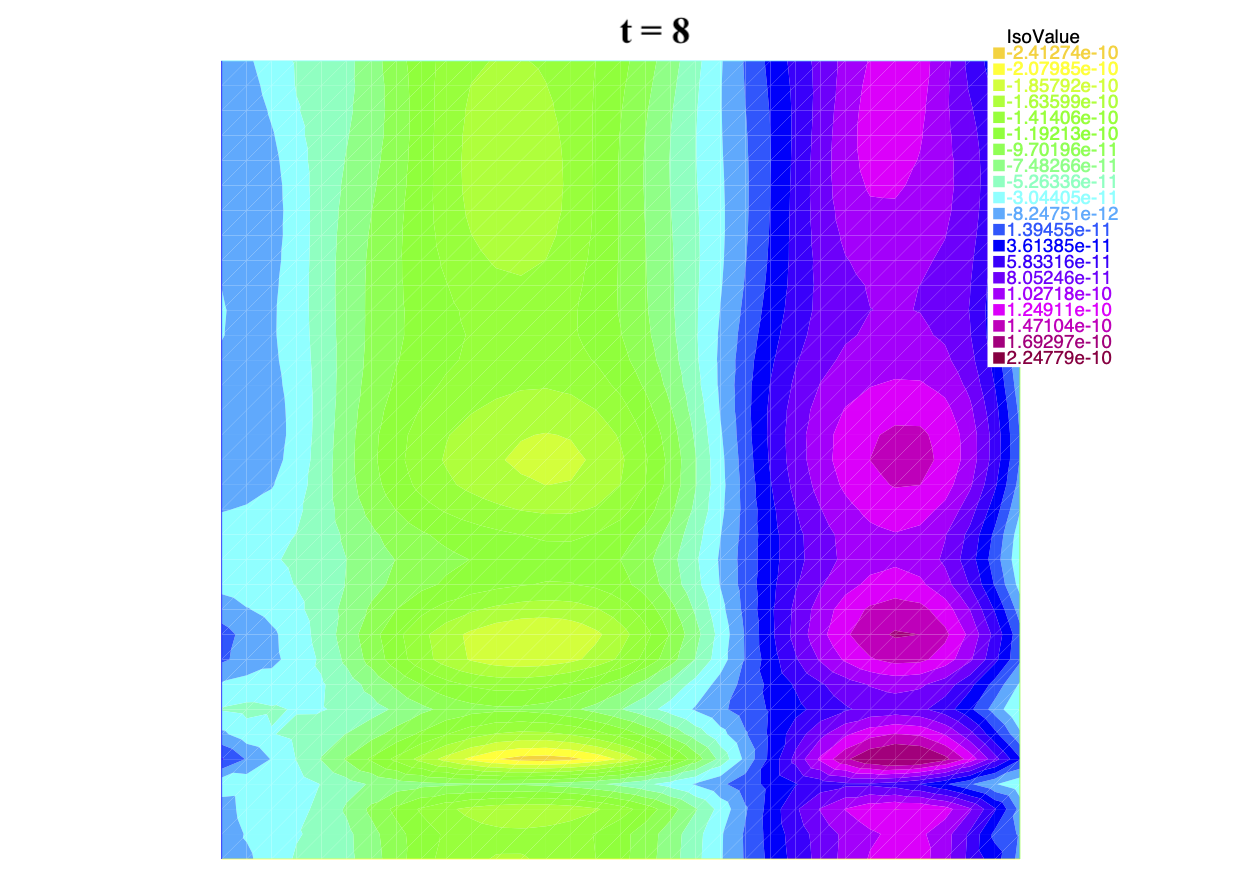}
\includegraphics[scale=0.25]{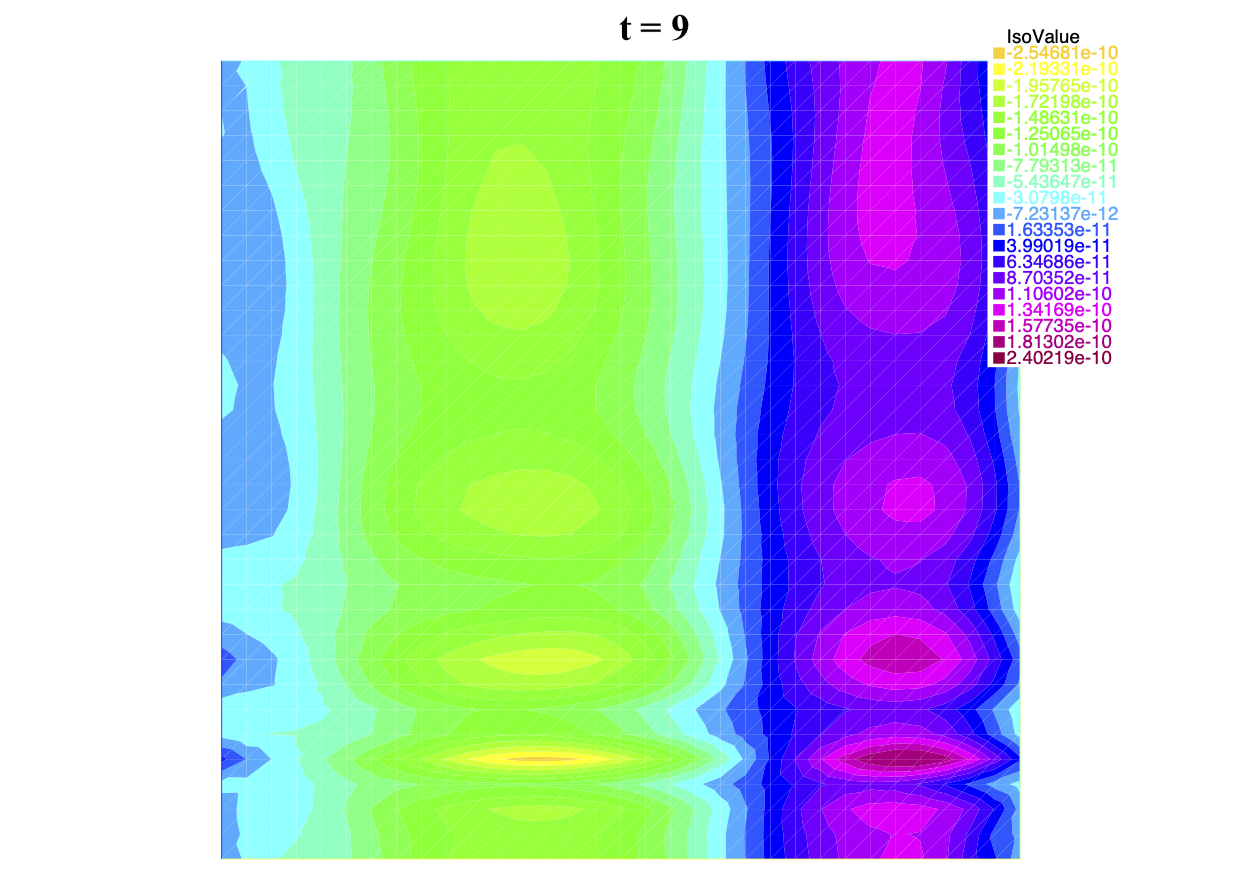} 
\includegraphics[scale=0.25]{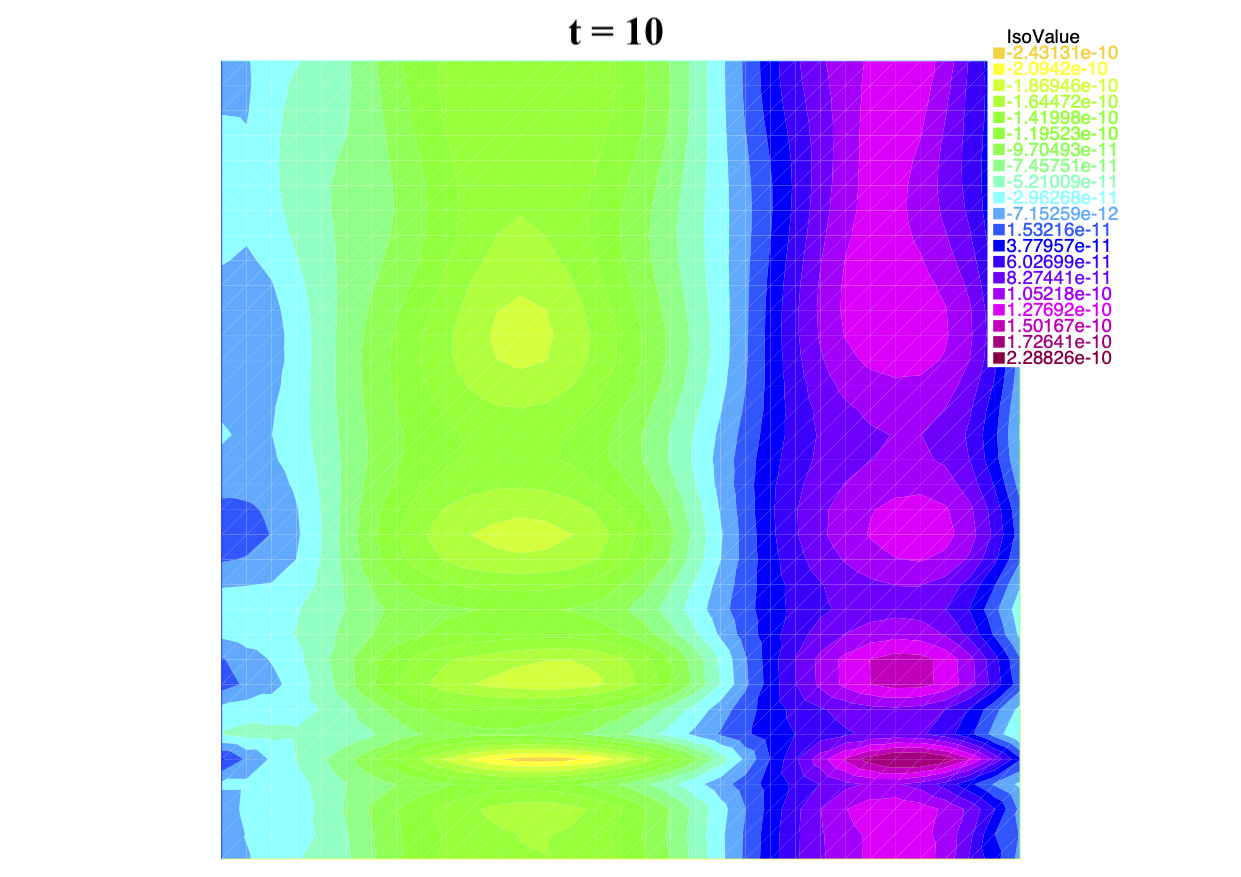} 
\includegraphics[scale=0.25]{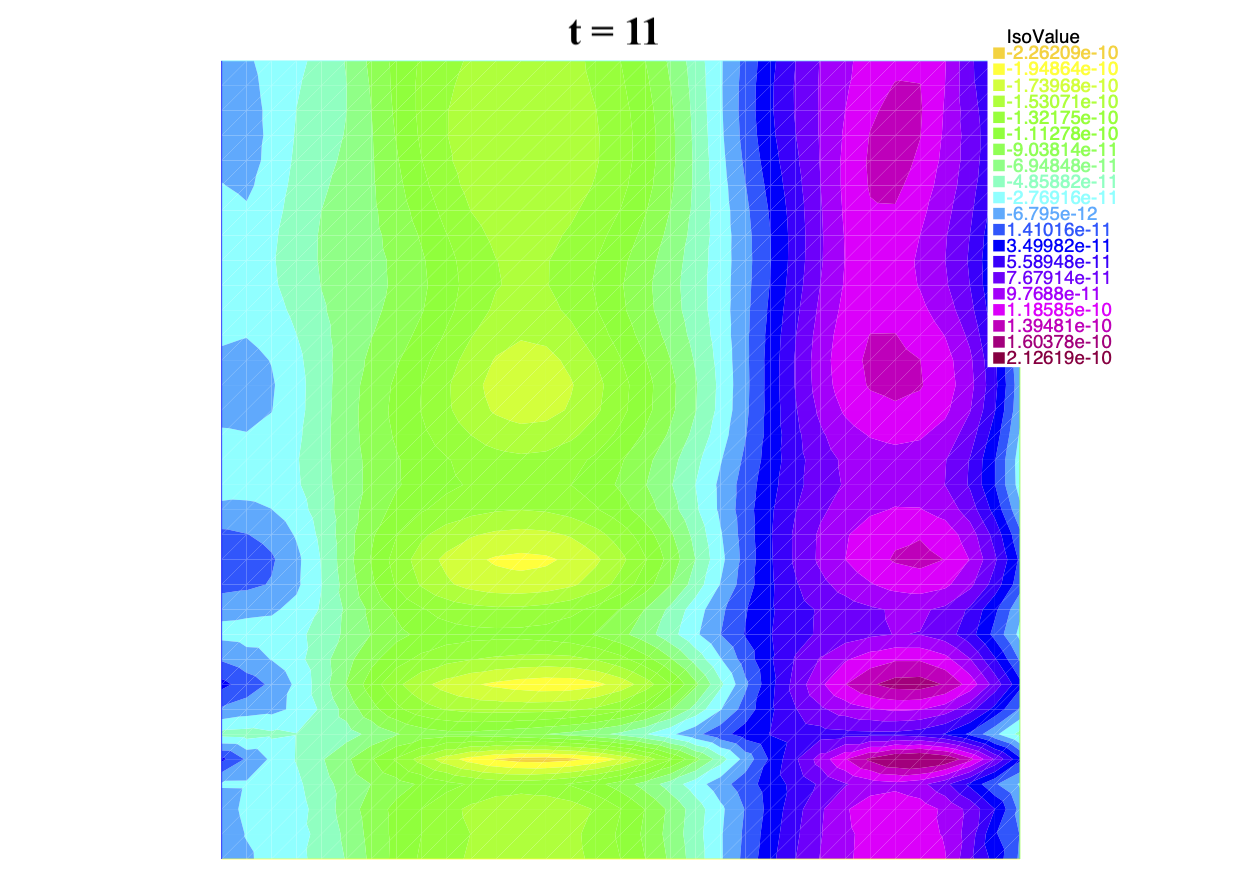} 
\includegraphics[scale=0.25]{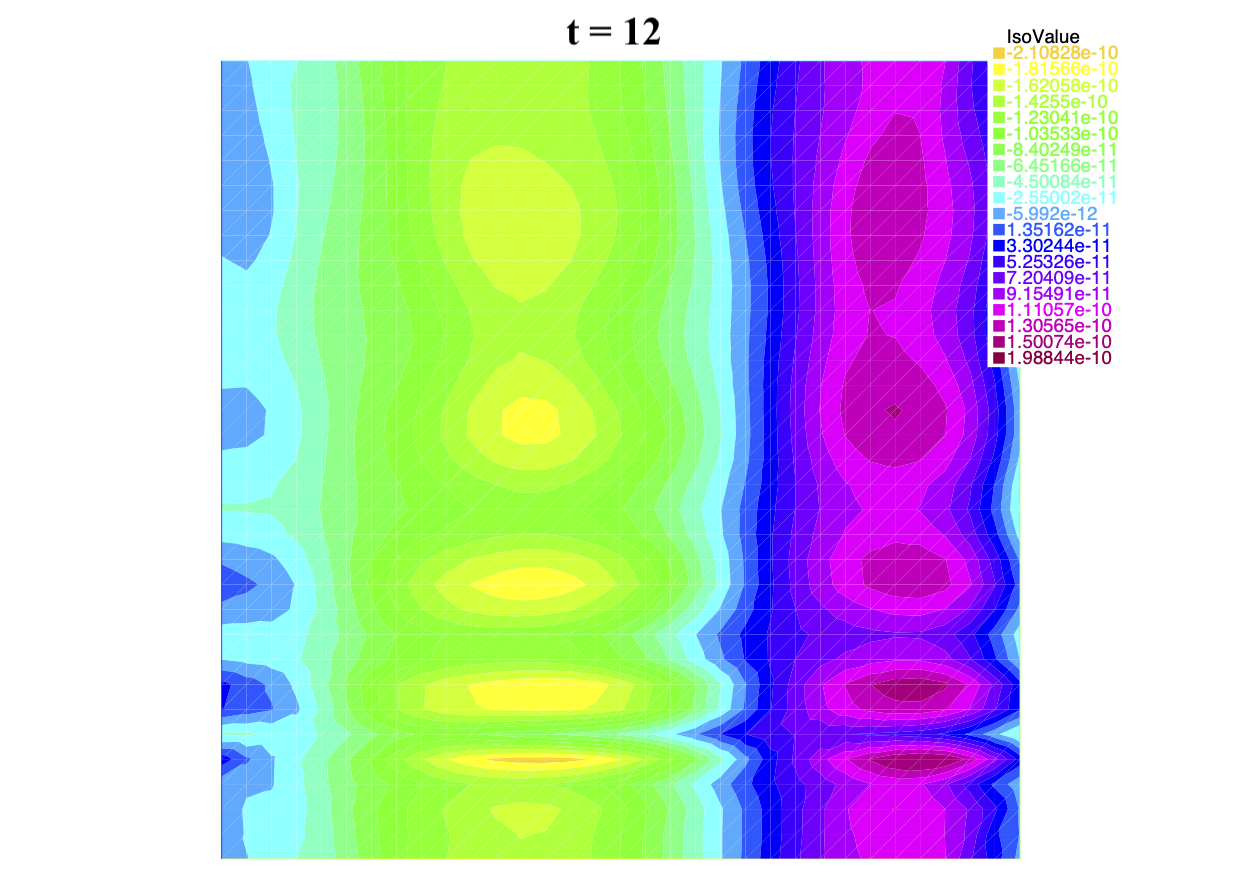} 
\includegraphics[scale=0.25]{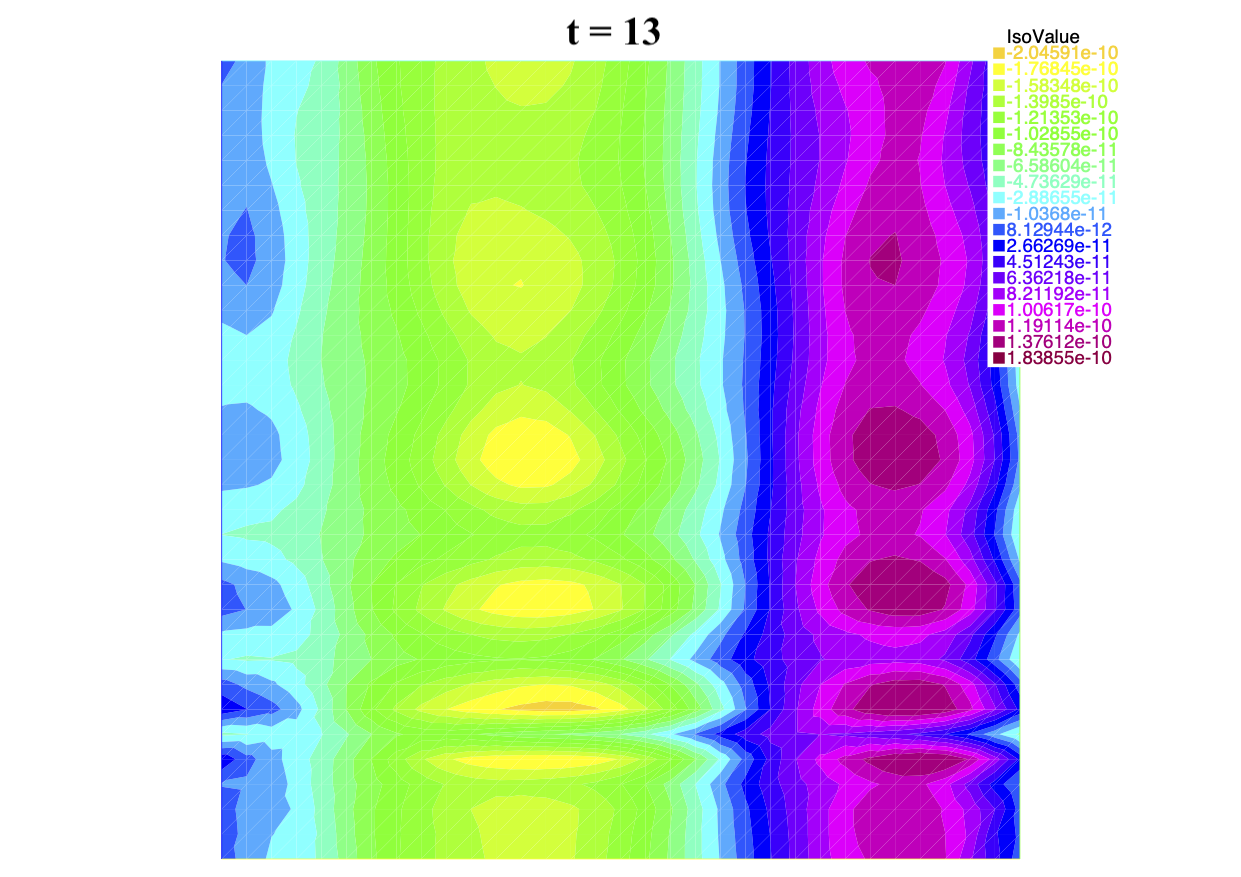} 
\includegraphics[scale=0.25]{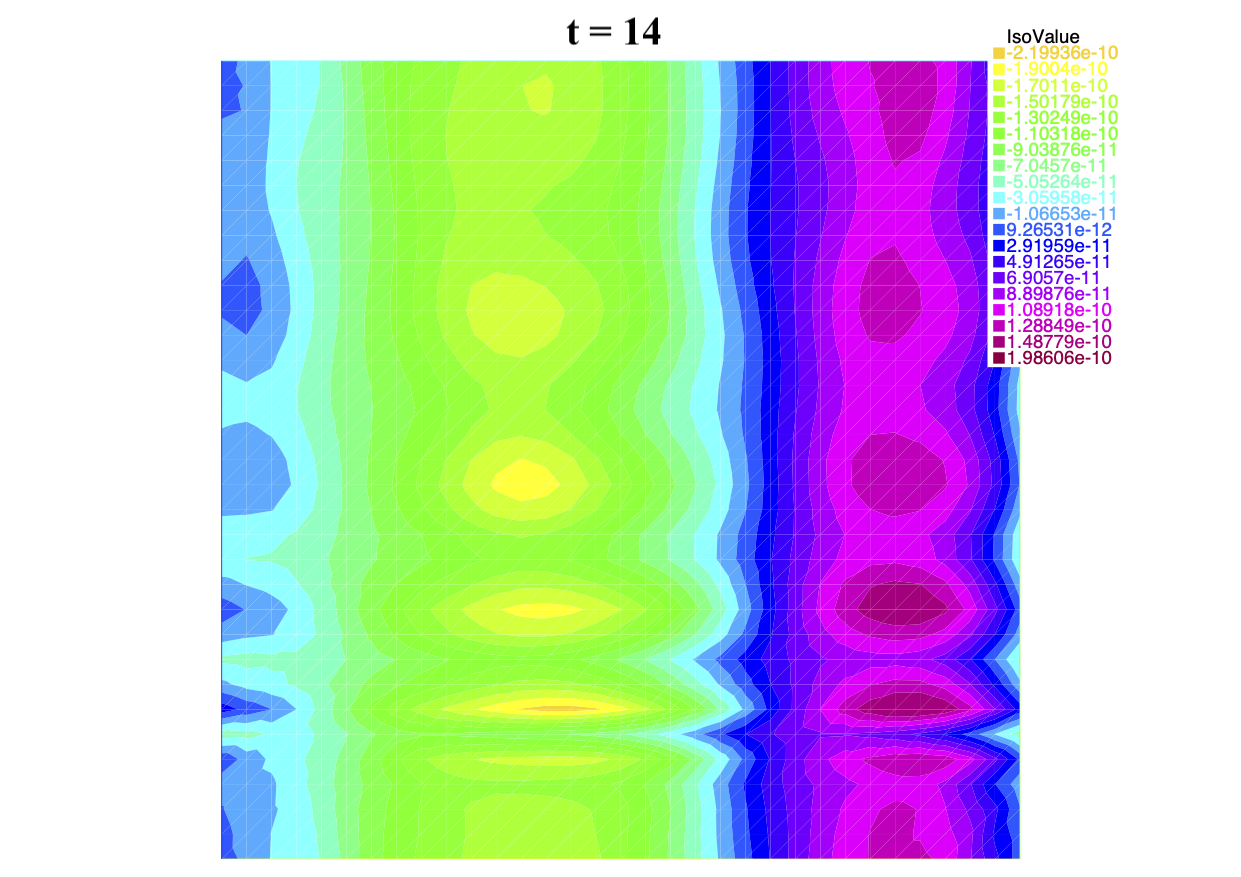} 
\caption{\it \small Time evolution of solution $u$ of $\eqref{HMC-Disc-Comp}$ for Test 4 using Algorithm \ref{alg:HMC-Modified Newton}, with $\tau = 0.1$, and a $33 \times 33$ grid on  $\Omega = [0,\pi] \times [0,\pi]$.}\label{fig:rand-32}
\end{figure}

\begin{figure}[H]
\centering
\includegraphics[scale=0.25]{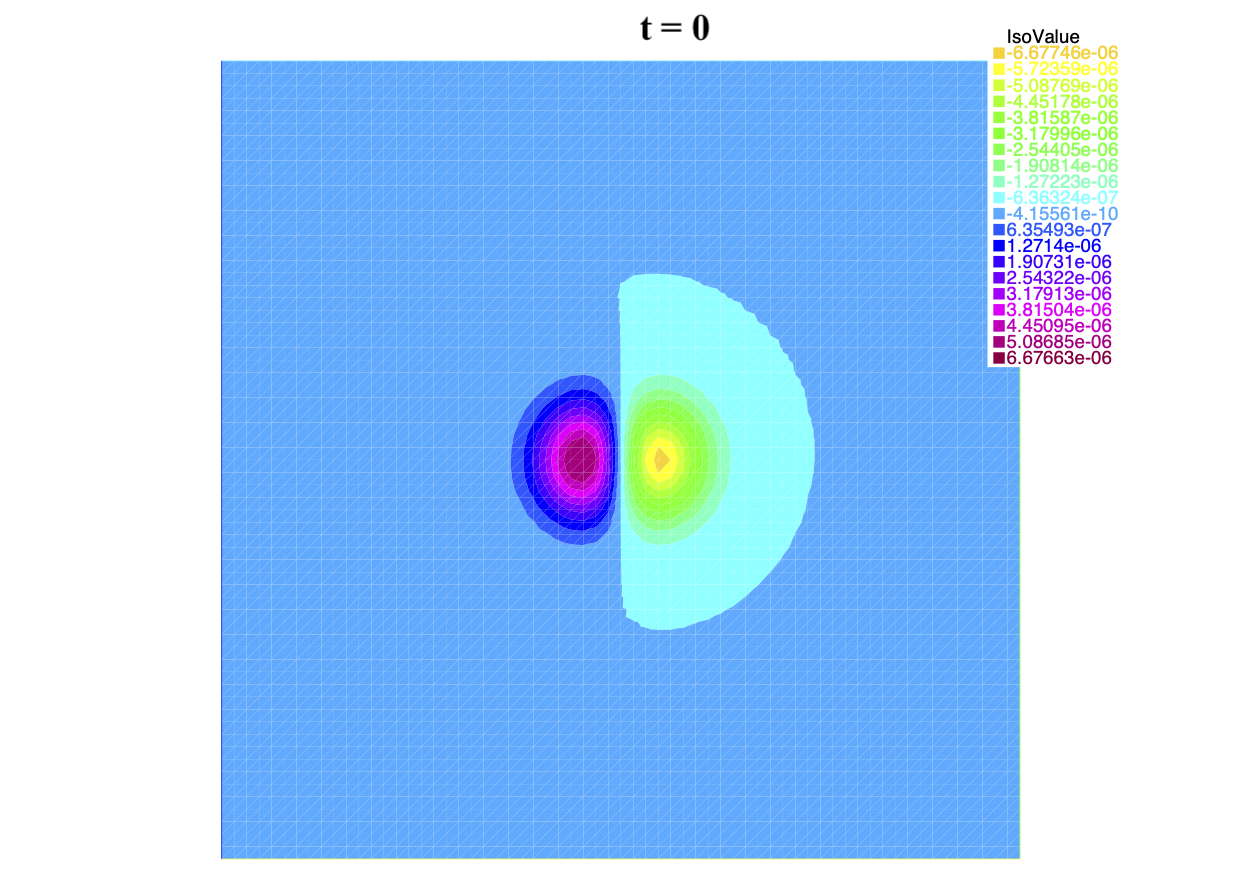}
\includegraphics[scale=0.25]{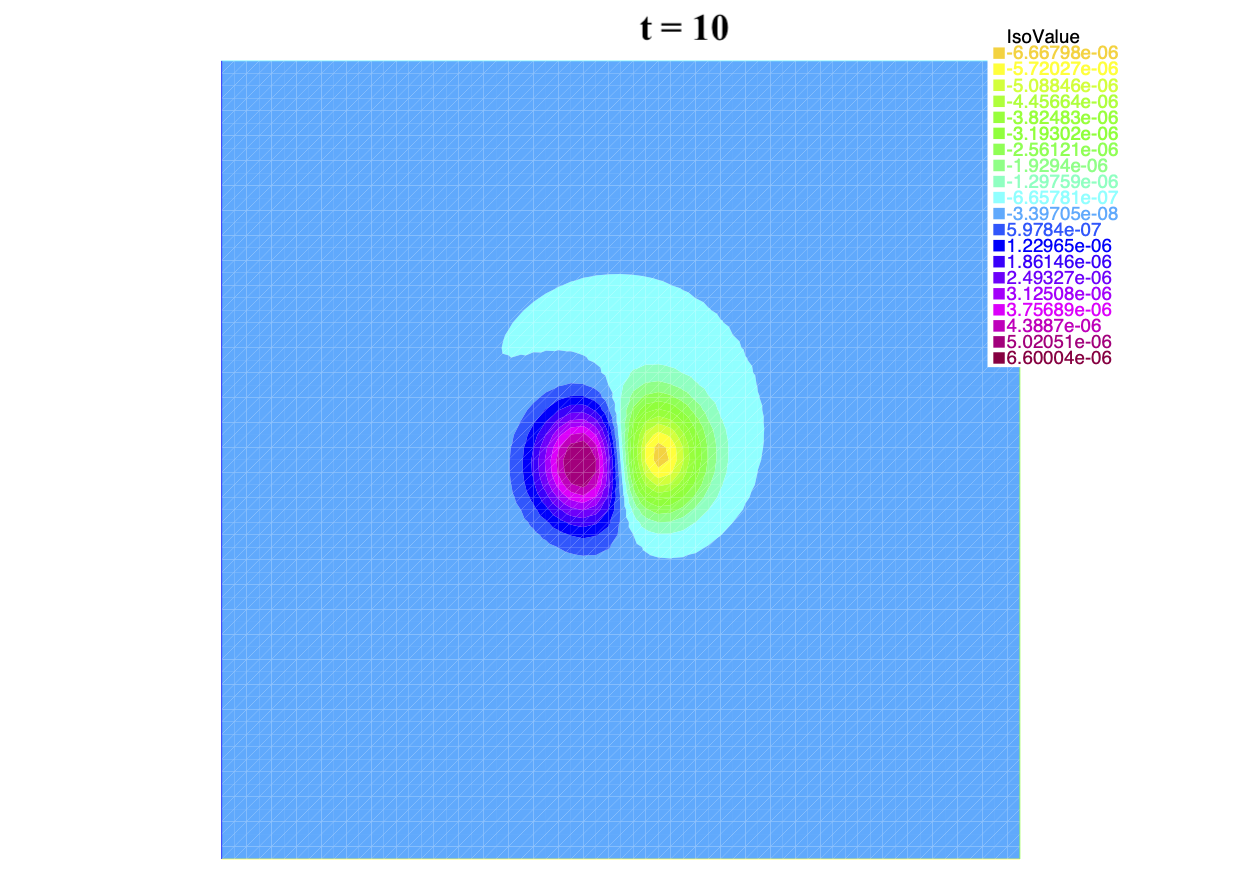}
\includegraphics[scale=0.25]{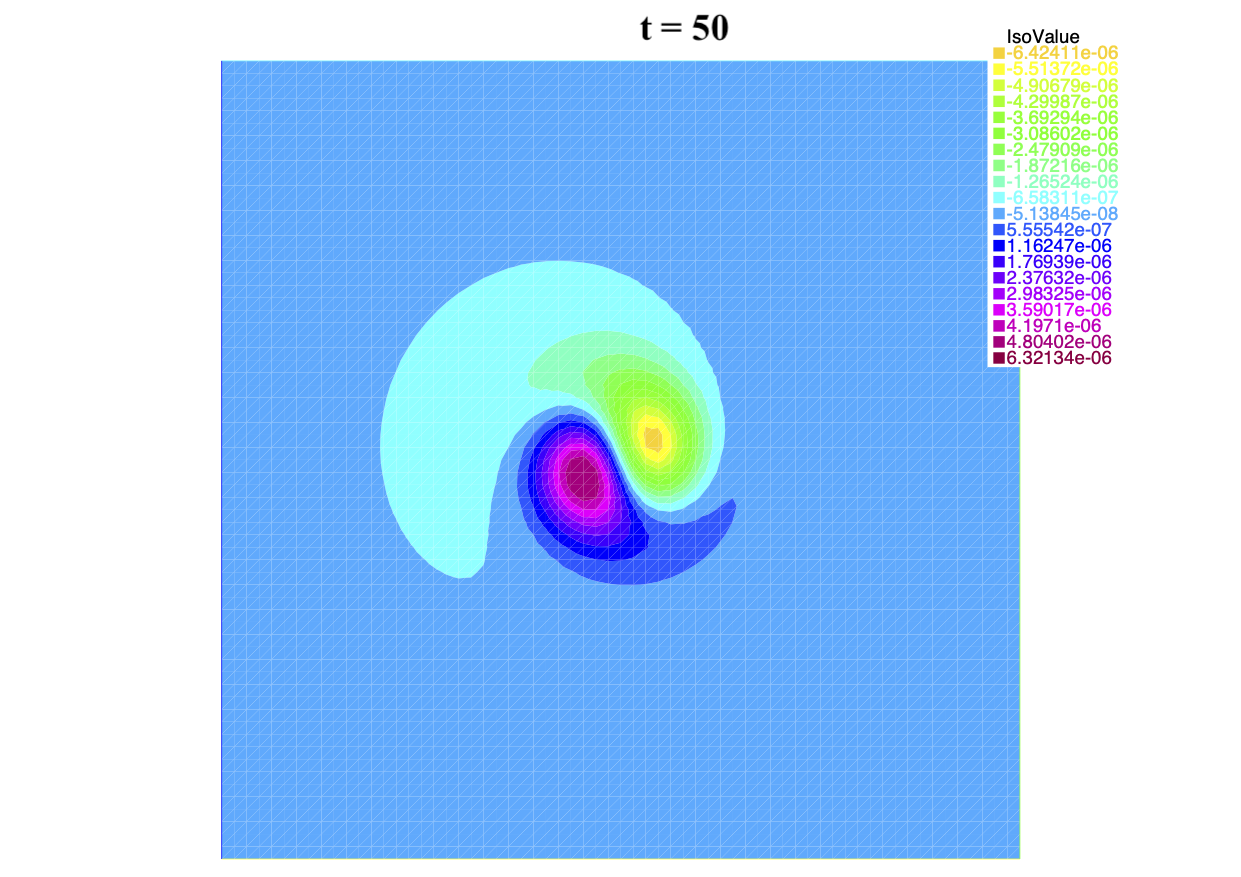}
\includegraphics[scale=0.25]{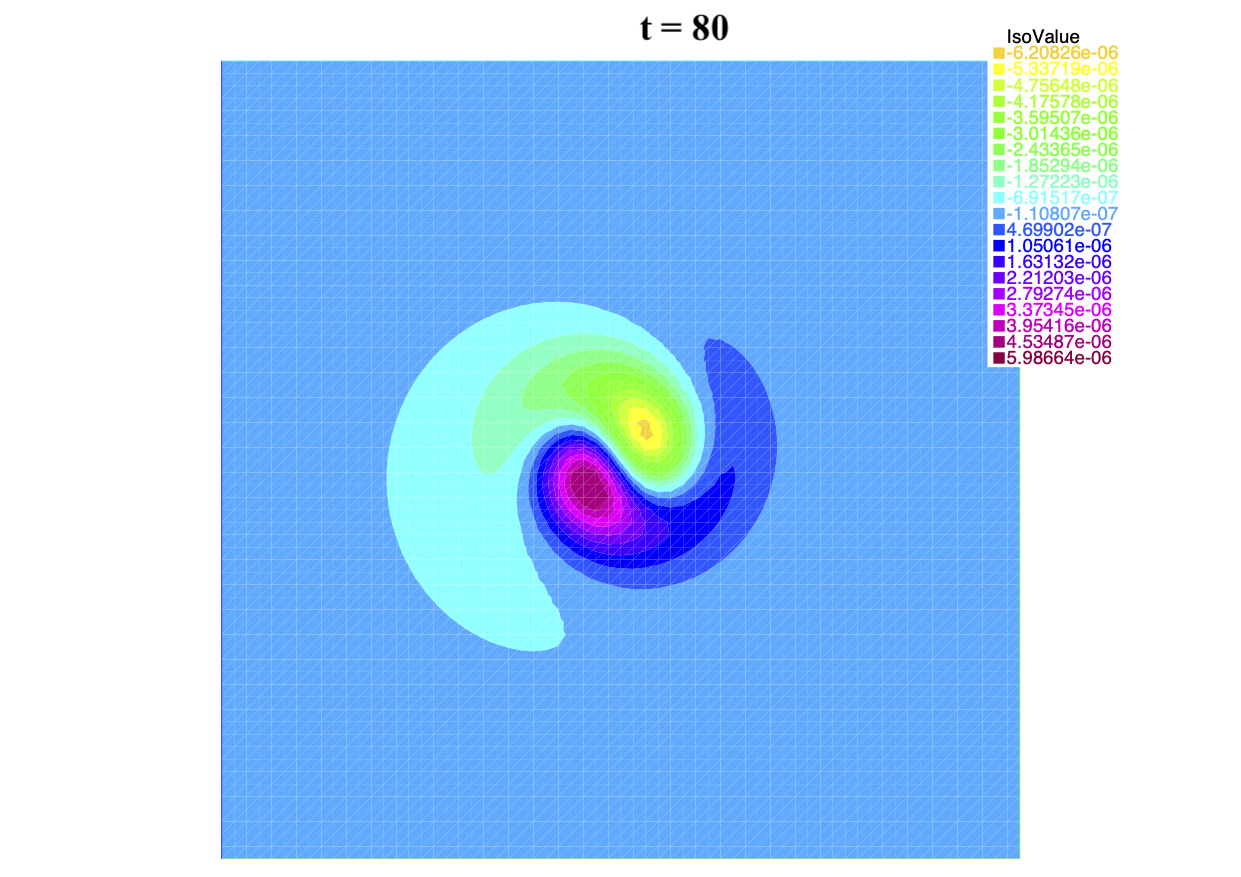}
\includegraphics[scale=0.25]{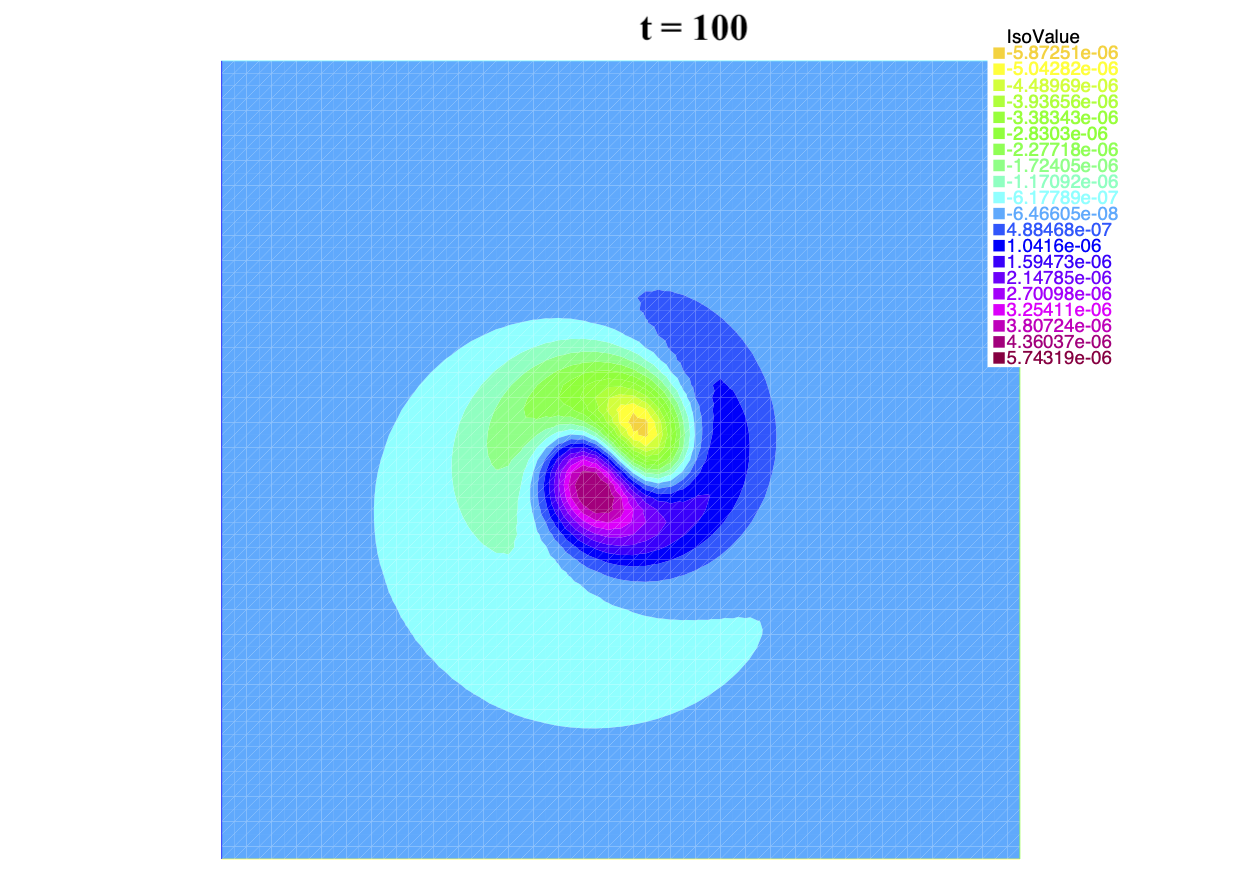}
\includegraphics[scale=0.25]{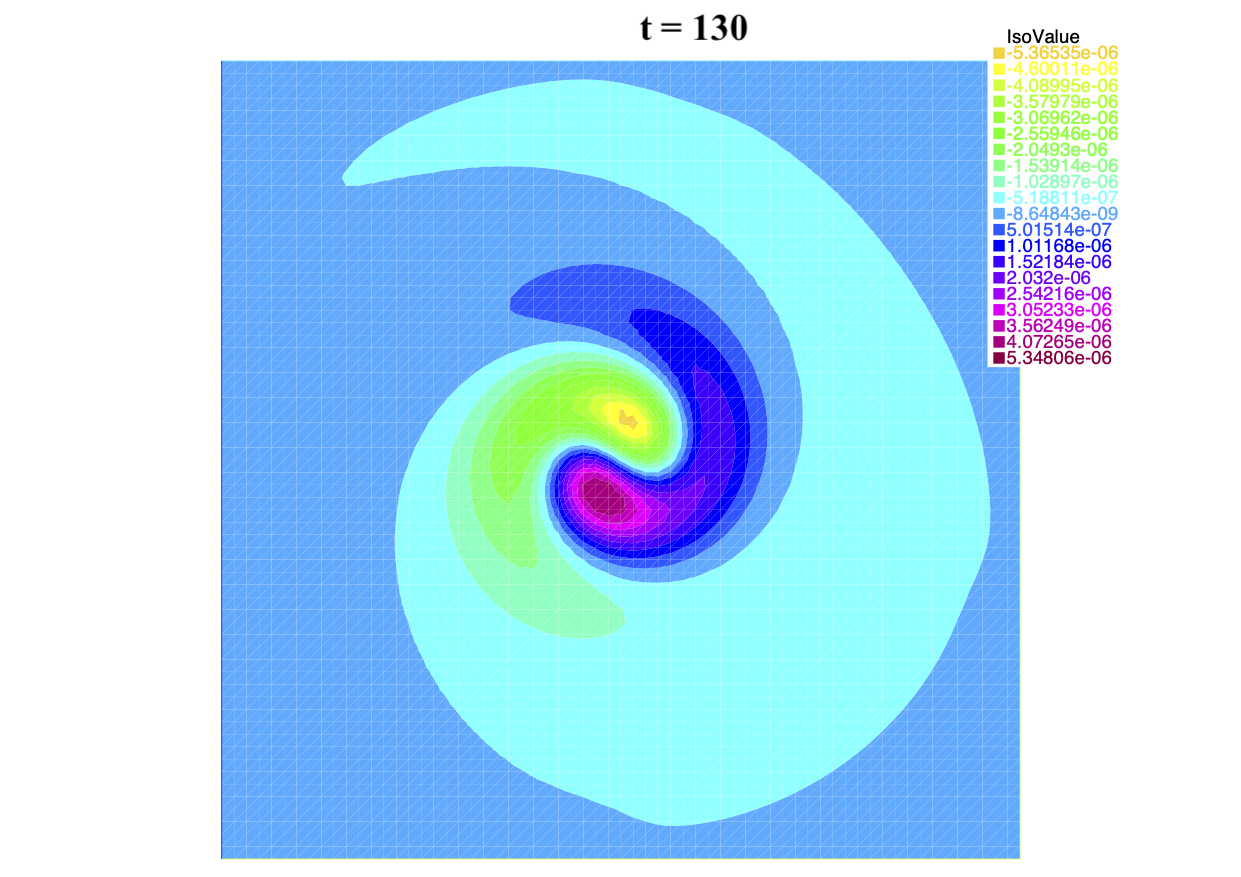}
\includegraphics[scale=0.25]{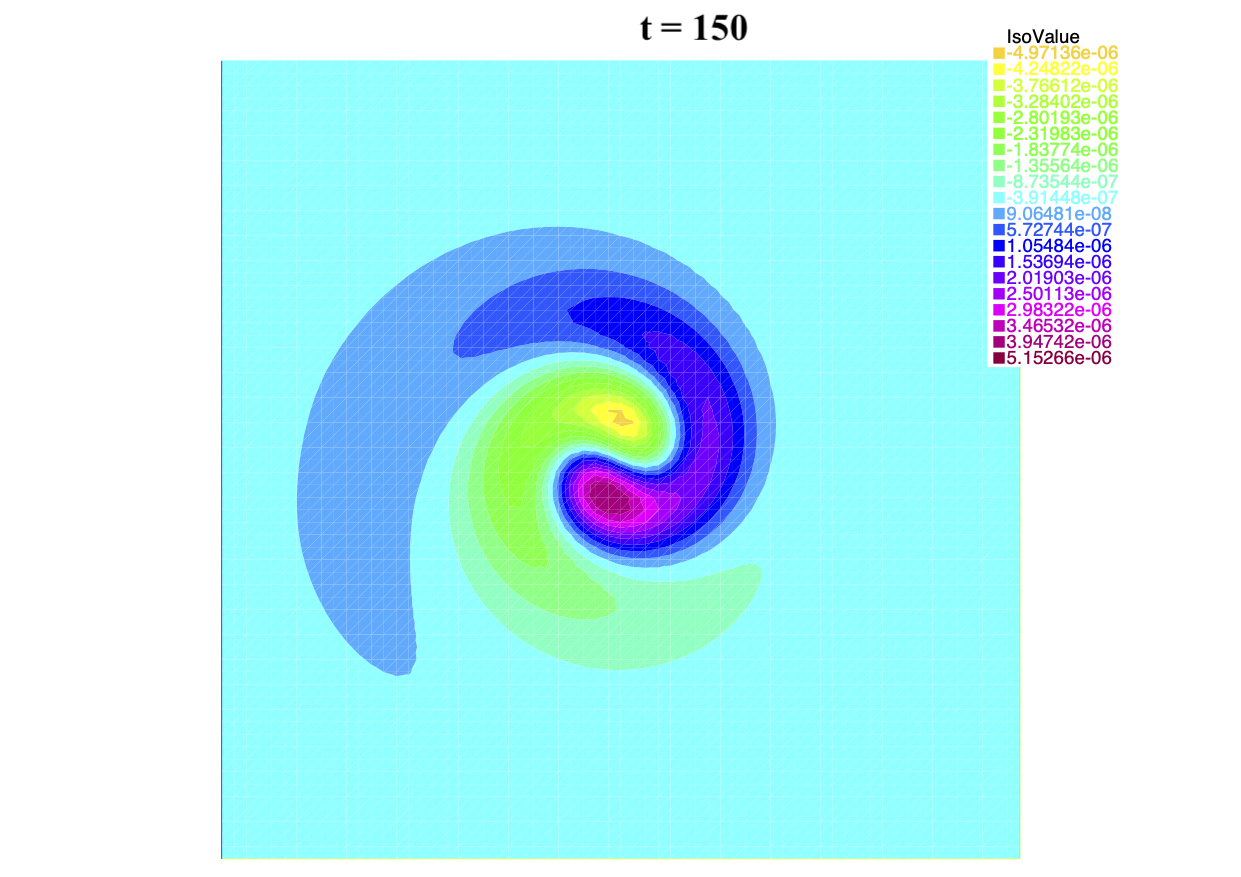}
\includegraphics[scale=0.25]{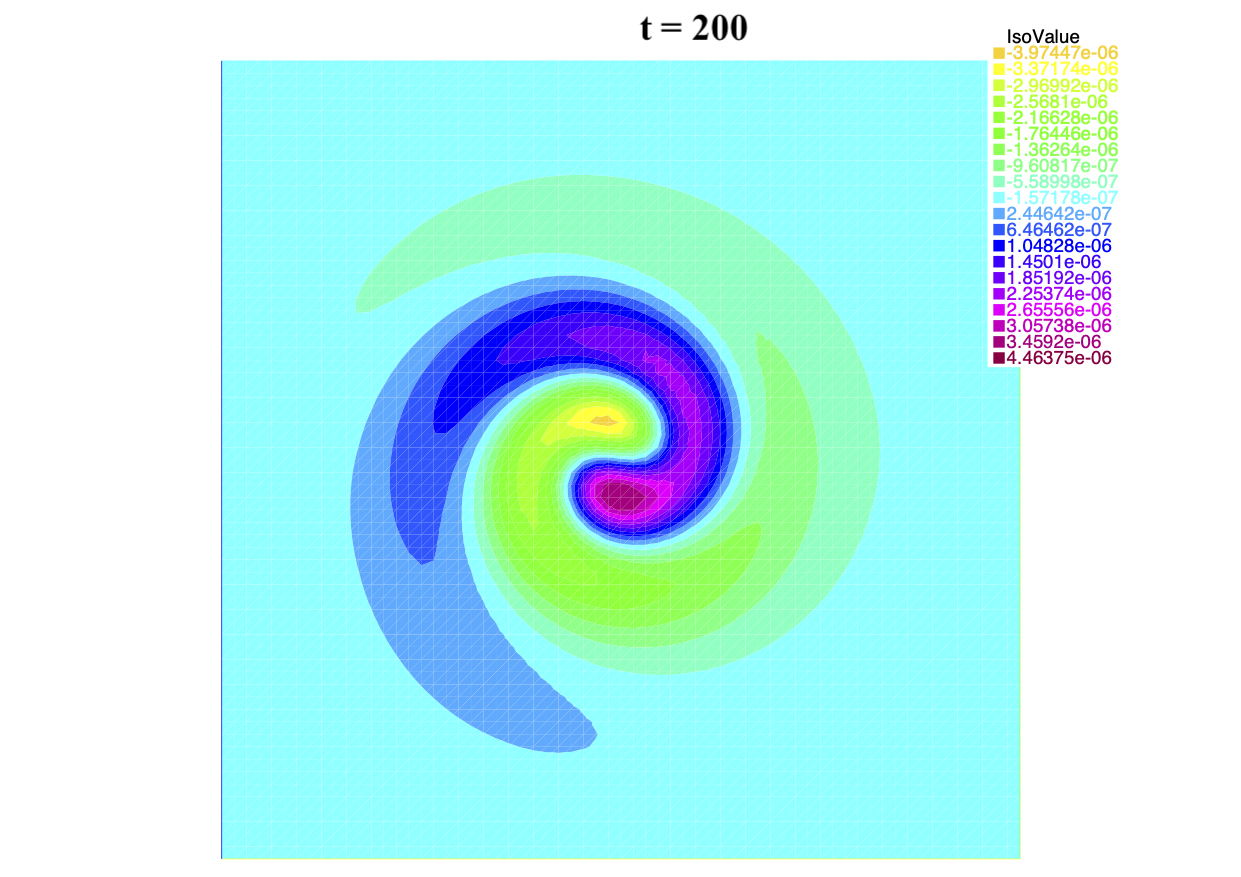}
\includegraphics[scale=0.25]{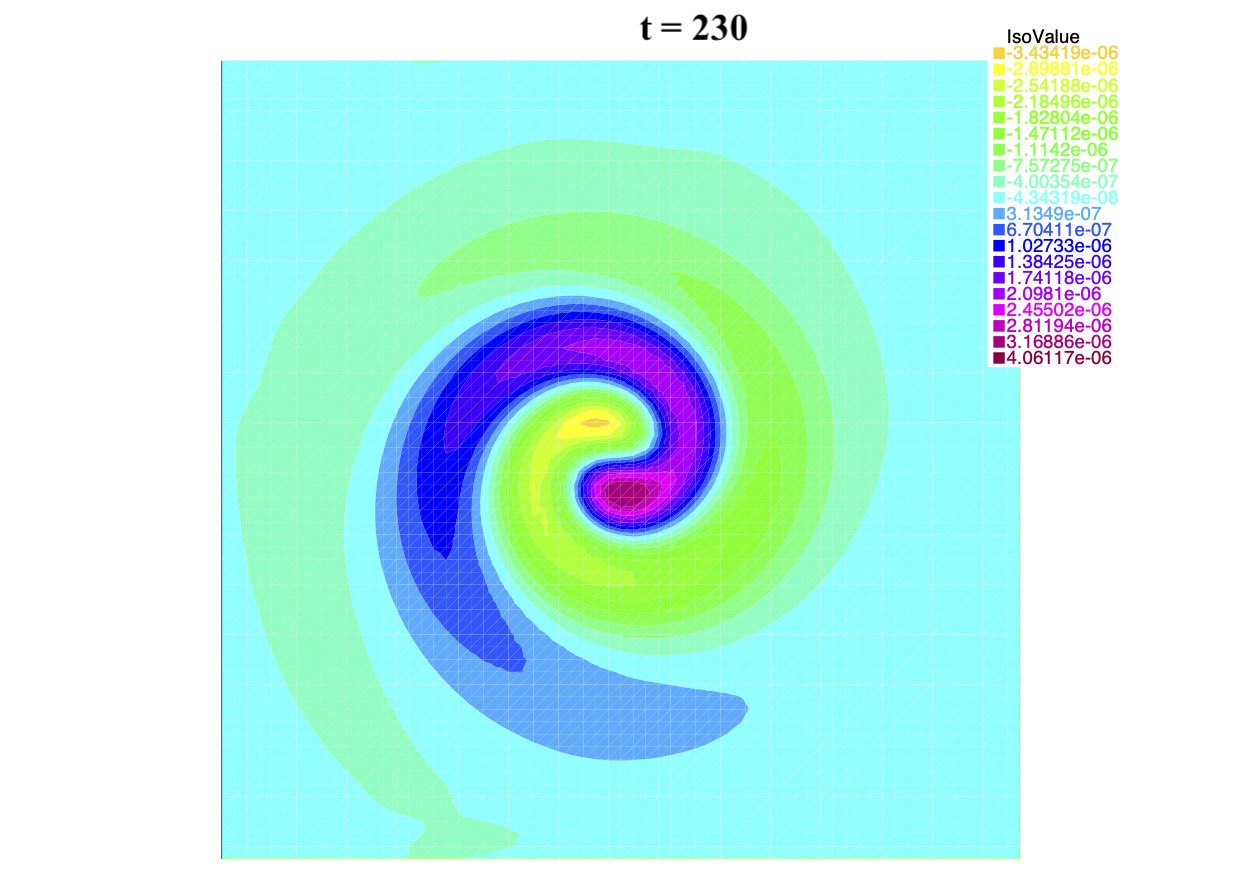}
\includegraphics[scale=0.25]{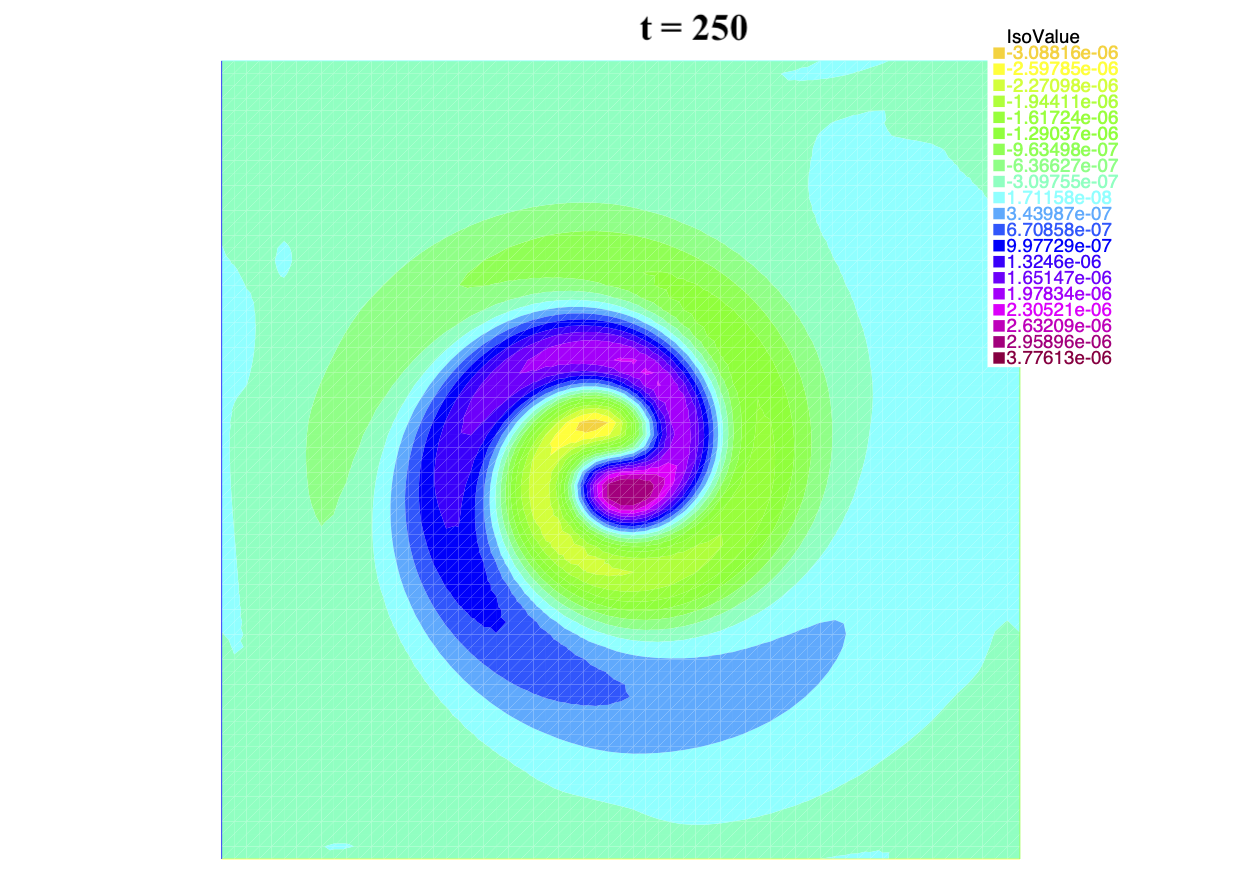} 
\includegraphics[scale=0.25]{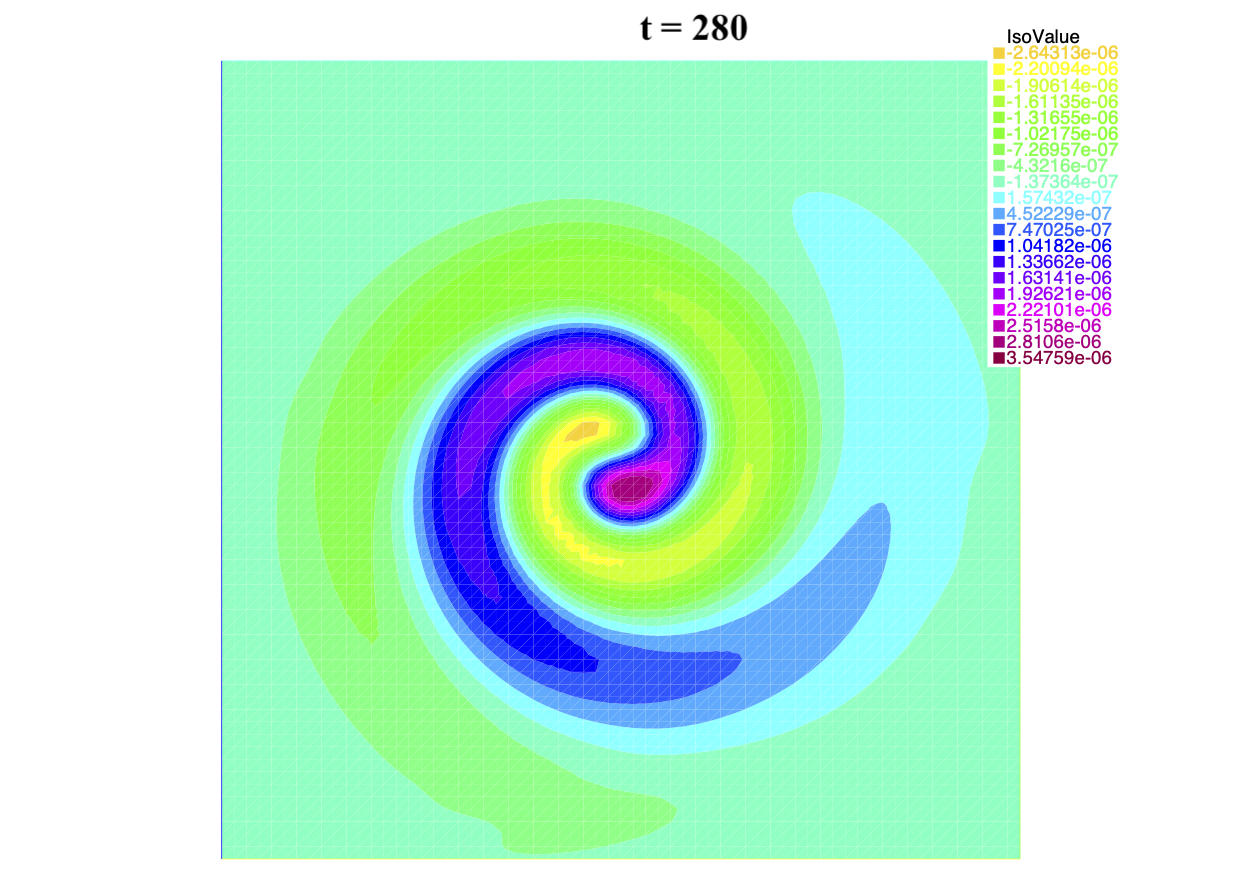} 
\includegraphics[scale=0.25]{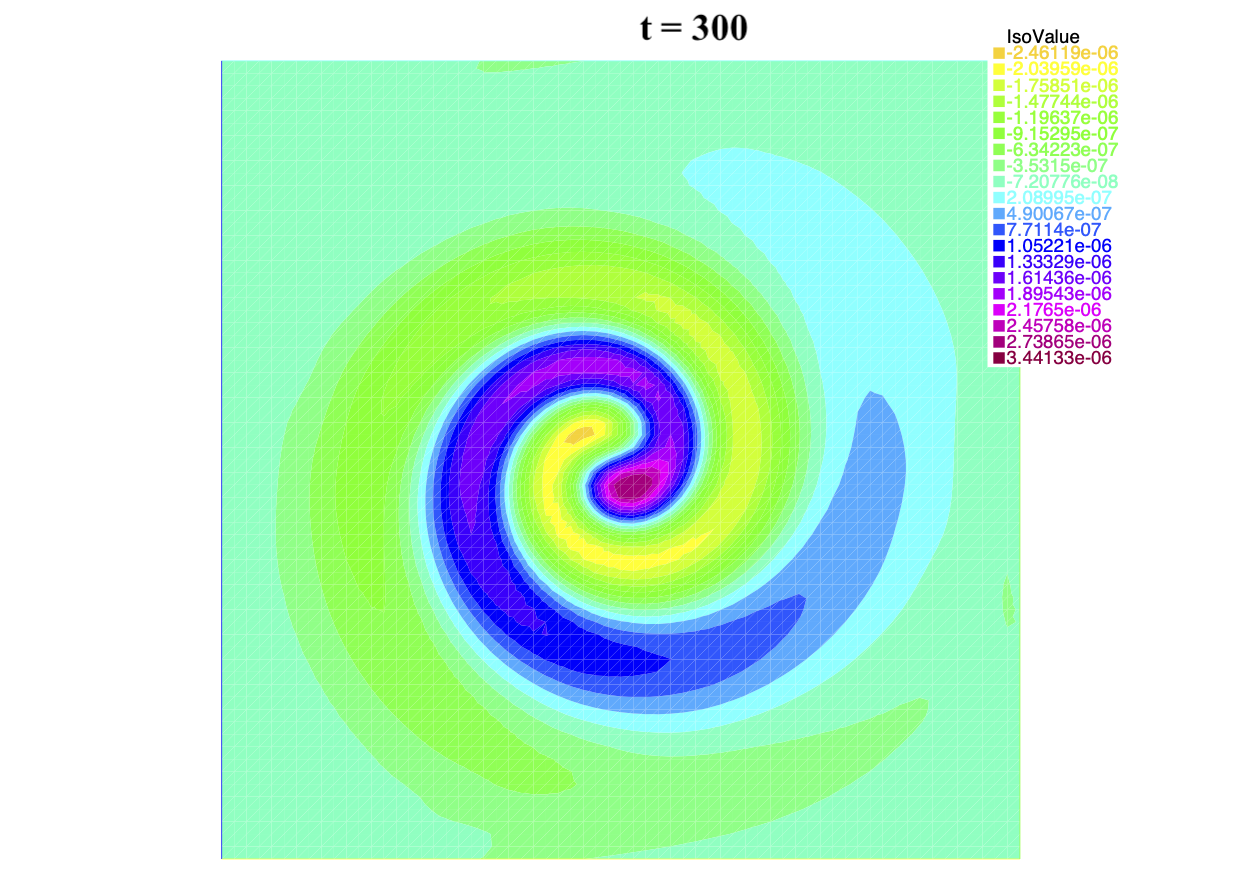} 
\includegraphics[scale=0.25]{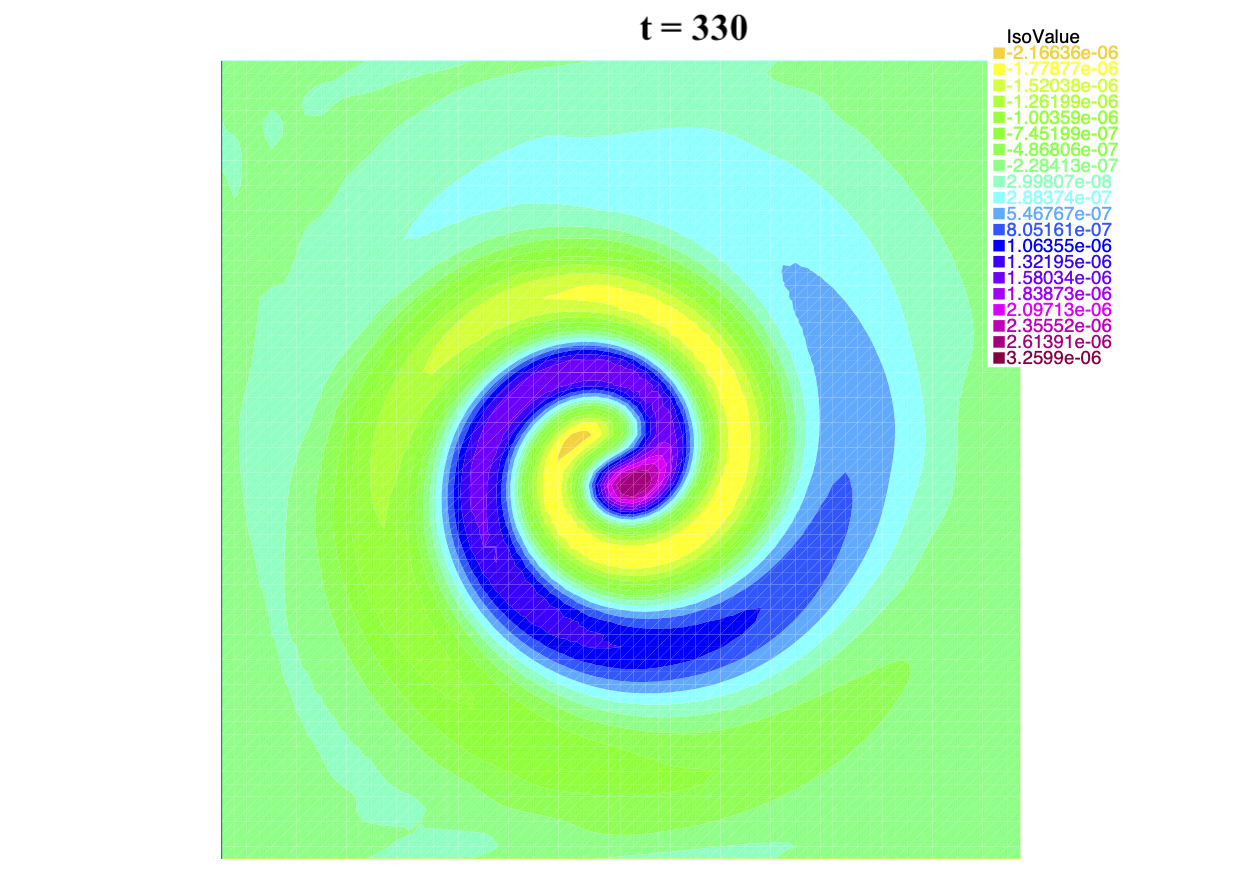} 
\includegraphics[scale=0.25]{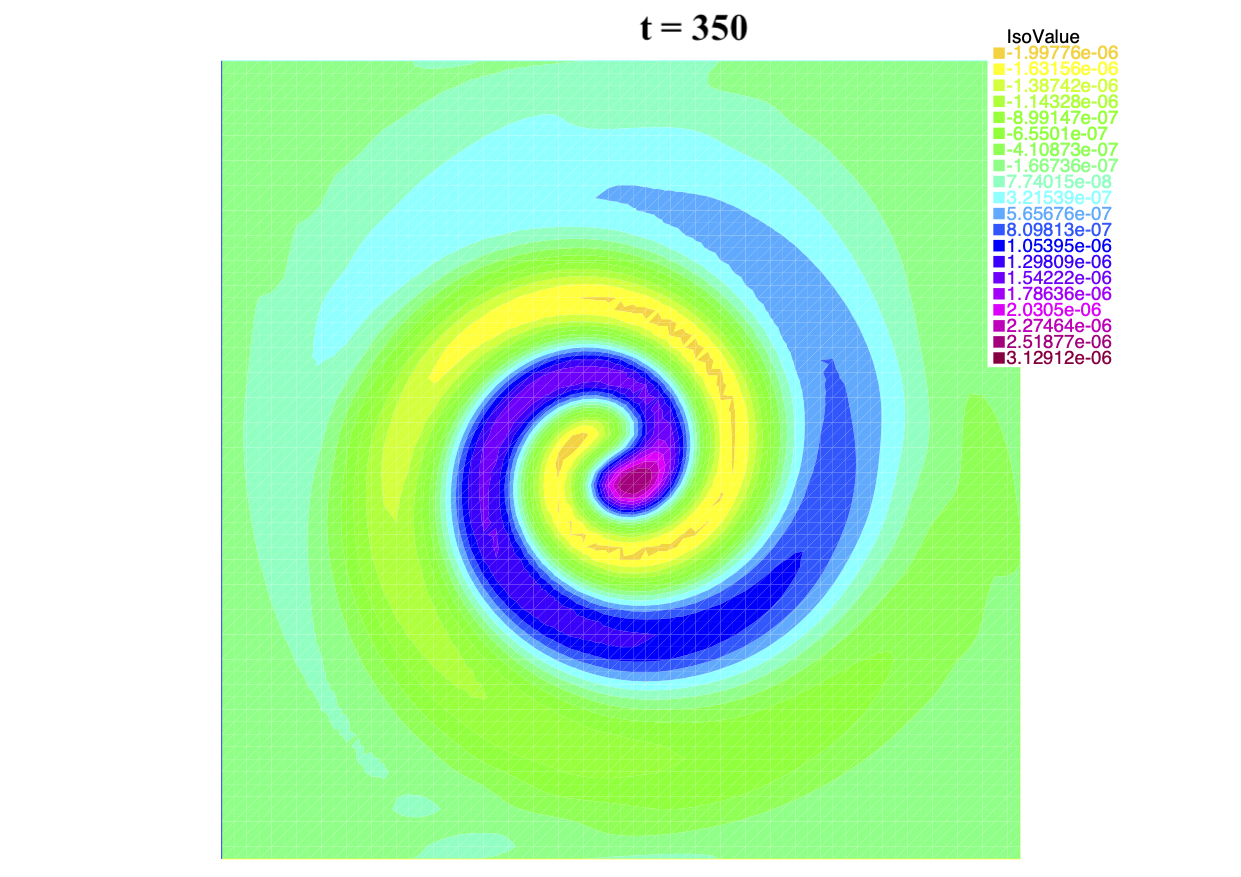} 
\includegraphics[scale=0.25]{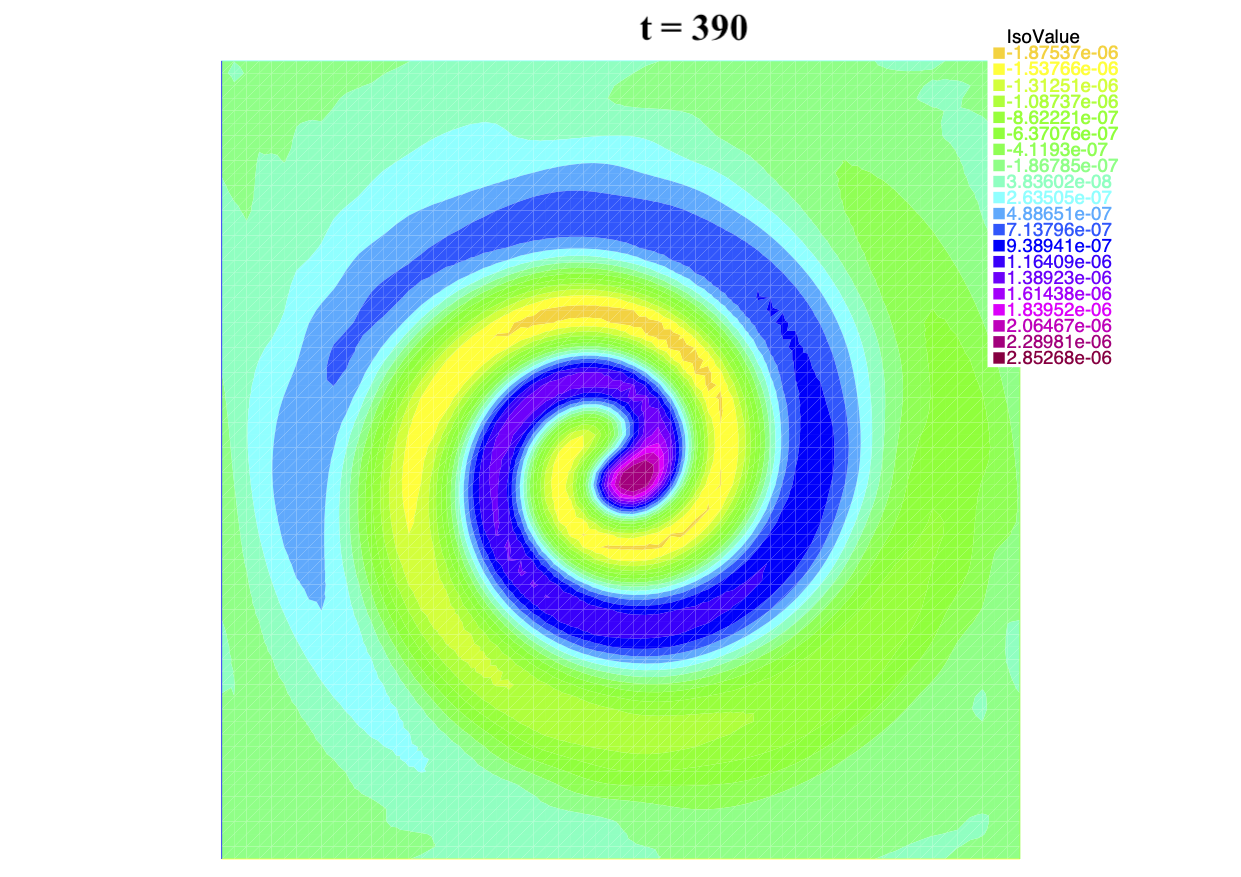} 
\caption{\it \small Time evolution of solution $u$ of $\eqref{HMC-Disc-Comp}$ for Test 5 using Algorithm \ref{alg:HMC-Modified Newton}, $\tau = 0.1$, and a $65 \times 65$ grid on  $\Omega = [0,20] \times [0,20]$.}\label{fig:gauss-64}
\end{figure}

Apart from runtime, it should be noted that the evolution of the solution with  respect to time was similar for the three methods. Thus, in Figures \ref{fig:sin9piy}-\ref{fig:gauss-64} we show the evolution of the solution for the fastest method, Modified Newton (Algorithm \ref{alg:HMC-Modified Newton}), for $n = 33$ or $65$ partition points in each direction, and $\tau = 0.1$. 

For Tests 1-4, $p(x,y) = 12x$ is of the form $Ax+b$, where $ p_x = 12$ and $p_y = 0$. Thus, the solution is expected to be a traveling wave in the y-direction for a nonzero $A$. The speed of the motion and its direction depend on the magnitude and sign of $A$ respectively. If $A>0$, then the motion is in downwards, whereas if $A<0$ the motion is upwards, and for $A = 0$ no motion. The larger the the magnitude of $A$, the faster the motion.\\ For example in Tests 1 and 2 (Figures \ref{fig:sin9piy}, and \ref{fig:sin3y32}) where $u_0$ is a $\sin$ function in the y-direction, the traveling wave effect in the y-direction is clear. 

Whereas in Test 3 (Figure \ref{fig:sin3x-32}), it may seem that the solution is stationary. However, the motion is along the y-direction, and the initial solution which is a $\sin$ function in the x-direction, is fixed for $x = a$, $a\in \mathbb{R}$. Thus, the traveling wave effect is not visible. But if the  $\sin$ function is multiplied by other function so that it is no longer fixed for $x = a$, like Test 4 (Figure \ref{fig:rand-32} ), then the motion is visible again. Note that if we set $p_x = 0$ and $p_y = 12$, then the solution will be moving in the x-direction in a similar manner, where for Tests 1 and 2, the solution will appear stationary and Test 3 and 4 will be traveling in the x-direction.

\newpage
As for Test 5 (Figure \ref{fig:gauss-64}), the solution is expected to have a circular motion around the center of the domain $(10,10)$, since $\nabla p = [-(x-10)/32, -(y-10)/32]$, which is observed. 
 Note that if $\nabla p = [(x-10)/32, (y-10)/32]$ for the same initial conditions, then the solution will be moving in the opposite direction. 

As for the comparision between the Newton-type methods for solving \eqref{HMC-Disc-Comp} and the semi-linear approach \eqref{HMC-Disc-Comp-semi} introduced in \cite{FEHM}, we note the following. The corresponding expected behavior was observed in all methods with one main difference.
In the Newton-type methods that were Tested for end time $T = 300$, the maximum entry in the solution vector remained $O(\max(U_0))$ and the algorithm was never stopped. Thus there was no need to put a cap on the amplitude of the solution. However, this was not the case for the semi-linear approach. For Test1 the solution grows with time to reach $||u||_{\infty} = 0.3$ at $t = 260.4$ when the algorithm is stopped. For 
Test2 the solution grows with time at a faster rate to reach $||u||_{\infty} = 0.3$ at $t = 9.6$. For 
 Test3  the solution remains unchanged up till $t=26$, and after it grows with time to reach $||u||_{\infty} = 0.3$ at $t = 42$.
 
 Thus, the Newton-type methods are numerically more stable and robust than the semi-linear approach. Moreover, the most competitive one is Modified Newton's method as it is the fastest with a similar runtime to that of  the semi-linear approach.
\section{Concluding Remarks}\label{sec:Conclude}
In this paper, we implement Newton-type methods for solving system \eqref{HMC-Disc-Comp}, specifically Newton, Chord and Modified Newton methods. Moreover, we justify the use of these methods by proving several results, in particular the convergence of the implemented methods. \vspace{2mm}

\noindent Although the sufficient conditions for proving our theorems restrict the time interval $\tau$ to be of the order of $h^{2.5}$ or $h^2$, yet in our computational implementations this restriction was lifted as we were able to use $\tau = O(h)$ without any difficulty.  Proving the mathematical validity of such choices remains an open question.   \vspace{2mm}

\noindent In terms of implementation, given a relative tolerance, all the methods converged in at most $k=2$ iterations per time-step, for all the tested cases. Moreover, the expected runtime behavior is observed, where Modified Newton's method is $\frac{T}{\tau}10^2$-times faster than Chord's method which is $k$-times faster than Newton's method.  \vspace{2mm}

\noindent On the other hand, the Newton-type methods are numerically more stable and robust than the semi-linear approach introduced in \cite{FEHM}, since there was no need to put a cap on the amplitude of the solution in the algorithm. Yet, the time evolution of the solution using the Newton-type methods followed the expected behavior for the corresponding cases. In addition, Modified Newton's method has a similar runtime to that of  the semi-linear approach. \vspace{2mm}
 
 \noindent Thus,  Modified Newton's method appears to be the most competent and robust version to be used for simulations of the Hasegawa-Mima plasma model.\vspace{2mm}

\noindent As for future avenues of research, these include principally the following.
\begin{enumerate}
\item
Proof of convergence of the solution to the nonlinear \eqref{HMC-Comp}, and \eqref{HMC-Disc-Comp} schemes as $\tau$ and $h$ go to zero, which is currently being investigated.
\item
 Another interesting problem for which these methods can be applied is 
 the {Modon Traveling Waves Solutions} to \eqref{HMC2}. These solutions are obtained by  considering the pair of variables $(\xi,\eta)$ given by $\xi=x\,\,\eta=y-ct$, one  looks for solutions to (\ref{HMC2}) in the form $u(x,y,t)=\phi(\xi,\eta)=\phi(x,y-ct)\mbox{ and }w(x,y,t)= \psi(\xi,\eta)=\psi(x,y-ct)$. By defining $\forall t\in(0,T):\,\Omega_t=\{\xi,\eta\,| 0<\xi<L,\,-ct<\eta<L-ct\}, $ then in terms of $\phi$ and $\psi$, the system (\ref{HMC2}) reduces to be solved on $\Omega_0=\Omega$. Thus, with $\nabla=\nabla_{\xi,\eta}$, one seeks $\{\phi,\psi\}:\overline{\Omega}\to\mathbb{R}^2$, such that: 
\begin{equation}\label{HMCTravel}
\left\{\begin{array}{lll}
-c\psi_\eta + \vec{V}(\phi) \cdot \nabla\psi = k\phi_\eta & \mbox{on }  \Omega  &\\
-\Delta \phi+\phi=\psi &\mbox{on }  \Omega &\\
\mbox{PBC's on } \phi,\, \phi_\xi,\, \phi_\eta,\,\psi & \mbox{on }  \partial \Omega & 
\end{array}\right.
\end{equation}
Undergoing research is being carried out on this problem. 
\end{enumerate}

\bibliographystyle{ieeetr}
\bibliography{references}

\begin{thebibliography}{1}

\bibitem{kn}
H.~Karakazian and N.~Nassif, ``Local existence of an $h^3_p$ solution to the
  hasegawa-mima plasma equation,'' {\em Submitted. arXiv:1712.05524}, 2019.

\bibitem{FEHM}
H.~Karakazian, S.~Moufawad, and N.~Nassif, ``A finite-element model for the
  hasegawa–mima wave equation,'' {\em Applied Mathematics and Computation},
  vol.~412, p.~126550, 2022.

\bibitem{hm77}
A.~Hasegawa and K.~Mima, ``Stationary spectrum of strong turbulence in
  magnetized nonuniform plasma,'' {\em Physics of Fluids}, vol.~39,
  pp.~205--208, Jul 1977.

\bibitem{hm78}
A.~Hasegawa and K.~Mima, ``Pseudo-three-dimensional turbulence in magnetized
  nonuniform plasma,'' {\em Physics of Fluids}, vol.~21, pp.~87--92, Jan 1978.

\bibitem{ciarlet}
P.~G. Ciarlet, {\em The Finite Element Method for Elliptic Problems}.
\newblock SIAM, 1979.

\bibitem{FEHMX}
H.~Karakazian, S.~Moufawad, and N.~Nassif, ``A finite-element model for the
  hasegawa-mima wave equation,'' 2021.

\bibitem{MR3043640}
F.~Hecht, ``New development in freefem++,'' {\em J. Numer. Math.}, vol.~20,
  no.~3-4, pp.~251--265, 2012.

\end{thebibliography}

\end{document}